\documentclass[english]{article}
\usepackage{amssymb}

\usepackage[T1]{fontenc}
\usepackage[latin1]{inputenc}
\usepackage{babel}

\makeatletter
\def \nz{\normalsize}

 \makeatother
\input{tcilatex}

\begin{document}

\title{Small time Edgeworth-type expansions for weakly convergent nonhomogeneous  Markov chains.\thanks{
 This research was supported by grant 436RUS113/467/81-2
from the Deutsche Forschungsgemeinschaft and by grants 05-01-04004
and 04-01-00700 from the Russian Foundation of Fundamental
Researches. The first author worked on the paper during a visit at
the Laboratory of Probability Theory and Random Models of the
University Paris VI in 2006. He is grateful for the hospitality
during his stay. This is the long version of a paper that has been
submitted to Probability Theory and Related Fields for
publication. We would like to thank Stephane Menozzi, two referees
and the associate editor for helpful comments.
 }\\}
\author{Valentin KONAKOV \\
  \nz Central Economics Mathematical Institute, Academy of Sciences\\
  \nz Nahimovskii av. 47, 117418 Moscow, Russia\\ \nz E mail:
  kv24@mail.ru\\
  Enno MAMMEN \\
  \nz Department of Economics, University of Mannheim \\
  \nz L 7, 3-5, 68131 Mannheim, Germany\\
  \nz E mail: emammen@rumms.uni-mannheim.de}

 \date{\today{}}
\maketitle

\bigskip

\begin{abstract} \noindent
We consider triangular arrays of Markov chains that converge weakly
to a diffusion process. Second order Edgeworth type expansions for
transition densities are proved. The paper differs from recent
results in two respects. We allow nonhomogeneous diffusion limits
and we treat transition densities with time lag converging to zero.
Small time asymptotics are motivated by statistical applications and
by resulting approximations for the joint density of diffusion
values at an increasing grid of points.
\end{abstract}


\noindent \textsl{Short title:} Edgeworth-type expansions for
 Markov chains

\pagebreak

\section{Introduction.}

Recently, there was some activity on Edgeworth-type expansions for
dependent data. In most approaches higher order expansions have been
derived by application of classical Edgeworth expansions for
independent data. The approaches differ in their main idea how the
dependence structure can be reduced to the case of independent data.
For sums of independent random variables and for functionals of such
sums the theory of Edgeworth expansions is classical and well
understood in a very general setting (see Bhattcharya and Rao (1976)
and G\"otze (1989)). For models with dependent variables three
approaches have been developed where the expansion is derived from
models with sums of independent random variables. In the first
approach mixing properties are used to approximate the Markov chain
by a sum of independent random variables and it is shown that their
Edgeworth expansion carries over to the Markov chain up to a certain
accuracy. The mixing approach was first used by G\"otze and Hipp
(1983) and it was further applied to continuous time processes in
Kusuoka and Yoshida (2000) and Yoshida (2004). Under appropriate
conditions Markov chains can be splitted at regeneration times into
a sequence of i.i.d.\ variables. This fact has been used in
Bolthausen (1980, 1982) to get Berry-Esseen bounds for Markov
chains. For the statement of Edgeworth expansions the regenerative
method has been used in Malinovskii (1987), Jensen (1989), Bertail
and Clemencon (2004) and Fukasawa (2006a). The higher order
Edgeworth expansions have been used to show higher order accuracy of
different bootstrap schemes, see Mykland (1992),Bertail and
Clemencon (2006) and Fukasawa (2006b).

Both approaches, the mixing method and the regenerative method only
have been used for Markov chains with a Gaussian limit. In this
paper we study Markov chains that converge weakly to a diffusion
limit. For the treatment of this case we make use of the parametrix
method. In this approach the transition density is represented as a
nested sum of functionals of densities of sums of independent
variables. Plugging Edgeworth expansions into this representation
will result in an expansion for the transition density. Thus as in
the mixing method and in the regenerative method the expansion is
reduced to models with sums of independent random variables.

The parametrix method permits to obtain tractable representations of
transition densities of diffusions and of Markov chains. For
diffusions the parametrix expansion is based on Gaussian densities,
see Lemma 1 below, and  standard references for the parametrix
method are the books of Friedman (1964) and Ladyzenskaja, Solonnikov
and Ural'ceva (1968) on parabolic PDE [see also McKean and Singer
(1967)]. For a short exposition of the parametrix method, see
Section \ref{sec:Para} and Konakov and Mammen (2000). Similar
representations hold for discrete time Markov
chains $%
X_{k,h}$, see Lemma 3 below.  The parametrix method for Markov
chains was developed in Konakov and Mammen (2000) and it is exposed
in Section \ref{sec:Parchain}. In Konakov and Mammen (2002) the
approach was used to state Edgeworth-type expansions for Euler
schemes for stochastic differential equations. Related treatments of
Euler schemes can be found in Bally and Talay (1996 a,b), Protter
and Talay (1997), Jacod and Protter (1998), Jacod (2004), Jacod,
Kurtz, Meleard and Protter (2005) and Guyon (2006).

In this paper we study triangular arrays of Markov chains $X_{k,h}$
$(k\geq 0)$ that converge weakly to a diffusion process $Y_{s}$
$(s\geq 0)$ for $n\longrightarrow \infty$. We consider the Markov
chains for the time interval $(0\leq k\leq n)$. The corresponding
time interval of the diffusion is $(0\leq s\leq T)$. The term
$h=T/n$ denotes the discretization step. We allow that $T$ depends
on $n$. In particular, we consider the case that $T\rightarrow 0$
for $n\rightarrow \infty$. Furthermore, we allow nonhomogeneous
diffusion limits.

Weak convergence of the distribution of scaled discrete time Markov
processes to diffusions has been extensively studied in the
literature ( see Skorohod (1965) and Stroock and Varadhan (1979)).
Local limit theorems for Markov chains were given in Konakov and
Molchanov (1984) and Konakov and Mammen (2000, 2002). In Konakov and
Mammen (2000) it was shown that the transition density of a Markov
chain converges with rate $O(n^{-1/2})$ to the transition density in
the diffusion model. For the proof there an analytical approach was
chosen that made essential use of the parametrix method.

The main result of this paper will give Edgeworth type expansions
for the transition densities of the Markov chains $X_{k,h}$ $(0\leq
k\leq n)$. The first order term of the expansion is the transition
density of the  diffusion process $Y_{s}$ $(0\leq s\leq T)$. The
order of the expansion is $o(h^{-1-\delta })$ with $\delta
>0$. Related results were shown in Konakov and Mammen (2005). The work of this paper
generalizes the results in Konakov and Mammen (2005) in two
directions. The time horizon $T$ is allowed to converge to $0$ and
also cases are treated with nonhomogeneous diffusion limit. Small
time asymptotics is done for two reasons.  First of all it allows
approximations for the joint density of values of the Markov chain
at an increasing grid of points. Secondly, it is motivated by
statistical applications. In statistics, diffusion models are used
as an approximation to the truth. They can be motivated by a high
frequency Markov chain that is assumed to run in the background on
a very fine time grid and is only observed on a coarser grid. If
the number of time steps between two observed values of the
process converges to infinity this allows diffusion approximations
(under appropriate conditions). This asymptotics reflects a set up
occurring in the high frequency statistical analysis for financial
data where diffusion approximations are used only for coarser time
scales. For the finest scale discrete pattern in the price
processes become transparent and do not allow diffusion
approximations. The statistical implications of our result will be
discussed elsewhere. The mathematical treatment of nonhomogeneous
diffusion limits with time horizon $T$ going to zero contributes
some additional qualitatively new problems. In this case some
additional terms appear that explode for $T \to 0$ and for this
reason these terms need a qualitatively different treatment as in
the case with fixed $T$. The nonhomogeneity adds an additional
term in the Edgeworth expansion. See also below for more details.

The paper is organized as follows. In the next section we will
present our model for the Markov chain and state our main result
that gives an Edgeworth-type expansion for Markov chains.
Connections with previously known results are also discussed in
Section \ref{sec:Int}. In Section \ref{sec:Pardiff} we will give a
short introduction into the parametrix method for diffusions. In
Section \ref{sec:Parchain} we will recall the parametrix approach
developed in Konakov and Mammen (2000) for Markov chains. Technical
discussions, auxiliary results and proofs are given in Sections
\ref{sec:tools} and \ref{sec:mainres}.

\section{\protect\bigskip {}\label{sec:Int}The main result: an Edgeworth-type
expansion for Markov chains converging to diffusions.}

We consider a family of Markov processes in $\Bbb{R}^{d}$ that have
the following form
\begin{equation}
X_{k+1,h}=X_{k,h}+m\left( kh,X_{k,h}\right) h+\sqrt{h}\xi
_{k+1,h},\;X_{0,h}=x\in \Bbb{R}^{d},\;k=0,...,n-1.  \label{eq:001}
\end{equation}
The innovation sequence $\left( \xi _{i,h}\right) _{i=1,...,n}$ is
assumed to satisfy the Markov assumption: the conditional
distribution of $\xi _{k+1,h}$ given the past
$X_{k,h}=x_{k},...,X_{0,h}=x_{0}$ depends only on the last value
$X_{k,h}=x_{k}$ and has a conditional density $q\left(
kh,x_{k},\cdot \right)$. The conditional covariance matrix
corresponding to this density is denoted by $\sigma (kh,x_{k})$ and the conditional $%
\nu -th$ cumulant by $\chi _{\nu }(kh,x_{k})$. The transition
densities of $\left( X_{i,h}\right) _{i=1,...,n}$ are denoted by
$p_{h}\left( 0,kh,x,\cdot \right)$. The time horizon $T=T(n)\leq 1$
is allowed to depend on $n$ and $h=T/n$ is the discretization step.

We make the following assumptions.

\begin{description}
\item[(A1)] It holds that $\int_{\Bbb{R}^{d}}yq\left(
t,x,y\right) dy=0$ for $0\leq t\leq 1,\ x\in \Bbb{R}^{d}$.

\item [(A2)] There exist positive constants $\sigma _{\star }$ and
$\sigma ^{\star }$ such that the covariance matrix $\sigma \left(
t,x\right) =\int_{\Bbb{R}^{d}}yy^{T}q\left( t,x,y\right) dy$
satisfies
\[
\sigma _{\star }\leq \theta ^{T}\sigma \left( t,x\right) \theta
\leq \sigma ^{\star } \] for all $\left\| \theta \right\| =1$ and
$t\in [0,1]$ and $x\in \Bbb{R}^{d}.$

\item [(A3)] There exist a positive integer
$S^{\prime }$ and a real nonnegative function $\psi \left( y\right)
,$ $y\in \Bbb{R}^{d}$ satisfying $\sup_{y\in \Bbb{R}^{d}}\psi \left(
y\right) <\infty $ and
$%
\int_{\Bbb{R}^{d}}\left\| y\right\| ^{S}\psi \left( y\right)
dy<\infty $ with $S=(S^{\prime }+2)d+4$ such that
\[
\left| D_{y}^{\nu }q\left( t,x,y\right) \right| \leq \psi \left(
y\right) ,\;t\in \lbrack 0,1],\;x,y\in \Bbb{R}^{d}\;\left| \nu
\right| =0,1,2,3,4
\]
and
\[
\left| D_{x}^{\nu }q\left( t,x,y\right) \right| \leq \psi \left(
y\right) ,\;t\in \lbrack 0,1],\;x,y\in \Bbb{R}^{d}\;\left| \nu
\right| =0,1,2.
\]
Moreover, for all $x,y\in R^{d},$ $h>0,0\leq t,t+jh\leq 1,j\geq
j_{0} $, with a bound $j_{0}$ that does not depend on $x,t$,
\[
\left| D_{x}^{\nu }q^{(j)}\left( t,x,y\right) \right| \leq
Cj^{-d/2}\psi \left( j^{-1/2}y\right) ,\left| \nu \right| =0,1,2,3
\]
for a constant $C<\infty$. Here $q^{(j)}(t,x,y)$ denotes the
$j$-fold convolution of $q$ for fixed $x$ as a function of $y$:
\[
q^{(j)}(t,x,y)=\int q^{(j-1)}(t,x,u)q(t+(j-1)h,x,y-u)du,
\]
$q^{(1)}(t,x,y)=q(t,x,y)$. \end{description} Note that the last
condition is motivated by
(A2) and the classical local limit theorem. Note also that for $1\leq j\leq j_{0}$%
\[
\int \left\| y\right\| ^{S}q^{(j)}(t,x,y)dy\leq C(j_{0\;},S).
\]

\begin{description}
\item [(B1)] The functions $m\left( t,x\right) $ and
$\sigma \left( t,x\right) $ and their first and second derivatives
w.r.t.\ $t$ and their derivatives up to the order six w.r.t. $x$ are
continuous and bounded uniformly in $t$ and $x.$ All these functions
are Lipschitz continuous with respect to $x$ with a Lipschitz
constant that does not depend on $t.$ The functions $\chi _{\nu
}(t,x)$, $\left| \nu \right| =3,4,$ are Lipschitz continuous with
respect to $t$ with a Lipschitz constant that does not depend on
$x.$ A sufficient condition for this is the following inequality
\[
\int_{\Bbb{R}^{d}}(1+\left\| z\right\| ^{4})\left| q\left(
t,x,z\right) -q\left( t^{\prime },x,z\right) \right| dz\leq
C\left| t-t^{\prime }\right| ,0\leq t,t^{\prime }\leq 1,
\]
with a constant that does not depend on $x\in \Bbb{R}^{d}.$
Furthemore, $%
D_{x}^{\nu }\sigma \left( t,x\right) $ exist for $\left| \nu
\right| \leq 6$ and are Holder continuous w.r.t. $x$ with a
positive exponent and a constant that does not depend on $t.$

\item [(B2)] There exists $\varkappa
<\frac{1}{5}$ such that $\liminf_{n\rightarrow \infty
}T(n)n^{\varkappa }>0$.
\end{description}

The Markov chain $X_{k,h}$, see (\ref{eq:001}), is an approximation
to the
following stochastic differential equation in $\Bbb{R}^{d}:$%
\[
dY_{s}=m\left( s,Y_{s}\right) ds+\Lambda \left( s,Y_{s}\right)
dW_{s},\;Y_{0}=x\in \Bbb{R}^{d},\;s\in \lbrack 0,T],
\]
where $\left( W_{s}\right) _{s\geq 0}$ is the standard Wiener
process
and $%
\Lambda $ is a symmetric positive definite $d\times d$ matrix such
that
$%
\Lambda \left( s,y\right) \Lambda \left( s,y\right) ^{T}=\sigma
\left( s,y\right) .$ The conditional density of $Y_{t},$ given
$Y_{0}=x$ is denoted by $p\left( 0,t,x,\cdot \right) $. We will use
the following differential operators $L$ and $\widetilde{L}$ :
\[
Lf(s,t,x,y)=\frac{1}{2}\sum_{i,j=1}^{d}\sigma
_{ij}(s,x)\frac{\partial ^{2}f(s,t,x,y)}{\partial x_{i}\partial
x_{j}}+\sum_{i=1}^{d}m_{i}(s,x)\frac{%
\partial f(s,t,x,y)}{\partial x_{i}},
\]
\begin{equation}
\tilde{L}f(s,t,x,y)=\frac{1}{2}\sum_{i,j=1}^{d}\sigma _{ij}(s,y)\frac{%
\partial ^{2}f(s,t,x,y)}{\partial x_{i}\partial x_{j}}%
+\sum_{i=1}^{d}m_{i}(s,y)\frac{\partial f(s,t,x,y)}{\partial
x_{i}}. \label{eq:001a}
\end{equation}
To formulate our main result we need also the following operators
\[
L^{\prime }f(s,t,x,y)=\frac{1}{2}\sum_{i,j=1}^{d}\frac{\partial
\sigma _{ij}(s,x)}{\partial s}\frac{\partial ^{2}f\left(
s,t,x,y\right) }{\partial x_{i}\partial
x_{j}}+\sum_{i=1}^{d}\frac{\partial m_{i}(s,x)}{\partial
s}%
\frac{\partial f\left( s,t,x,y\right) }{\partial x_{i}}
\]
\begin{equation}
\widetilde{L}^{\prime
}f(s,t,v,z)=\frac{1}{2}\sum_{i,j=1}^{d}\frac{\partial \sigma
_{ij}(s,y)}{\partial s}\frac{\partial ^{2}f\left( s,t,x,y\right)
}{%
\partial x_{i}\partial x_{j}}+\sum_{i=1}^{d}\frac{\partial
m_{i}(s,y)}{%
\partial s}\frac{\partial f\left( s,t,x,y\right) }{\partial x_{i}}.
\label{eq:001b}
\end{equation}
and the convolution type binary operation $\otimes :$%
\[
f\otimes g\left( s,t,x,y\right) =\int_{s}^{t}du\int_{R^{d}}f\left(
s,u,x,z\right) g\left( u,t,z,y\right) dz.
\]

Konakov and Mammen (2000) obtained a nonuniform rate of convergence
for the difference $p_{h}\left( 0,T,x,\cdot \right) -p\left(
0,T,x,\cdot \right) $ as $n\rightarrow \infty $ in the case $T\asymp
1.$  Edgeworth type expansions for the case $T\asymp 1$ and
homogenous diffusions were obtained in Konakov and Mammen (2005).
The goal of the present paper is to obtain an Edgeworth type
expansion for nonhomegenous case which remains valid for the both
cases $T\asymp 1$ or $T=o\left( 1\right) $. The following theorem
contains our main result. It gives Edgeworth type expansions for
$p_{h}$. For the statement of the theorem we introduce the following
differential operators
\[
\mathcal{F}_{1}[f](s,t,x,y)=\sum_{\left| \nu \right| =3}\frac{\chi
_{\nu }(s,x)}{\nu !}D_{x}^{\nu }f(s,t,x,y),
\]
\[
\mathcal{F}_{2}[f](s,t,x,y)=\sum_{\left| \nu \right| =4}\frac{\chi
_{\nu }(s,y)}{\nu !}D_{x}^{\nu }f(s,t,x,y).
\]
Furthermore, we introduce two terms corresponding to the classical
Edgeworth expansion (see Bhattacharya and Rao (1976))
\begin{eqnarray}
\widetilde{\pi }_{1}(s,t,x,y)&=&(t-s)\sum_{\left| \nu \right| =3}\frac{%
\overline{\chi }_{\nu }(s,t,y)}{\nu !}D_{x}^{\nu
}\widetilde{p}(s,t,x,y), \label{eq:pi1}
\\
\widetilde{\pi }_{2}(s,t,x,y)&=&(t-s)\sum_{\left| \nu \right| =4}\frac{%
\overline{\chi }_{\nu }(s,t,y)}{\nu !}D_{x}^{\nu
}\widetilde{p}(s,t,x,y) \nonumber \\
&& +\frac{1}{2}(t-s)^{2}\left\{ \sum_{\left| \nu \right|
=3}\frac{\overline{%
\chi }_{\nu }(s,t,y)}{\nu !}D_{x}^{\nu }\right\}
^{2}\widetilde{p}(s,t,x,y), \label{eq:pi2}
\end{eqnarray}
where
\[
\overline{\chi }_{\nu }(s,t,y)=\frac{1}{t-s}\int_{s}^{t}\chi _{\nu
}(u,y)du
\]
and $\chi _{\nu }(t,x)$ is the $\nu -th$ cumulant of the density of
the innovations $q(t,x,\cdot ).$ The gaussian transition densities
$\widetilde{%
p}(s,t,x,y)$ are defined in (\ref{eq:002}). Note, that in the
homogenous case $\chi _{\nu }(u,y)\equiv \chi _{\nu }(y)$ and
$\overline{\chi }_{\nu }(s,t,y)\equiv \chi _{\nu }(y),$ where $\chi
_{\nu }(y)$ is the $\nu -th$ cumulant of the density $q(y,\cdot ).$
\bigskip

\noindent \textbf{Theorem 1. } \textit{Assume (A1)-(A3), (B1)-(B2).
Then there exists a
constant }%
$\delta >0$\textit{ such that the following expansion holds:}
\[
\sup_{x,y\in R^{d}} T^{d/2}\left( 1+\left\| \frac{y-x}{\sqrt{T}}%
\right\| ^{S^{\prime }}\right) \times \Big|
p_{h}(0,T,x,y)-p(0,T,x,y)
\]
\[
 -h^{1/2}\pi _{1}(0,T,x,y)-h\pi _{2}(0,T,x,y)\Big|
 =O(h^{1+\delta }),
\]
\textit{where S}$^\prime $\textit{is defined in Assumption (A3) and
where}
\begin{eqnarray*}
\pi _{1}(0,T,x,y)&=&(p\otimes \mathcal{F}_{1}[p])(0,T,x,y),
\\
\pi _{2}(0,T,x,y)&=&(p\otimes
\mathcal{F}_{2}[p])(0,T,x,y)+p\otimes \mathcal{F} _{1}[p\otimes
\mathcal{F}_{1}[p]](0,T,x,y) \\ && +\frac{1}{2}p\otimes (L_{\star
}^{2}-L^{2})p(0,T,x,y)-\frac{1}{2}p\otimes (L^{\prime
}-\widetilde{L}^{\prime })p(0,T,x,y). \end{eqnarray*} \textit{Here
}$p(s,t,x,y)$\textit{ is the transition density of the limiting
diffusion }$Y_{s}$\textit{and the operator L}$_{\star }$\textit{
is defined as }$\widetilde{L}$\textit{, but with the coefficients
``frozen'' at the point }$x$\textit{. The norm }$\left\| \cdot
\right\| $\textit{ is the usual Euclidean norm.}
\bigskip

\noindent \textbf{Remark 1. } The terms of the Edgeworth expansion
have subgaussian tails and are of order $n^{-1/2}$ or $n^{-1}$,
respectively:
\begin{eqnarray*}
\left| h^{1/2}\pi _{1}(0,T,x,y)\right| &\leq&
C_{1}n^{-1/2}T^{-d/2}\exp \left[ -C_{2}\left\|
\frac{y-x}{\sqrt{T}}\right\| ^{2}\right] ,
\\
\left| h\pi _{2}(0,T,x,y)\right| &\leq& C_{1}n^{-1}T^{-d/2}\exp
\left[ -C_{2}\left\| \frac{y-x}{\sqrt{T}}\right\| ^{2}\right] ,
\end{eqnarray*}
with some positive constants $C_{1}$ and $C_{2}.$
\bigskip

\noindent \textbf{Remark 2. } If the innovation density
$q(t,x,\cdot )$ and the conditional mean $m(t,x)$ do not depend on
$x$ then we are in the classical case of independent non
identically distributed random vectors. We now show that then the
Edgeworth expansion of Theorem 1 coincides with the first two
terms of the classical Edgeworth
expansion $h^{1/2}\widetilde{\pi }_{1}(0,T,x,y)+h\widetilde{\pi }%
_{2}(0,T,x,y)$.  Note first that in this case $L_{\star
}=L,L^{\prime }=\widetilde{L}^{\prime }$ and
$p(s,t,x,y)=%
\widetilde{p}(s,t,x,y)$ where $\widetilde{p}$ is defined in
(\ref{eq:002}) with $\sigma (s,t,y)=\sigma (s,t)=\int_{s}^{t}\sigma
(u)du$ and \
$%
m(s,t,y)=m(s,t)=\int_{s}^{t}m(u)du.$ This gives \begin{eqnarray*}
&&\pi _{1}(0,T,x,y)=\int_{0}^{T}ds\int
\widetilde{p}(0,s,x,v)\sum_{\left| \nu \right| =3}\frac{\chi _{\nu
}(s)}{\nu !}D_{v}^{\nu }\widetilde{p}(s,T,v,y)dv \\ && \qquad
=-\sum_{\left| \nu \right| =3}\int_{0}^{T}\frac{\chi _{\nu }(s)}{\nu
!}%
dsD_{y}^{\nu }\int \widetilde{p}(0,s,x,v)\widetilde{p}(s,T,v,y)dv
\\
&&\qquad =-\sum_{\left| \nu \right| =3}\frac{T}{\nu !}\overline{\chi
}_{\nu }(0,T)D_{y}^{\nu }\widetilde{p}(0,T,x,y) \\ && \qquad
=\sum_{\left| \nu \right| =3}\frac{T}{\nu !}\overline{\chi }_{\nu
}(0,T)D_{x}^{\nu }\widetilde{p}(0,T,x,y)=\widetilde{\pi
}_{1}(0,T,x,y), \\ && \widetilde{p}\otimes
\mathcal{F}_{1}[\widetilde{p}](s,T,z,y)=\int_{s}^{T}du%
\int \widetilde{p}(s,u,z,w)\sum_{\left| \nu \right| =3}\frac{\chi
_{\nu
}(u)%
}{\nu !}D_{w}^{\nu }\widetilde{p}(u,T,w,y)dw
\\
&&\qquad
=-\sum_{\left| \nu \right| =3}\int_{s}^{T}\frac{\chi _{\nu
}(u)}{\nu
!}du%
\text{ }D_{y}^{\nu }\widetilde{p}(s,T,z,y)=(T-s)\sum_{\left| \nu
\right| =3}%
\frac{\overline{\chi }_{\nu }(s,T)}{\nu !}D_{z}^{\nu }\widetilde{p}%
(s,T,z,y), \\ &&
\mathcal{F}_{1}[\widetilde{p}\otimes \mathcal{F}_{1}[\widetilde{p}%
]](s,T,z,y)=(T-s)\sum_{\left| \nu \right| =3}\frac{\chi _{\nu
}(s)}{\nu
!}%
D_{z}^{\nu }\left[ \sum_{\left| \nu ^{\prime }\right|
=3}\frac{\overline{%
\chi }_{\nu ^{\prime }}(s,T)}{\nu ^{\prime }!}D_{z}^{\nu ^{\prime }}%
\widetilde{p}(s,T,z,y)\right] \\ && \qquad =(T-s)\sum_{\left| \nu
\right| =3,\left| \nu ^{\prime }\right| =3}\frac{\chi _{\nu
}(s)}{\nu !}\frac{\overline{\chi }_{\nu ^{\prime }}(s,T)}{\nu
^{\prime }!}D_{z}^{\nu +\nu ^{\prime }}\widetilde{p}(s,T,z,y),
\\
&& \widetilde{p}\otimes
\mathcal{F}_{2}[\widetilde{p}](0,T,x,y)+\widetilde{p}%
\otimes \mathcal{F}_{1}[\widetilde{p}\otimes
\mathcal{F}_{1}[\widetilde{p}%
]](0,T,x,y)=T\sum_{\left| \nu \right| =4}\frac{\overline{\chi
}_{\nu
}(0,T)}{%
\nu !}D_{x}^{\nu }\widetilde{p}(0,T,x,y) \\ && \qquad \qquad
+\int_{0}^{T}ds\int \widetilde{p}(0,s,x,z)(T-s)\sum_{\left| \nu
\right|
=3,\left| \nu ^{\prime }\right| =3}\frac{\chi _{\nu }(s)}{\nu !}\frac{%
\overline{\chi }_{\nu ^{\prime }}(s,T)}{\nu ^{\prime }!}D_{z}^{\nu
+\nu ^{\prime }}\widetilde{p}(s,T,z,y)dz \\ && \qquad =T\sum_{\left|
\nu \right| =4}\frac{\overline{\chi }_{\nu }(0,T)}{\nu
!}%
D_{x}^{\nu }\widetilde{p}(0,T,x,y)+\sum_{\left| \nu \right|
=3,\left| \nu
^{\prime }\right| =3}\frac{1}{\nu !}\frac{1}{\nu ^{\prime }!}%
\int_{0}^{T}\chi _{\nu }(s)\left( \int_{s}^{T}\chi _{\nu ^{\prime
}}(u)du\right) dsD_{x}^{\nu +\nu ^{\prime }}\widetilde{p}(s,T,x,y).
\end{eqnarray*}
For $\nu =\nu ^{\prime }$ we have
\begin{eqnarray*}
&&\int_{0}^{T}\chi _{\nu }(s)\left( \int_{s}^{T}\chi _{\nu ^{\prime
}}(u)du\right) ds=\frac{1}{2}\int_{0}^{T}\int_{0}^{T}\chi _{\nu
}(s)\chi _{\nu }(u)dsdu=\frac{T^{2}}{2}\overline{\chi }_{\nu
}(0,T)\overline{\chi }%
_{\nu }(0,T).
\end{eqnarray*}
For $\nu \neq \nu ^{\prime }$ we get
\begin{eqnarray*}
&&\int_{0}^{T}\chi _{\nu }(s)\left( \int_{s}^{T}\chi _{\nu ^{\prime
}}(u)du\right) ds+\int_{0}^{T}\chi _{\nu ^{\prime }}(s)\left(
\int_{s}^{T}\chi _{\nu }(u)du\right) ds
\\ && \qquad
=\int_{0}^{T}\int_{s}^{T}\left[ \chi _{\nu }(s)\chi _{\nu ^{\prime
}}(u)+\chi _{\nu ^{\prime }}(s)\chi _{\nu }(u)\right] dsdu
\\ && \qquad =\frac{1}{2}%
\int_{0}^{T}\int_{0}^{T}\left[ \chi _{\nu }(s)\chi _{\nu ^{\prime
}}(u)+\chi _{\nu ^{\prime }}(s)\chi _{\nu }(u)\right] dsdu
\\ && \qquad
=\frac{T^{2}}{2}\overline{\chi }_{\nu }(0,T)\overline{\chi }_{\nu
^{\prime }}(0,T)+\frac{T^{2}}{2}\overline{\chi }_{\nu ^{\prime
}}(0,T)\overline{\chi }%
_{\nu }(0,T).
\end{eqnarray*}
From these equations we obtain
\begin{eqnarray*}&&
\widetilde{p}\otimes
\mathcal{F}_{2}[\widetilde{p}](0,T,x,y)+\widetilde{p}%
\otimes \mathcal{F}_{1}[\widetilde{p}\otimes
\mathcal{F}_{1}[\widetilde{p}%
](0,T,x,y)
\\ && \qquad
=T\sum_{\left| \nu \right| =4}\frac{\overline{\chi }_{\nu
}(0,T)}{\nu
!}%
D_{x}^{\nu }\widetilde{p}(0,T,x,y)+\frac{T^{2}}{2}\left\{
\sum_{\left| \nu \right| =3}\frac{\overline{\chi }_{\nu }(0,T)}{\nu
!}D_{x}^{\nu }\right\} ^{2}\widetilde{p}(0,T,x,y)\\ && \qquad
=\widetilde{\pi }_{2}(0,T,x,y).
\end{eqnarray*}
This shows the claim that we get for this case the first two terms
of the classical Edgeworth expansion.
\bigskip

\noindent \textbf{Remark 3. } If $\chi _{\nu }(t,x)=0$ for $\left|
\nu \right| =3$ and for $t\in \lbrack 0,T]\times R^{d}$ then it
holds that $\
\mathcal{F}%
_{1}\equiv 0.$ The Theorem 1 holds with
\begin{eqnarray*}
\pi _{1}(0,T,x,y)&=&0,
\\ \pi _{2}(0,T,x,y)&=&(p\otimes
\mathcal{F}_{2}[p])(0,T,x,y) +\frac{1}{2}p\otimes (L_{\star
}^{2}-L^{2})p(0,T,x,y)-\frac{1}{2}p\otimes (L^{\prime
}-\widetilde{L}^{\prime })p(0,T,x,y).
\end{eqnarray*}
If in addition $\chi _{\nu }(t,x)=0$ for $\left| \nu \right| =4$
then the first four moments of the innovations coincide with the
first four moments of a normal distribution with zero mean and
covariance matrix
$%
\sigma (t,x).$ In this case we have $\mathcal{F}_{2}=0$ and we
have
\begin{eqnarray*}
\pi _{1}(0,T,x,y)&=&0,
\\
\pi _{2}(0,T,x,y)&=&\frac{1}{2}p\otimes (L_{\star
}^{2}-L^{2})p(0,T,x,y)-\frac{%
1}{2}p\otimes (L^{\prime }-\widetilde{L}^{\prime })p(0,T,x,y)
\end{eqnarray*}
and the first two terms of the Edgeworth expansion do not depend on
the innovation density. In particular, it holds that $\chi _{\nu
}(t,x)=0$ for $\left| \nu \right| =3,4$ for Markov chains that are
defined by Euler approximations to diffusions. Thus, an Edgeworth
expansion for the Euler scheme holds with the same $\pi _{1}$ and
$\pi _{2}$ as just defined. For the homogenous case we have that
$L^{\prime }=\widetilde{L}^{\prime }=0$ and we obtain for the Euler
scheme in this case
\begin{eqnarray*}
\pi _{1}(0,T,x,y)&=&0,
\\
\pi _{2}(0,T,x,y)&=&\frac{1}{2}p\otimes (L_{\star
}^{2}-L^{2})p(0,T,x,y).
\end{eqnarray*}
This result for $T=[0,1]$ under Hormander's condition on a diffusion
matrix was obtained by Bally and Talay (1996).

\bigskip

\noindent \textbf{Remark 4. } We now shortly discuss an application
of our result to statistics. Assume that one observes a Markov
process $X_{1,h},...,X_{nk,h}$ at time points $k, 2k,...,nk$. That
means we assume that a high frequency Markov chain runs in the
background on a very fine time grid but that it is only observed on
a coarser grid.  This asymptotics reflects a set up occurring in the
high frequency statistical analysis for financial data where
diffusion approximations are used only for coarser time scales. For
the finest scale discrete pattern in the price processes become
transparent that could not be modeled by diffusions. The joint
distribution of the observed values of the Markov process is denoted
by $P_h$. We assume that this joint distribution can be approximated
by the distribution of $(Y_1,...,Y_n)$ where $Y_1,...,Y_n$ are the
values of a diffusion on the equidistant grid $kh,2kh,...,nkh$. The
joint distribution of $(Y_1,...,Y_n)$ is denoted by $Q_h$. According
to our theorem the one-dimensional marginal distributions of $P_h$
can be approximated by the one-dimensional marginal distributions of
$Q_h$. Under appropriate conditions the L$_1$-norm of this
difference is of order $k^{-1/2}$. This implies that the L$_1$-norm
of the difference between the joint distributions $P_h$ and $Q_h$ is
of order $nk^{-1/2}$. That means the diffusion approximation is only
accurate if $k \gg n^2$, i.e.\ only if the grid of observed points
is very coarse in comparison to the grid on which the Markov process
lives. Only in this case it can be guarantied that a statistical
inference that is based on the diffusion model is accurate. Or put
it in another way, data that come from the Markov model could not be
asymptotically statistically distinguished from diffusion
observations. Our results help to analyze what may go wrong if $k
\gg n^2$ does not hold. The (signed) transition densities $p+h^{1/2}
\pi_1 + h \pi_2$ given in the statement of Theorem 1 define a joined
(signed) measure $R_h$. According to Theorem 1, the marginal
distributions of $R_h$ approximate the one-dimensional marginal
distributions of $P_h$ with order $o(k^{-1-\delta})$. One may
conjecture that under some regularity assumptions the exact order is
$k^{-3/2}$. This implies that $\|P_h-R_h\|_1$ is of order
$nk^{-3/2}$. Thus, this approximation is appropriate as long as $k
\gg n^{2/3}$. This is a much more acceptable assumption. Now,  one
can check which statistical procedures behave differently under the
models $Q_h$ and $R_h$. These procedures may lead to erroneous
conclusions for the Markov data.

\section{\label{sec:Para} The parametrix method.}
\subsection{\label{sec:Pardiff} The parametrix method for diffusions.}

We now give a short overview on the parametrix method for
diffusions. For any $s\in \lbrack 0,T],$ $x,y\in \Bbb{R}^{d}$ we
consider the following
 family of ''frozen'' diffusion processes
\[
d\widetilde{Y}_{t}=m\left( t,y\right) dt+\Lambda \left( t,y\right)
dW_{t},\;%
\widetilde{Y}_{s}=x,\;s\leq t\leq T.
\]
Let $\widetilde{p}^{y}\left( s,t,x,\cdot \right) $ be the
conditional density of $\widetilde{Y}_{t},$ given
$\widetilde{Y}_{s}=x$. In the sequel for any $z$ we will denote
$\widetilde{p}\left( s,t,x,z\right)
=\widetilde{p%
}^{z}\left( s,t,x,z\right) ,$ where the variable $z$ acts here
twice: as the argument of the density and as defining quantity of
the process
$\widetilde{Y}%
_{t}.$

The transition densities $\widetilde{p}$ can be computed
explicitly
\begin{eqnarray}
\widetilde{p}\left( s,t,x,y\right) &=&\left( 2\pi \right)
^{-d/2}\left( \det
\sigma \left( s,t,y\right) \right) ^{-1/2}  \nonumber \\
&&\times \exp \left( -\frac{1}{2}\left( y-x-m\left( s,t,y\right)
\right) ^{T}\sigma ^{-1}\left( s,t,y\right) \left( y-x-m\left(
s,t,y\right) \right) \right) ,  \label{eq:002}
\end{eqnarray}
where
\[
\sigma \left( s,t,y\right) =\int_{s}^{t}\sigma \left( u,y\right)
du,\;\;\;\;\;m\left( s,t,y\right) =\int_{s}^{t}m\left( u,y\right)
du.
\]
Note that the following differential operators $L$ and $\tilde{L}$
correspond to the infinitesimal operators of $Y$ or of the frozen
process $%
\widetilde{Y}$, respectively, i.e.
\[
Lf(s,t,x,y)=\lim_{h\rightarrow 0}h^{-1}\{E[f(s,t,Y(s+h),y)\mid
Y(s)=x]-f(s,t,x,y)\},
\]
\[
\tilde{L}f(s,t,x,y)=\lim_{h\rightarrow 0}h^{-1}\{E[f(s,t,\tilde{Y}%
(s+h),y)\mid \widetilde{Y}(s)=x]-f(s,t,x,y)\}.
\]
We put
\[
H=(L-\tilde{L})\tilde{p}.
\]
Then
\begin{eqnarray*}
H\left( s,t,x,y\right) &=&\frac{1}{2}\sum_{i,j=1}^{d}\left( \sigma
_{ij}\left( s,x\right) -\sigma _{ij}\left( s,y\right) \right)
\frac{\partial ^{2}\widetilde{p}\left( s,t,x,y\right) }{\partial
x_{i}\partial x_{j}}
\\
&&+\sum_{i,j=1}^{d}\left( m_{i}\left( s,x\right) -m_{i}\left(
s,y\right) \right) \frac{\partial \widetilde{p}\left(
s,t,x,y\right) }{\partial x_{i}}.
\end{eqnarray*}
In the following lemmas the $k$- fold convolution of $H$ is denoted
by $%
H^{\left( k\right) }.$ The following results have been proved in
Konakov and Mammen (2000).
\bigskip

\noindent \textbf{Lemma 1. \label{lem1}}\textit{Let }$0\leq s<t\leq
T.$ \textit{It holds}
\emph{%
\label{km-1}
\[
p(s,t,x,y)=\sum_{r=0}^{\infty }\widetilde{p}\otimes
H^{(r)}(s,t,x,y).
\]
}
\bigskip

\noindent \textbf{Lemma 2. \label{km-2}} \textit{Let }$0\leq
s<t\leq T$\textit{. There are constants }$C$\textit{ and
}$C_{1}$\textit{\ such that}\emph{\
\[
\left| H(s,t,x,y)\right| \leq C_{1}\rho ^{-1}\phi _{C,\rho }(y-x)
\]
}

\textit{and}\emph{\
\[
\left| \widetilde{p}\otimes H^{(r)}(s,t,x,y)\right| \leq
C_{1}^{r+1}\frac{%
\rho ^{r}}{\Gamma (1+\frac{r}{2})}\phi _{C,\rho }(y-x),
\]
}

\textit{where }$\rho ^{2}=t-s,$\textit{ }$\phi _{C,\rho }(u)=\rho
^{-d}\phi _{C}(u/\rho )$\textit{and}\emph{\
\[
\phi _{C}(u)=\frac{\exp (-C\left\| u\right\|^{2} )}{\int \exp
(-C\left\| v\right\|^{2} dv)}.
\]
}

\subsection{\label{sec:Parchain}The parametrix method for Markov chains.}

\bigskip We now give a short overview on the parametrix method for
Markov chains. This theory was developed in Konakov and Mammen
(2000). For any $0\leq jh\leq T,$ $x,y\in \Bbb{R}^{d}$ we consider
an additional family of ''frozen'' Markov chains defined for $jh\leq
ih\leq T$ as
\begin{equation}
\widetilde{X}_{i+1,h}=\widetilde{X}_{i,h}+m\left( ih,y\right)
h+\sqrt{h}%
\widetilde{\xi }_{i+1,h},\;\widetilde{X}_{j,h}=x\in
\Bbb{R}^{d},\,\,j\leq i\leq n,  \label{Froz001}
\end{equation}
where $\widetilde{\xi }_{j+1,h},...,\widetilde{\xi }_{n,h}$ is an
innovation sequence such that the conditional density of
$\widetilde{\xi }_{i+1,h}$ given the past
$\widetilde{X}_{i,h}=x_{i},...,\widetilde{X}_{0,h}=x_{0}$ equals
to $q\left( ih,y,\cdot \right) .$ Let us introduce the
infinitesimal operators corresponding to Markov chains
(\ref{eq:001}) and (\ref{Froz001}) respectively,
\[
L_{h}f\left( jh,kh,x,y\right) =h^{-1}\left( \int p_{h}\left(
jh,\left( j+1\right) h,x,z\right) f\left( \left( j+1\right)
h,kh,z,y\right) dz-f\left( \left( j+1\right) h,kh,x,y\right)
\right)
\]
and
\[
\widetilde{L}_{h}f\left( jh,kh,x,y\right) =h^{-1}\left( \int
\widetilde{p}%
_{h}^{y}\left( jh,\left( j+1\right) h,x,z\right) f\left( \left(
j+1\right) h,kh,z,y\right) dz-f\left( \left( j+1\right)
h,kh,x,y\right) \right) ,
\]
where $\widetilde{p}_{h}^{y}\left( jh,j^{\prime }h,x,\cdot \right)
$ denotes the conditional density of $\widetilde{X}_{j^{\prime
},h}$ given
$\widetilde{%
X}_{j,h}=x.$ Similarly as above, for brevity for any $z$ we write $\widetilde{p}%
_{h}\left( jh,j^{\prime }h,x,z\right) =\widetilde{p}_{h}^{z}\left(
jh,j^{\prime }h,x,z\right) ,$ where the variable $z$ acts here
twice: as the argument of the density and as defining quantity of
the process
$\widetilde{X%
}_{i,h}.$ For technical convenience the terms $f\left( \left(
j+1\right) h,kh,z,y\right) $ on the right hand side of $L_{h}f$
and $\widetilde{L}_{h}f$ appear instead of $f\left(
jh,kh,z,y\right) .$

In analogy with the definition of $H$ we put, for $k>j,$%
\[
H_{h}\left( jh,kh,x,y\right) =\left(
L_{h}-\widetilde{L}_{h}\right) \widetilde{p}_{h}\left(
jh,kh,x,y\right) .
\]
We also shall use the convolution type binary operation $\otimes
_{h}$ which is a discrete version of $\otimes $:
\[
g\otimes _{h}f\left( jh,kh,x,y\right) =\sum_{i=j}^{k-1}h\int_{\Bbb{R}%
^{d}}g\left( jh,ih,x,z\right) f\left( ih,kh,z,y\right) dz,
\]
where $0\leq j<k\leq n.$ We write $g\otimes _{h}H_{h}^{\left(
0\right) }=g$ and $g\otimes _{h}H_{h}^{\left( r\right) }=\left(
g\otimes _{h}H_{h}^{\left( r-1\right) }\right) \otimes _{h}H_{h}$
for $r\geq 1.$ For the higher order convolutions we use the
convention $\sum_{i=j}^{l}=0$ for $l<j.$ One can show the following
analog of the ''parametrix '' expansion for $p_{h}$
{[}%
see Konakov and Mammen (2000){]}.
\bigskip

\noindent \textbf{Lemma 3. }\textit{Let }$0\leq jh<kh\leq T.$ \
\textit{It holds} \emph{\label{km-3}
\[
p_{h}(jh,kh,x,y)=\sum_{r=0}^{k-j}\widetilde{p}_{h}\otimes
_{h}H_{h}^{(r)}(jh,kh,x,y),
\]
} \textit{where}\emph{\
\[
\widetilde p_{h}(jh,jh,x,y)=p_{h}(kh,kh,x,y)=\delta (y-x)
\]
}\textit{and }$\delta $\textit{ is the Dirac delta symbol.}

\section{\label{sec:tools} Some technical tools.}
\subsection{\label{sec:bonedg} Plugged in Edgeworth expansions for independent observations.}

In this Section we will develop some tools that are helpful for
the comparison of the expansion of $p$ (see Lemma 1) and the
expansion of $p_{h} $ ( see Lemma 3). These expansions are simple
expressions in $\tilde{p}$ or $\tilde{p}_{h}$, respectively.
Recall that $\tilde{p}$ is a Gaussian density, see (\ref{eq:002}),
and that $\tilde{p}_{h}$ is the density of a sum of independent
variables. The densities $\tilde{p}$ and $\tilde{p}_{h}$ can be
compared by application of the classical Edgeworth expansions.
This is done in Lemma 5 and this is the essential step for the
comparison of the expansions of $p$ and $p_{h}$. Lemmas 4 and 7
contain technical tools that will be used below.  Lemma 7 contains
bounds on derivatives of $\tilde{p}_{h}$ that will be used at
several places in the proof of Theorem 1. Its proof makes use of
 Lemma 6 that is a generalisation of a result in Konakov
and Molchanov (1984) (Lemma 4 on page 68). Lemma 5 is a higher
order extension of the results from Section 3.3 in Konakov and
Mammen (2000).

For the formulation of the lemmas we need some additional
notations. Suppose that $X\in \Bbb{R}^{d}$ is a random vector
having a
density $%
q(\mathbf{x}),\mathbf{x\in }\Bbb{R}^{d},$ $EX=0,$$Cov(X,X)=\Sigma
,$
where $%
\Sigma $ be a positively definite $d\times d$ matrix . Denote
$A=\left\| a_{ij}\right\| =\Sigma ^{-1/2}$ and let $\chi _{\nu }(Z)$
be a cumulant of the order $\nu =(\nu _{1},...,\nu _{d})$ of a
random vector $Z\in \Bbb{R}^{d} $, $\phi (x)$ denotes a function in
$\Bbb{R}^{d}$ such that $D_{x}^{\nu
}\phi (x)$ exist and are continuous for $\left| \nu \right| =4$, and $%
A^{-1}=\left\| a^{ij}\right\| =\Sigma ^{1/2}.$
\bigskip

\noindent \textbf{Lemma 4.} \label{km-4}\textit{The following
relation holds for
 }$s=3$%
\textit{ and for }$s=4$%
\[
\sum_{\left| \nu \right| =s}\frac{\chi _{\nu }(AX)D_{z}^{\nu }\phi
(z)}{\nu !%
}=\sum_{\left| \nu \right| =s}\frac{\chi _{\nu }(X)D_{x}^{\nu
}\phi
(Ax)}{%
\nu !}
\]
\textit{where }$z=Ax.$

Denote
\begin{equation}
\mu
_{j,k}(y)=h\sum_{i=j}^{k-1}m(ih,y),V_{j,k}(y)=h\sum_{i=j}^{k-1}\sigma
(ih,y).  \label{eq:003a}
\end{equation}
\bigskip

\emph{\textsc{Proof of Lemma 4.}} For $\left| \nu \right| =3,\nu
=(\nu _{1},...,\nu _{d}),$ each cumulant $\chi _{\nu }(AX)$ is a
linear combination of $\chi _{\mu }(X)$ with $\left| \mu \right| =3$
and with coefficients depending only on $a_{ij}$. It follows from
the following relation
\[
\chi _{\nu }(AX)=\mu _{\nu }(AX)=\int
(a_{11}x_{1}+...+a_{1d}x_{d})^{\nu
_{1}}\times ...\times (a_{11}x_{1}+...+a_{1d}x_{d})^{\nu _{d}}q(\mathbf{x})d%
\mathbf{x}.
\]
Analogously, from the usual differentiation rule of a composite
function and
from the relation $\phi (z)=\phi (Ax)$ , $x=A^{-1}z$, it follows that $%
D_{z}^{\nu }\phi (z)=D_{z}^{\nu }\phi (Ax)$ is a linear combination of $%
D_{x}^{\nu }\phi (Ax)$ with coefficients depending only on $a^{ij}$.
As a result of such substitutions we obtain that
\begin{eqnarray*}
&&\sum_{\left| \nu \right| =3}\frac{\chi _{\nu }(AX)D_{z}^{\nu }\phi (z)}{\nu !%
}=\frac{1}{3!}\sum_{j=1}^{d}\left[ \sum_{\left| \mu \right|
=3}\frac{3!}{\mu _{1}!...\mu _{d}!}a_{1j}^{\mu _{1}}...a_{dj}^{\mu
_{d}}\chi _{\mu }(X)\right] \\&&\qquad \times \left[ \sum_{\left|
\mu ^{\prime }\right| =3}\frac{3!}{\mu _{1}^{\prime }!...\mu
_{d}^{\prime }!}(a^{j1})^{\mu _{1}^{\prime }}...(a^{jd})^{\mu
_{d}^{\prime }}D_{x}^{\mu ^{\prime }}\phi (AX)\right]
\\&&\qquad
+\frac{1}{2!1!}\sum_{\{i\neq ,j\}}\left[ \sum_{l=1}^{d}\sum_{\left|
\mu \right| =2}\frac{2!}{\mu _{1}!...\mu _{d}!}a_{1j}^{\mu
_{1}}...a_{dj}^{\mu _{d}}a_{il}\chi _{\mu +e_{l}}(X)\right]
\\&&\qquad
\times \left[ \sum_{l^{\prime }=1}^{d}\sum_{\left| \mu ^{\prime }\right| =2}%
\frac{2!}{\mu _{1}^{\prime }!...\mu _{d}^{\prime }!}(a^{j1})^{\mu
_{1}^{\prime }}...(a^{jd})^{\mu _{d}^{\prime }}a^{il^{\prime
}}D_{x}^{\mu ^{\prime }+e_{l^{\prime }}}\phi (AX)\right]
\\&&\qquad
+\frac{1}{3!}\sum_{\{i\neq j\neq k\}}\left[
\sum_{l,q=1}^{d}\sum_{\left| \mu \right| =1}\frac{1}{\mu _{1}!...\mu
_{d}!}a_{1j}^{\mu _{1}}...a_{dj}^{\mu _{d}}a_{il}a_{kq}\chi _{\mu
+e_{l}+e_{q}}(X)\right]
\end{eqnarray*} \begin{eqnarray*}
\times \left[ \sum_{l^{\prime },q^{\prime }=1}^{d}\sum_{\left| \mu
^{\prime }\right| =1}\frac{1}{\mu _{1}^{\prime }!...\mu _{d}^{\prime
}!}(a^{j1})^{\mu _{1}^{\prime }}...(a^{jd})^{\mu _{d}^{\prime
}}a^{il^{\prime }}a^{kq^{\prime }}D_{x}^{\mu ^{\prime }+e_{l^{\prime
}}+e_{q^{\prime }}}\phi (AX)\right]
\end{eqnarray*}
where $\sum_{\{i\neq ,j\}}$ ( $\sum_{\{i\neq j\neq k\}}$) denotes
the sum over all different pairs (triples) of $i,j\in \{1,2,...,d\}$
( of $i,j,k\in \{1,2,...,d\}$) and $e_{i}\in \Bbb{R}^{d}$ denotes
the vector whose $i-$th coordinate is equal to 1 and other
coordinates are zero. Collecting the similar terms in the last
equation we obtain that for $\nu =3e_{k}$, $\nu ^{\prime }=3e_{l}$
the coefficient before $\chi _{\nu }(X)D_{x}^{\nu
^{\prime }}\phi (AX)$ is equal to $\frac{1}{3!}%
(a_{1k}a^{l1}+...+a_{dk}a^{ld})^{3}=\frac{1}{3!}\delta _{kl}$, for
$\nu =e_{q}+2e_{r},\text{$\nu ^{\prime }=e_{l}+2e_{n}$ , }q\neq r,$
the coefficient before $\chi _{\nu }(X)D_{x}^{\nu ^{\prime }}\phi
(AX)$ \ is
equal to $\frac{1}{2!}%
(a_{1q}a^{l1}+...+a_{dq}a^{ld})(a_{1r}a^{n1}+...+a_{dr}a^{nd})^{2}=
\frac{1}{2!}\delta _{ql}\delta _{rn}$, \ in particular , for $l=n$ \
the last expression is equal to zero. For \ \ $\ \nu
=e_{q}+e_{r}+e_{n},\nu ^{\prime }=e_{q^{\prime }}+e_{r^{\prime
}}+e_{n^{\prime }}$ \ $q\neq r,q\neq n,r\neq n,$ the coefficient before $%
\chi _{\nu }(X)D_{x}^{\nu ^{\prime }}\phi (AX)$ \ is equal to $%
(a_{1q}a^{q^{\prime }1}+...+a_{dq}a^{q^{\prime }d})\times
(a_{1r}a^{r^{\prime }1}+...+a_{dr}a^{r^{\prime }d})\times $ $%
(a_{1n}a^{n^{\prime }1}+...+a_{dn}a^{n^{\prime }d})=\delta
_{qq^{\prime
}}\delta _{rr^{\prime }}\delta _{nn^{\prime }}$ $.$ This proves lemma for $%
\left| \nu \right| =3.$ The proof for $\left| \nu \right| =4$ is
quite similar. For this case we use the relation which enabes to
express a cumulant $\chi _{\nu }(AX)$ as $\mu _{\nu }(AX)$ plus a
second order polynomial of the moments $\mu _{\nu ^{\prime
}}(AX)$, $\left| \nu ^{\prime
}\right| =2.$ A necessary correction term for $\mu _{\nu }(X)$ to get a $%
\chi _{\nu }(X)$ comes from the derivation of $D_{z}^{\nu }\phi
(z)$. This completes the proof of the lemma.

\bigskip

\noindent \textbf{Lemma 5.} \emph{\label{km-5} }\textit{The
following bound holds with a constant }$C$\textit{ for }$\nu =(\nu
_{1},...\nu _{p})^{T}$\textit{ with }$0\leq \left| \nu \right| \leq
6$
\begin{eqnarray*}
&&\left| D_{z}^{\nu }\widetilde{p}_{h}(jh,kh,x,y)-D_{z}^{\nu
}\widetilde{p}%
(jh,kh,x,y) -\sqrt{h}D_{z}^{\nu }\widetilde{\pi
}_{1}(jh,kh,x,y)-hD_{z}^{\nu
}%
\widetilde{\pi }_{2}(jh,kh,x,y)\right| \\ && \qquad \leq
Ch^{3/2}\rho ^{-3}\zeta _{\rho }^{S-\left| \nu \right| }(y-x)
\end{eqnarray*}
\emph{f}\textit{or all }$j<k,x$\textit{and }$y$\textit{. Here
}$D_{z}^{\nu
}$\textit{denotes the partial differential operator of order }$\nu $%
\textit{ with respect to }$z=V_{j,k}^{-1/2}(y)(y-x-\mu
_{j,k}(y))$\textit{. The quantity }$\rho $\textit{ denotes again the
term }$\rho
=[h(k-j)]^{1/2}$%
\textit{ and the functions }$\widetilde{\pi
}%
_{1} $\textit{ and }$\widetilde{\pi }_{2}$\textit{ are defined in
(\ref{eq:pi1}) and (\ref{eq:pi2}). We write }$\zeta _{\rho
}^{k}(\cdot )=\rho ^{-d}\zeta ^{k}(\cdot /\rho )$\textit{ where}
\[
\zeta ^{k}(z)=\frac{[1+\left\| z\right\| ^{k}]^{-1}}{\int
[1+\left\| z^{\prime }\right\| ^{k}]^{-1}dz^{\prime }}.
\]
\bigskip

\textsc{Proof of Lemma 5.} We note first that
$\tilde{p}_{h}(jh,kh,x,\cdot )$ is the density of the vector
\[
x+\mu _{j,k}(y)+h^{1/2}\sum_{i=j}^{k-1}\widetilde{\xi }_{i+1,h},
\]
where, as above in the definition of the {}``frozen'' Markov chain $\tilde{Y}%
_{n}$, $\widetilde{\xi }_{i+1,h}$ is a sequence of independent
variables with densities $q(ih,y,\cdot ),$ $\mu
_{j,k}(y)=\sum_{i=j}^{k-1}hm(ih,y).$ Let $f_{h}(\cdot )$ be the
density of the normalized sum
\[
h^{1/2}\left[ V_{j,k}(y)\right] ^{-1/2}\sum_{i=j}^{k-1}\widetilde{\xi }%
_{i+1,h}.
\]
Clearly, we have
\[
\tilde{p}_{h}(jh,kh,x,\cdot )=\det \left[ V_{j,k}(y)\right] ^{-1/2}f_{n}\{%
\left[ V_{j,k}(y)\right] ^{-1/2}[\cdot -x-\mu _{j,k}(y)]\}.
\]
We now argue that an Edgeworth expansion holds for $f_{h}$. This
implies the
following expansion for $\tilde{p}_{h}(jh,kh,x,\cdot )$%
\begin{equation}
\tilde{p}_{h}(jh,kh,x,\cdot )  \label{eq:10}
\end{equation}
\[
=\det \left[ V_{j,k}(y)\right] ^{-1/2}[\sum_{r=0}^{S-3}(k-j)^{-r/2}P_{r}(-%
\phi :\{\bar{\chi}_{\beta ,r}\})\{\left[ V_{j,k}(y)\right]
^{-1/2}[\cdot -x-\mu _{j,k}(y)]\}
\]
\[
+[k-j]^{-(S-2)/2}O([1+\left\| \{[V_{j,k}(y)]^{-1/2}[\cdot -x-\mu
_{j,k}(y)]\}\right\| ^{S}]^{-1})]
\]
with standard notations, see Bhattacharya and Rao (1976), p. 53. In
particular, $P_{r}$ denotes a product of a standard normal density
with a polynomial that has coefficients depending only on cumulants
of order $\leq r+2$. Expansion (\ref{eq:10}) follows from Theorem
19.3 in Bhattacharya and Rao (1976). This can be seen as in the
proof of Lemma 3.7 in Konakov and Mammen (2000).

It follows from (\ref{eq:10}) and Condition (A3) that
\[
\left| \tilde{p}_{h}(jh,kh,x,y)-\tilde{p}(jh,kh,x,y)\right.
\left. -h^{1/2}\widehat{\pi }_{1}(jh,kh,x,y)-h\widehat{\pi }%
_{2}(jh,kh,x,y)\right|
\]
\begin{equation}
\leq Ch^{3/2}\rho ^{-3}\zeta _{\rho }^{S-|\nu |}(y-x),
\label{eq:11}
\end{equation}
where
\begin{eqnarray*}
&&\tilde{p}(jh,kh,x,y)=\det \left[ V_{j,k}(y)\right] ^{-1/2}(2\pi
)^{-p/2}
\\&& \qquad
\exp \{-\frac{1}{2}(y-x-\mu _{j,k}(y))^{T}\left[ V_{j,k}(y)\right]
^{-1}(y-x-\mu _{j,k}(y))\},
\\&&
\widehat{\pi }_{1}(jh,kh,x,y)=-\rho ^{-1}\det \left[
V_{j,k}(y)\right] ^{-1/2}\sum_{\left| \nu \right|
=3}\frac{\overline{\chi }_{\nu ,j,k}(y)}{\nu !}D_{z}^{\nu }\phi
\left\{ \left[ V_{j,k}(y)\right] ^{-1/2}(y-x-\mu _{j,k}(y))\right\}
,
\\&&
\widehat{\pi }_{2}(jh,kh,x,y)=\rho ^{-2}\det \left[ V_{j,k}(y)\right] ^{-1/2}%
\left[ \sum_{\left| \nu \right| =4}\frac{\overline{\chi }_{\nu
,j,k}(y)}{\nu !}D_{z}^{\nu }\phi \left\{ \left[ V_{j,k}(y)\right]
^{-1/2}(y-x-\mu _{j,k}(y))\right\} \right.
\\&& \qquad
\left. +\frac{1}{2}\left\{ \sum_{|\nu |=3}\frac{\overline{\chi
}_{\nu
,j,k}(y)}{\nu !}D_{z}^{\nu }\right\} ^{2}\phi \left\{ \left[ V_{j,k}(y)%
\right] ^{-1/2}(y-x-\mu _{j,k}(y))\right\} \right] ,
\end{eqnarray*}
where $\ \overline{\chi }_{\nu
,j,k}(y)=\frac{1}{k-j}\sum_{i=j}^{k-1}\chi _{\nu ,j,k,i}(y),$ $\chi
_{\nu ,j,k,i}(y)=\nu -$th cumulant of $\ \rho \left[
V_{j,k}(y)\right] ^{-1/2}\widetilde{\xi }_{i+1,h}=\rho ^{\left| \nu
\right|
}\times \{\nu -$th cumulant of $\ \left[ V_{j,k}(y)\right] ^{-1/2}\widetilde{%
\xi }_{i+1,h}\},$ and $D_{z}^{\nu }\phi (z)$ denotes the $\nu -$th \
derivative of $\phi $ with respect to $z=\left[ V_{j,k}(y)\right] ^{-1/2}$ $%
(y-x-\mu _{j,k}(y))$\ . It follows from the (conditional) independence of \ $%
\widetilde{\xi }_{i+1,h},i=j,...,k-1,$ that $\ \overline{\chi }_{\nu
,j,k}(y)=\frac{\rho ^{\left| \nu \right| }}{k-j}h^{-\left| \nu
\right| /2}\times \chi _{\nu }(AX),$ where $A=h^{1/2}\left[
V_{j,k}(y)\right]
^{-1/2}=\Sigma ^{-1/2},\Sigma =Cov(X,X),X=\sum_{i=j}^{k-1}\widetilde{\xi }%
_{i+1,h}.$ By Lemma 4 for $s=3,4$
\begin{eqnarray} \nonumber
\sum_{\left| \nu \right| =s}\frac{\overline{\chi }_{\nu ,j,k}(y)}{\nu !}%
D_{z}^{\nu }\phi (z)&=&\rho ^{s}\frac{1}{k-j}\sum_{\left| \nu \right| =s}\frac{%
\chi _{\nu }(AX)}{\nu !}D_{h^{1/2}z}^{\nu }\phi _{h}(h^{1/2}z)
\\ \nonumber &=&(-1)^{s}\rho ^{s}\sum_{\left| \nu \right| =s}\frac{\overline{\chi
}_{\nu }(X)}{\nu !}D_{x}^{\nu }\phi _{h}(A(y-x-\mu _{j,k}(y)))
\\&=&(-1)^{s}\rho ^{s}\sum_{\left| \nu \right| =s}\frac{\overline{\chi
}_{\nu }(X)}{\nu !}D_{x}^{\nu }\phi (\left[ V_{j,k}(y)\right]
^{-1/2}(y-x-\mu _{j,k}(y))),  \label{eq:12}
\end{eqnarray}
where we put $\phi _{h}(z)=\phi (h^{-1/2}z),\overline{\chi }_{\nu }(X)=\frac{%
1}{k-j}\sum_{i=j}^{k-1}\chi _{\nu }(ih,y).$ It follows from
(\ref{eq:12}) and the condition \textbf{B1 }that up to the error
term in the right hand
side of \ (\ref{eq:11}) the \ functions $\widehat{\pi }_{1}$ and $\widehat{%
\pi }_{2}$ coincide with the functions $\widetilde{\pi }_{1}$ and \ $%
\widetilde{\pi }_{2}$ given at the beginning of Section 4. For $\nu
=0$ the statement of the lemma immediately follows from
(\ref{eq:11}). For $\nu >0$ one proceeds similarly. See the remark
at the end of the proof of Lemma 3.7 in Konakov and Mammen (2000).
 \bigskip

\noindent \textbf{Lemma 6. }\textit{Let }$L(d)$\textit{ be the set
of symmetric matrices, and for $0<\lambda ^{-}<\lambda ^{+}<\infty $
let }$D_{\lambda ^{+},\lambda ^{-}}\subset L(d)$\textit{ be the open
subset of $L(d)$ that contains all }$\Lambda \in L(d)$\textit{ with
} $\lambda ^{-}I< \Lambda < \lambda
^{+}I$. %
\textit{ For $\Lambda \in L(d)$ define }$A=A(\Lambda )$\textit{
as the symmetric solution of the equation }$%
A^{2}=\Lambda $\textit{.  Then for any }$k,l,i,j\leq d$\textit{ and
}$\Lambda \in D_{\lambda ^{+},\lambda ^{-}}$\textit{ we have that
with a constant $C_m$ depending on $m$}\emph{}
\begin{equation} \label{eqadd1}
\left| \frac{\partial^m a_{ij}(\Lambda )}{(\partial \lambda
_{kl})^m}\right| \leq C_m (\lambda ^{-})^{-(2m-1)/2}. \end{equation}
\textit{Here $a_{ij}(\Lambda )$ are the elements of $A=A(\Lambda )$.
}
\bigskip

\noindent \textsc{Proof of Lemma 6.}  For $m=1$ the lemma was proved
in Konakov and Molchanov (1984) (see Lemma 4). Suppose now that
(\ref{eqadd1}) holds for $m\leq l$. From the equality $A A =
\Lambda$ we obtain for $m=l+1$
$$
d^{l+1} (A A) = (d^{l+1} A) A + \left(
                                  \begin{array}{c}
                                    l+1 \\
                                    1 \\
                                  \end{array}
                                \right) (d^{l} A) dA+ ...+ \left(
                                  \begin{array}{c}
                                    l+1 \\
                                    l \\
                                  \end{array}
                                \right)  dA (d^{l} A)
                              + A(d^{l+1} A)=0, $$
where $d$ denotes elementwise differentiation of a matrix with
respect to a fixed element of $\Lambda$. This implies
 \begin{equation} \label{eqadd2}
 (d^{l+1} A) A + A(d^{l+1} A) = - \left(
                                  \begin{array}{c}
                                    l+1 \\
                                    1 \\
                                  \end{array}
                                \right) (d^{l} A) dA- ...- \left(
                                  \begin{array}{c}
                                    l+1 \\
                                    l \\
                                  \end{array}
                                \right)  dA (d^{l} A).
 \end{equation}
 Denote the symmetric matrix in the right hand side of
 (\ref{eqadd2}) by $\widetilde \Lambda$. Then equality
 (\ref{eqadd2}) determines a linear operator $\ell$ mapping
$d^{l+1} A$ to $\widetilde \Lambda$. In the linear space of
symmetric $d\times d$ matrices we introduce the scalar product
$\langle X , Y \rangle = \mbox{trace} (XY)$. The operator $\ell$
determines a quadratic form
$$\langle \ell X , X\rangle = \mbox{trace}[(XA+AX)X] = 2
\mbox{trace}[XAX]\geq 2 \sqrt {\lambda^-} \mbox{trace}[XX] = 2 \sqrt
{\lambda^-} \langle  X , X\rangle,$$ where in the inequality we have
used that $ A-\sqrt {\lambda^-}I$ positive definite implies that
$X(A-\sqrt {\lambda^-}I)X=XAX-\sqrt {\lambda^-}XX$ is positive
definite. Similarly, we get $\langle \ell X , X\rangle \leq  2 \sqrt
{\lambda^+} \langle  X , X\rangle$. Hence,
$$2 \sqrt {\lambda^-} \leq \|\ell\| = \sup_{X\not = 0} {{\|\ell X\|}
\over {\|X\|}} \leq 2 \sqrt {\lambda^+}$$ and $${1 \over 2 \sqrt
{\lambda^+}} \leq \|\ell^{-1}\| \leq {1 \over 2 \sqrt
{\lambda^-}}.$$ We obtain $$ \|d^{l+1}A\| \leq {1 \over 2 \sqrt
{\lambda^-}} \|\widetilde \Lambda\|.$$ Using the induction
hypothesis we get from (\ref{eqadd2})
$$ \|d^{l+1}A\| \leq C_{l+1} ({\lambda^-})^{(2l+1)/2}.$$ This completes the proof.

\bigskip

From Lemmas 5 and 6 we get the following corollary. The statement of
\ the next lemma is an extension of Lemma 3.7 in Mammen and Konakov
(2000) where the result has been shown for $0\leq |b|\leq 2,a=0$.
\bigskip

\noindent \textbf{Lemma 7. }\textit{The following bound
holds:}\emph{}
\[
\left| D_{y}^{a}D_{x}^{b}\widetilde{p}_{h}(jh,kh,x,y)\right| \leq
C\rho ^{-\left| a\right| -\left| b\right| }\zeta _{\rho
}^{S-\left| a\right| }(y-x)
\]
\textit{for all }$j<k,$\textit{ for all }$x$\textit{ and
}$y$\textit{ and for all }$a,b$\textit{ with }$0\leq \left|
a\right| +\left| b\right| \leq 6.$\textit{ Here, }$\rho
=[(k-j)h]^{1/2}$\textit{. The constant }$S$\textit{ has been
defined in Assumption {(A3)}. }
\bigskip

\noindent \textsc{Proof of Lemma 7.} For two matrices $A$ and $B$
with elements $ a_{ij}$ or $ b_{kl} $, respectively where
$a_{ij}(B)$ are smooth functions of
$%
b_{kl} $ we write $\left| \frac{\partial A}{\partial B}\right| \leq
C$ if $\left| \frac{\partial a_{ij}}{\partial
b_{kl}}%
\right| \leq C$ for all $1\leq i,j\leq d,1\leq k,l\leq d.$ To
obtain the
assertion of the lemma we have to estimate the derivatives $%
D_{y}^{a}D_{x}^{b}z,$ where \emph{$z=V_{j,k}^{-1/2}(y)(y-x-\mu
_{j,k}(y)).$
} Note that \emph{$z=z(V_{j,k}^{-1/2},\mu _{j,k},x,y),$ }where $%
V_{j,k}^{-1/2}=V_{j,k}^{-1/2}(y) $ and \emph{}$\mu _{j,k}=\mu
_{j,k}(y)$. For $l=1,...,6$ it follows from condition ({B1}) and
(\ref{eq:003a}) that
\begin{equation}
\left| \frac{\partial^l \mu _{j,k}(y)}{(\partial y)^l}\right| \leq
C\rho ^{2},\left| \frac{\partial^l V_{j,k}(y)}{(\partial
y)^l}\right| \leq C\rho ^{2}. \label{eq:12a}
\end{equation}
It follows from Lemma 6 that
\begin{equation}
\left| \frac{\partial^l V_{j,k}^{1/2}}{(\partial V_{j,k})^l}\right|
\leq C \rho^{-(2l-1)/2}. \label{eq:12b}
\end{equation}
From inequalities (3.16) in Konakov and Mammen (2000) and from the
representation of an inverse matrix in terms of cofactors divided by
the determinant we obtain that
\begin{equation}
\left| \frac{\partial^l V_{j,k}^{-1/2}}{(\partial
V_{j,k}^{1/2})^l}\right| \leq C\rho ^{-(l+1)}.  \label{eq:12c}
\end{equation}
From (\ref{eq:12a})-(\ref{eq:12c}) and from the chain rule we get
\begin{equation}
\left| \frac{\partial^l V_{j,k}^{-1/2}(y)}{(\partial y)^l}\right|
\leq C \rho ^{-l}. \label{eq:12d}
\end{equation}
Now, Lemma 5 implies the assertion of Lemma 7.

\subsection{\label{sec:bonker}Bounds on operator kernels used in the
parametrix expansions.}

In this Section we will present bounds for operator kernels
appearing in the expansions based on the parametrix method. In
Lemma 8 we compare the infinitesimal operators $L_{h}$ and
$\tilde{L}_{h}$ with the differential operators $L$ and
$\tilde{L}$. We give an approximation for the error if, in the
definition of $H_{h}=(L_{h}-\tilde{L}_{h})\tilde{p}_{h}$, the
terms
$%
L_{h}$ and $\tilde{L}_{h}$ are replaced by $L$ or $\tilde{L}$,
respectively. We show that this term can be approximated by
$K_{h}+M_{h}$, where
$K_{h}=(L-%
\tilde{L})\tilde{p}_{h}$ and where $M_{h}$ is defined in Remark 5
after Lemma 8 . The bounds obtained in Lemma 9 will be used in the
proof of our theorem to show that in the expansion of $p_{h}$\ the
terms $\widetilde{p}_{h}\otimes _{h}H_{h}^{(r)}$\ can be replaced
by $\widetilde{p}_{h}\otimes _{h}(K_{h}+M_{h})^{(r)}.$
\bigskip

\noindent \textbf{Lemma 8.} \textit{The following bound holds with a
constant }$C$
\begin{eqnarray*}
&&\left| H_{h}(jh,kh,x,y)-K_{h}^{\prime }(jh,kh,x,y)-M_{h}^{\prime
}(jh,kh,x,y)-R_{h}(jh,kh,x,y)\right|
\\ && \qquad
 \leq Ch^{3/2}\rho ^{-1}\zeta
_{\rho }^{S}(y-x)
\end{eqnarray*}
\textit{with }$\zeta _{\rho }^{S}$\textit{ as in Lemma 5 for all
}$j<k$%
\textit{, }$x$\textit{ and }$y$\textit{. For $j {<} k-1$ we
define}
\begin{eqnarray*}
&& K_{h}^{\prime }(jh,kh,x,y)=(L-\widetilde{L})\lambda
(x),M_{h}^{\prime }(jh,kh,x,y) \\ && \qquad =M_{h,1}(jh,kh,x,y)
+M_{h,2}(jh,kh,x,y)+M_{h,3}^{\prime }(jh,kh,x,y),
\\ &&
M_{h,1}(jh,kh,x,y)=h^{1/2}\sum_{\left| \nu \right|
=3}\frac{D_{x}^{\nu }\lambda (x)}{\nu !}(\chi _{\nu }(jh,x)-\chi
_{\nu }(jh,y)),
\\ &&
M_{h,2}(jh,kh,x,y)=h\sum_{\left| \nu \right| =4}\frac{D_{x}^{\nu
}\lambda (x)%
}{\nu !}(\chi _{\nu }(jh,x)-\chi _{\nu }(jh,y)),
\\ &&
M_{h,3}^{\prime }(jh,kh,x,y)=\frac{h}{2}(L_{\star }^{2}-\widetilde{L}%
^{2})\lambda (x),
\\ &&
R_{h}(jh,kh,x,y)=h^{3/2}\sum_{\left| \nu \right| =4}\frac{D_{x}^{\nu
}\lambda (x)}{\nu !}\sum_{r=1}^{d}\nu _{r}[m_{r}(jh,x)\mu _{\nu
-e_{r}}(jh,x)-m_{r}(jh,y)\mu _{\nu -e_{r}}(jh,y)]
\\ && \qquad
+5\sum_{\left| \nu \right| =5}\frac{1}{\nu !}%
\sum_{k=1}^{d}(m_{k}(jh,x)-m_{k}(jh,y))\left\{ \nu _{k}\int
q(jh,x,\theta )%
\widetilde{h}^{\nu -e_{k}}(\theta )\right.
\\ && \qquad
\left. \times \left[ \int_{0}^{1}(1-u)^{4}D^{\nu }\lambda
(x+u\widetilde{h}%
(\theta ))du\right] d\theta +\int q(jh,x,\theta
)\widetilde{h}^{\nu }(\theta )\left[ \int_{0}^{1}(1-u)^{4}uD^{\nu
+e_{k}}\lambda (x+u\widetilde{h}(\theta ))du\right] d\theta
\right\}
\\ && \qquad
+h^{2}\sum_{\left| \nu \right| =4}\frac{D_{x}^{\nu }\lambda (x)}{\nu
!}%
\sum_{\left| \nu ^{\prime }\right| =2}\nu !N(\nu ,\nu ^{\prime
})[m^{\nu ^{\prime }}(jh,x)\mu _{\nu -\nu ^{\prime }}(jh,x)-m^{\nu
^{\prime }}(jh,y)\mu _{\nu -\nu ^{\prime }}(jh,y)].
\end{eqnarray*}
\textit{Here }$L_{\star }$\textit{ is defined as }$%
\widetilde{L}$\textit{ but with the coefficients ''frozen'' at the
point x, e}$_{r}$\textit{ denotes a }$d$\textit{-dimensional
vector with the $r$-th element equal to }$1$\textit{ and with all
other elements equal to }$0.$\textit{ Furthermore, for }$\left|
\nu \right| =4,\left| \nu ^{\prime }\right| =2$\textit{ we define}
\[
N(\nu ,\nu ^{\prime })=2^{\chi \lbrack \nu ^{\prime }!=1]+\chi
\lbrack (\nu -\nu ^{\prime })!=1]-2},
\]
\textit{where }$\chi (\cdot )$\textit{ is the indicator function. We
put }%
$m(x)^{\nu }=m_{1}(x)^{\nu _{1}}\cdot ...\cdot m_{d}(x)^{\nu
_{d}}$\textit{ and }$m(x)^{\nu }=0,$\textit{ }$\nu !=0$. \textit{
We define }$\mu_{\nu}(t,x) = \int z^{\nu} q(t,x,z) dz$ \textit{
and }$\mu _{\nu
}(t,x)=0$%
\textit{ if at least one of the coordinates of }$\nu =(\nu
_{1},...,\nu _{d})$\textit{ is negative. We use also the following
definitions}
\begin{eqnarray*}
\lambda (x)&=&\widetilde{p}_{h}((j+1)h,kh,x,y),
\\
\widetilde{h}(\theta )&=&m(jh,y)h+\theta h^{1/2}.
\end{eqnarray*}
\textit{Here again }$\rho $\textit{ denotes the term }$\rho =\left[
h(k-j)%
\right] ^{1/2}.$\textit{ For }$j=k-1$\textit{ we define}
\[
K_{h}^{\prime
}(jh,kh,x,y)=R_{h}(jh,kh,x,y)=M_{h,2}(jh,kh,x,y)=M_{h,3}^{\prime
}(jh,kh,x,y)=0
\]
and $$M_{h,1}(jh,kh,x,y) = h^{-(d+2)/2} \left[
q\left\{jh,x,h^{-1/2}(y-x-m[jh,x]h)\right\} -
q\left\{jh,y,h^{-1/2}(y-x-m[jh,y]h)\right\}\right].$$
\bigskip

\noindent \textsc{Proof of Lemma 8. }As in the proof of Lemma 3.9 \
in Konakov and Mammen (2000) we have
\[
H_{h}(jh,kh,x,y)=H_{h}^{1}(jh,kh,x,y)-H_{h}^{2}(jh,kh,x,y),
\]
where
\begin{eqnarray}
&&H_{h}^{1}(jh,kh,x,y)=h^{-1}\int q(jh,x,\theta )[\lambda
(x+h(\theta ))-\lambda (x)]d\theta,   \label{eq:12g}
\\ &&
H_{h}^{2}(jh,kh,x,y)=h^{-1}\int q(jh,y,\theta )[\lambda
(x+\widetilde{h} (\theta ))-\lambda (x)]d\theta ,  \label{eq:12h}
\\ \nonumber
&&h(\theta )=m(jh,x)h+\theta h^{1/2},\widetilde{h}(\theta
)=m(jh,y)h+\theta h^{1/2}.
\end{eqnarray}
For $[\lambda (x+h(\theta ))-\lambda (x)]$ and $[\lambda
(x+\widetilde{%
h}(\theta ))-\lambda (x)]$ in (\ref{eq:12g}), (\ref{eq:12h}) we
use now the Taylor expansion up to order 5 with remaining term in
integral form. To pass from moments to cumulants we use the well
known relations (see e.g. relation (6.11) on page 46 in
Bhattacharya and Rao (1986)). After long but simple calculations
we come to the conclusion of the lemma. \bigskip

\noindent \textbf{Remark 5. }We show now that the function
$K_{h}^{\prime }(jh,kh,x,y)+M_{h,3}^{\prime }(jh,kh,x,y)$ in Lemma 8
is equal to $K_{h}(jh,kh,x,y)+\frac{h}{2}(L_{\star
}^{2}-2L\widetilde{L}+\widetilde{L}%
^{2})\lambda (x)+M_{h,3}^{\prime \prime }(jh,kh,x,y)$ where
\begin{eqnarray} \label{eq:012h}
M_{h,3}^{\prime \prime }(jh,kh,x,y)&=&-h^{2}\sum_{\left| \mu \right|
=2}\frac{%
m^{\mu }(jh,y)}{\mu !}(L-\widetilde{L})D^{\mu }\lambda (x)
\\ \nonumber &&-3\sum_{\left| \mu \right| =3}\int_{0}^{1}(1-\delta
)^{2}d\delta \int q(jh,y,\theta )\frac{\widetilde{h}(\theta )^{\mu
}}{\mu
!}(L-\widetilde{L}%
)D^{\mu }\lambda (x+\delta \widetilde{h}(\theta ))d\theta .
\end{eqnarray}
Thus in Lemma 8 we can replace $K_{h}^{\prime
}(jh,kh,x,y)+M_{h}^{\prime
}(jh,kh,x,y)$ by $K_{h}(jh,kh,x,y)+M_{h}(jh,kh,x,y)$ where $%
K_{h}(jh,kh,x,y)=(L-\widetilde{L})\widetilde{p}_{h}(jh,kh,x,y),$ $%
M_{h}(jh,kh,x,y)=\frac{h}{2}(L_{\star
}^{2}-2L\widetilde{L}+\widetilde{L}%
^{2})\lambda (x)+M_{h}^{\prime \prime },M_{h}^{\prime \prime
}=M_{h,1}(jh,kh,x,y)+M_{h,2}(jh,kh,x,y)+M_{h,3}^{\prime \prime
}(jh,kh,x,y)$ and
\[
\max \{\left| M_{h}^{\prime }(jh,kh,x,y)\right| ,\left|
M_{h}(jh,kh,x,y)\right| \}\leq C\rho ^{-1}\zeta _{\rho }(y-x),
\]
$\rho ^{2}=kh-jh.$ To show this we note that
\[
\widetilde{p}_{h}(jh,kh,x,y)=\int q(jh,y,\theta )\lambda
(x+\widetilde{h}%
(\theta ))d\theta,
\]
where $\widetilde{h}(\theta )=m(jh,y)h+h^{1/2}\theta .$ From the
Taylor expansion we get
\begin{eqnarray*}
\widetilde{p}_{h}(jh,kh,x,y)&=&\lambda (x)+h\widetilde{L}\lambda
(x)+h^{2}\sum_{\left| \mu \right| =2}\frac{m^{\mu }(jh,y)}{\mu
!}D^{\mu }\lambda (x)
\\&&
+3\sum_{\left| \mu \right| =3}\int_{0}^{1}(1-\delta )^{2}d\delta
\int q(jh,y,\theta )\frac{\widetilde{h}(\theta )^{\mu }}{\mu
!}D^{\mu }\lambda (x+\delta \widetilde{h}(\theta ))d\theta
\end{eqnarray*}
and, hence,
\begin{eqnarray} \nonumber
&&K_{h}^{\prime
}(jh,kh,x,y)=K_{h}(jh,kh,x,y)+(L-\widetilde{L})[\lambda
(x)-%
\widetilde{p}_{h}(jh,kh,x,y)]
\\ && \qquad
=K_{h}(jh,kh,x,y)+h(\widetilde{L}^{2}-L\widetilde{L})\lambda
(x)+M_{h,3}^{\prime \prime }(jh,kh,x,y).  \label{eq:12i}
\end{eqnarray}
From \begin{eqnarray*} h(\widetilde{L}^{2}-L\widetilde{L})\lambda
(x)+M_{h,3}^{\prime
}(jh,kh,x,y)&=&h(\widetilde{L}^{2}-L\widetilde{L})\lambda
(x)+\frac{h}{2}%
(L_{\star }^{2}-\widetilde{L}^{2})\lambda (x)
\\&=&\frac{h}{2}(L_{\star
}^{2}-2L\widetilde{L}+\widetilde{L}^{2})\lambda (x)
\end{eqnarray*}
and from the definitions of the operators $L,\widetilde{L}$ and
$%
L_{\star }$ and from the Lipschitz conditions on the coefficients
$m(t,x)$
and $%
\sigma (t,x)$ we obtain that
\begin{equation}
\left| \frac{h}{2}(L_{\star
}^{2}-2L\widetilde{L}+\widetilde{L}^{2})\lambda (x)\right| \leq
Ch\rho ^{-3}\zeta _{\rho }(y-x).  \label{eq:12j}
\end{equation}
Analogously, we have
\begin{eqnarray}
&&\left| h^{2}\sum_{\left| \mu \right| =2}\frac{m^{\mu }(jh,y)}{\mu
!}(L-%
\widetilde{L})D^{\mu }\lambda (x)\right| \leq Ch^{2}\rho
^{-3}\zeta _{\rho }(y-x),  \label{eq:12k}
\\ &&
\left| 3\sum_{\left| \mu \right| =3}\int_{0}^{1}(1-\delta
)^{2}d\delta \int q(jh,y,\theta )\frac{\widetilde{h}(\theta )^{\mu
}}{\mu
!}(L-\widetilde{L}%
)D^{\mu }\lambda (x+\delta \widetilde{h}(\theta ))d\theta \right|
\leq Ch^{3/2}\rho ^{-4}\zeta _{\rho }(y-x).  \label{eq:12l}
\end{eqnarray}
Now (\ref{eq:12i})-(\ref{eq:12l}) imply the assertion of this
remark.
\bigskip

\noindent \textbf{Lemma 9. }\textit{The following bound holds:}
\begin{eqnarray} \nonumber &&
\left| \sum_{r=0}^{n}\widetilde{p}_{h}\otimes
_{h}(K_{h}+M_{h}+R_{h})^{(r)}(0,T,x,y)-\sum_{r=0}^{n}\widetilde{p}%
_{h}\otimes _{h}(K_{h}+M_{h})^{(r)}(0,T,x,y)\right|
\\ \label{eq:13}&& \qquad
\leq C(\varepsilon )hn^{-1/2+\varepsilon }\zeta
_{\sqrt{T}}^{S}(y-x),
\end{eqnarray}
\textit{where} $\lim_{\varepsilon \downarrow 0}C(\varepsilon
)=+\infty .$\bigskip

\noindent
\textsc{Proof of Lemma 9. }For $r=1$ we will show that for any $%
\varepsilon >0$%
\begin{eqnarray}
&&\left| \widetilde{p}_{h}\otimes
_{h}(K_{h}+M_{h}+R_{h})(0,kh,x,y)-\widetilde{%
p}_{h}\otimes _{h}(K_{h}+M_{h})(0,kh,x,y)\right| \nonumber \\ &&
=\left| \widetilde{p}_{h}\otimes _{h}R_{h}(0,kh,x,y)\right| \leq
Ch^{3/2-\varepsilon }(kh)^{-1/2+\varepsilon
}B(\frac{1}{2},\varepsilon )\zeta _{\rho }^{S}(y-x),\rho ^{2}=kh.
\label{eq:14}
\end{eqnarray}
Clearly, to estimate $\widetilde{p}_{h}\otimes _{h}R_{h}(0,T,x,y)$ \
it is enough to estimate
\[
I_{1}=h^{3/2}\sum_{j=0}^{k-2}h\int \widetilde{p}%
_{h}(0,kh,x,z)(f(jh,z)-f(jh,y))D_{z}^{\nu }\widetilde{p}%
_{h}((j+1)h,kh,z,y)dz
\]
for $\nu ,$ $\left| \nu \right| =4,$ and
\begin{eqnarray*}
I_{2}&=&h^{2}\sum_{j=0}^{k-2}h\int \widetilde{p}%
_{h}(0,jh,x,z)(f(jh,z)-f(jh,y))\int q(jh,z,\theta
)\widetilde{h}^{\nu -e_{k}}(\theta )
\\ &&
\times \int_{0}^{1}(1-u)^{4}D_{z}^{\nu }\lambda
(z+u\widetilde{h}(\theta ))dud\theta dz
\end{eqnarray*}
for $\nu ,\left| \nu \right| =5,$ $1\leq k\leq d.$ Here $f(t,x)$ is
a function whose first and second derivatives with respect to $x$
are continuous and bounded uniformly in $t$ and $x.$ After
integration by parts we obtain \begin{eqnarray*}
I_{1}&=&-h^{3/2}\sum_{j=0}^{k-2}h\int D_{z}^{e_{l}}\widetilde{p}%
_{h}(0,jh,x,z)(f(jh,z)-f(jh,y))D_{z}^{\nu -e_{l}}\widetilde{p}%
_{h}((j+1)h,kh,z,y)dz
\\&&
+h^{3/2}\sum_{j=0}^{k-2}h\int D_{z}^{e_{s}}\widetilde{p}%
_{h}(0,jh,x,z)D_{z}^{e_{k}}f(jh,z)D_{z}^{\nu
-e_{l}-e_{s}}\widetilde{p}%
_{h}((j+1)h,kh,z,y)dz
\\&&
+h^{3/2}\sum_{j=0}^{k-2}h\int \widetilde{p}%
_{h}(0,jh,x,z)D_{z}^{e_{l}+e_{s}}f(jh,z)D_{z}^{\nu
-e_{l}-e_{s}}\widetilde{p}%
_{h}((j+1)h,kh,z,y)dz
\end{eqnarray*}
for $1\leq l,s\leq d$. Hence,
\begin{equation}
\left| I_{1}\right| \leq
Ch^{3/2}\sum_{j=0}^{k-2}h\frac{1}{\sqrt{jh}(kh-jh)}%
\zeta _{\rho }(y-x)\leq Ch^{3/2-\varepsilon
}(kh)^{-1/2+\varepsilon
}B(\frac{%
1}{2},\varepsilon )\zeta _{\sqrt{T}}^{S}(y-x).  \label{eq:15}
\end{equation}
In the same way after integration by parts we get with $1\leq
l,s\leq d.$
\begin{eqnarray} \nonumber
I_{2}&=&-h^{2}\sum_{j=0}^{k-2}h\int_{0}^{1}(1-u)^{4}du\int d\theta
(m(jh,y)h^{1/2}+\theta )^{\nu -e_{l}}\int D_{z}^{e_{l}}\widetilde{p}%
_{h}(0,jh,x,z)(f(jh,z)-f(jh,y))
\\ && \nonumber \qquad
\times q(jh,z,\theta )D_{z}^{\nu
-e_{l}}\widetilde{p}_{h}((j+1)h,kh,z+u%
\widetilde{h}(\theta
),y)dz+h^{2}\sum_{j=0}^{k-2}h\int_{0}^{1}(1-u)^{4}du\int d\theta
(m(jh,y)h^{1/2}+\theta )^{\nu -e_{l}}
\\ && \nonumber \qquad
\times \int
D_{z}^{e_{s}}[\widetilde{p}_{h}(0,jh,x,z)D^{e_{l}}f(jh,z)q(jh,z,%
\theta )]D_{z}^{\nu
-e_{l}-e_{s}}\widetilde{p}_{h}((j+1)h,kh,z+u\widetilde{h}%
(\theta ),y)dz
\\ && \nonumber
-h^{2}\sum_{j=0}^{k-2}h\int_{0}^{1}(1-u)^{4}du\int d\theta
(m(jh,y)h^{1/2}+\theta )^{\nu -e_{l}}\widetilde{p}%
_{h}(0,jh,x,z)(f(jh,z)-f(jh,y))
\\ && \qquad
\times D_{z}^{e_{l}}q(jh,z,\theta )D_{z}^{\nu -e_{l}}\widetilde{p}%
_{h}((j+1)h,kh,z+u\widetilde{h}(\theta ),y)dz.  \label{eq:16}
\end{eqnarray}
It follows from (\ref{eq:16}) that
\begin{equation}
\left| I_{2}\right| \leq Ch^{3/2-\varepsilon
}(kh)^{-1/2+\varepsilon
}B(%
\frac{1}{2},\varepsilon )\zeta _{\rho }^{S}(y-x).  \label{eq:17}
\end{equation}
Claim (\ref{eq:14}) follows now from (\ref{eq:15}) and
(\ref{eq:17}).
For $%
r\geq 2$ we use the identity
\begin{eqnarray} \nonumber
&& \widetilde{p}_{h}\otimes
_{h}(K_{h}+M_{h}+R_{h})^{(r)}(0,T,x,y)-\widetilde{p}%
_{h}\otimes _{h}(K_{h}+M_{h})^{(r)}(0,T,x,y)
\\ \nonumber
&& =\left[ \widetilde{p}_{h}\otimes
_{h}(K_{h}+M_{h}+R_{h})^{(r-1)}-\widetilde{p%
}_{h}\otimes _{h}(K_{h}+M_{h})^{(r-1)}\right] \otimes
_{h}(K_{h}+M_{h})(0,T,x,y) \\ \nonumber && \qquad \qquad
+\widetilde{p}_{h}\otimes _{h}(K_{h}+M_{h}+R_{h})^{(r-1)}\otimes
_{h}R_{h}(0,T,x,y)\\
 && =I+II. \label{eq:18}
\end{eqnarray}
For $r=2$ we obtain from (\ref{eq:14}) and simple estimate $\left|
(K_{h}+M_{h})(jh,kh,z,y)\right| \leq $ $C\rho _{2}^{-1}\zeta _{\rho
_{2}}^{S}(y-z),\rho _{2}^{2}=kh-jh,$%
\begin{eqnarray*}
\left| I\right| &=&\left| \left[ \widetilde{p}_{h}\otimes
_{h}(K_{h}+M_{h}+R_{h})-\widetilde{p}_{h}\otimes
_{h}(K_{h}+M_{h})\right] \otimes _{h}(K_{h}+M_{h})(0,kh,x,y)\right|
\\ &&
\leq C^{2}h^{3/2-\varepsilon }B(\frac{1}{2},\varepsilon
)\sum_{j=0}^{k-2}h(jh)^{-1/2+\varepsilon }(kh-jh)^{-1/2}\int \zeta
_{\rho _{1}}^{S}(z-x)\zeta _{\rho _{2}}^{S}(y-z)dz
\\ &&
\leq C^{2}h^{3/2-\varepsilon }B(\frac{1}{2},\varepsilon )B(\frac{1}{2}%
,\varepsilon +\frac{1}{2})(kh)^{\varepsilon }\zeta _{\rho }^{S}(y-x)
\end{eqnarray*}
with $\rho ^{2}=kh$. For $r\geq 3$ we obtain by induction
\begin{eqnarray} \nonumber
\left| I\right| &=&\left| \left[ \widetilde{p}_{h}\otimes
_{h}(K_{h}+M_{h}+R_{h})^{(r-1)}-\widetilde{p}_{h}\otimes
_{h}(K_{h}+M_{h})^{(r-1)}\right] \otimes
_{h}(K_{h}+M_{h})(0,kh,x,y)\right|
\\ \nonumber
&\leq& C^{r}h^{3/2-\varepsilon }B(\frac{1}{2},\varepsilon )B(\frac{1}{2}%
,\varepsilon +\frac{1}{2})...B(\frac{1}{2},\varepsilon +\frac{r-1}{2}%
)(kh)^{\varepsilon +(r-2)/2}\zeta _{\rho }^{S}(y-x)
\\
&\leq& \Gamma (\varepsilon )h^{3/2-\varepsilon }\frac{[C\Gamma
(1/2)]^{r}}{%
\Gamma (\varepsilon +\frac{r}{2})}(kh)^{\varepsilon +(r-2)/2}\zeta
_{\rho }^{S}(y-x) \label{eq:19}
\end{eqnarray}
with $\rho ^{2}=kh$. To estimate $II$  we use the following
estimates
\begin{eqnarray}
&&\left| D_{y}^{a}D_{x}^{b}\widetilde{p}_{h}(jh,kh,x,y)\right|
\leq C\rho ^{-\left| a\right| -\left| b\right| }\zeta _{\rho
}^{S-\left| a\right| }(y-x),\left|
D_{x}^{b}\widetilde{p}_{h}(jh,kh,x,x+v)\right|\leq C\zeta _{\rho
}^{S}(v), \label{eq:20}
\\ &&
\left| D_{x}^{b}(K_{h}+M_{h}+R_{h})(jh,kh,x+v,x)\right| \leq C\rho
^{-1}\zeta _{\rho }^{S}(v).  \label{eq:21}
\end{eqnarray}
The inequalities (\ref{eq:20}) and (\ref{eq:21}) are obtained by
using the same arguments as is the proof of Lemma 7. Using these
inequalities and mimicking the proof of Theorem 2.3 in Konakov and
Mammen (2002) we obtain the following bounds for $r=0,1,...$
\begin{eqnarray} \nonumber
&&\left| D_{x}^{b}D_{y}^{a}\widetilde{p}_{h}\otimes
_{h}(K_{h}+M_{h}+R_{h})^{(r)}(0,kh,x,y)\right|
\\ && \qquad \nonumber
\leq C^{r}(kh)^{-\left| a\right| -\left| b\right|
+r}B(\frac{1}{2},\frac{1}{2%
})B(1,\frac{1}{2})...B(\frac{r}{2},\frac{1}{2})\zeta _{\rho
}^{S-\left| a\right| }(y-x)
\\ && \qquad
\leq \frac{\lbrack C\Gamma (1/2)]^{r}}{\Gamma
(\frac{r+1}{2})}(kh)^{-\left| a\right| -\left| b\right| +r}\zeta
_{\rho }^{S-\left| a\right| }(y-x). \label{eq:22}
\end{eqnarray}
Inequality (\ref{eq:22}) allows us to estimate $II=[\widetilde{p}%
_{h}\otimes _{h}(K_{h}+M_{h}^{\prime \prime }+R_{h})^{(r-1)}]\otimes
_{h}R_{h}(0,kh,x,y)$. For this it is enough to estimate
\begin{eqnarray} \nonumber
&& h^{3/2}\sum h\int [\widetilde{p}_{h}\otimes
_{h}(K_{h}+M_{h}+R_{h})^{(r-1)}](0,jh,x,z)
\\ && \qquad
\times D^{\nu }\widetilde{p}_{h}((j+1)h,kh,z,y)(f(jh,z)-f(jh,y))dz
\label{eq:23}
\end{eqnarray}
for $r\geq 2$ ,$\left| \nu \right| =4,$ and
\begin{eqnarray} \nonumber &&
\sum_{j=0}^{n-2}h\int [\widetilde{p}_{h}\otimes
_{h}(K_{h}+M_{h}+R_{h})^{(r-1)}](0,jh,x,z)(f(jh,z)-f(jh,y))
\\ && \qquad
\times \int q(jh,z,\theta )\widetilde{h}^{\nu -e_{l}}(\theta
)\int_{0}^{1}(1-u)^{4}D^{\nu
}\widetilde{p}_{h}((j+1)h,kh,z+u\widetilde{h}%
(\theta ),y)dud\theta dz  \label{eq:23a}
\end{eqnarray}
for $r\geq 2,$ $\left| \nu \right| =5,$ $1\leq l\leq d.$ Here
$f(t,x)$ is a function whose first and second derivatives with
respect to $x$ are continuous and bounded uniformly in $t$ and $x$.
The upper bound for (\ref{eq:23}) follows from (\ref{eq:22}) by
integration by parts exactly in the same way as it was done to
obtain the upper bound for $I_{1}$, see (\ref{eq:15}). This gives
the estimate
\begin{eqnarray} \nonumber &&
h^{3/2}\Big| \sum h\int [\widetilde{p}_{h}\otimes
_{h}(K_{h}+M_{h}+R_{h})^{(r-1)}](0,jh,x,z)
\\  \nonumber && \qquad
 \times D^{\nu
}\widetilde{p}_{h}((j+1)h,kh,z,y)(f(jh,z)-f(jh,y))dz%
\Big|
\\  && \qquad
\leq \Gamma (\varepsilon )h^{3/2-\varepsilon }\frac{[C\Gamma
(1/2)]^{r}}{%
\Gamma (\frac{r+1}{2})}(kh)^{\varepsilon +(r-2)/2}\zeta _{\rho
}^{S}(y-x). \label{eq:24}
\end{eqnarray}
The upper bound for (\ref{eq:23a}) follows from (\ref{eq:22}) by
integration by parts in the same way as it was done to obtain an
upper bound for $I_{2}$, see (\ref{eq:17}). This gives for
(\ref{eq:23a}) the same estimate as in (\ref{eq:24}) and, hence,
\begin{equation}
\left| II\right| \leq C\Gamma (\varepsilon )h^{3/2-\varepsilon }\frac{%
[C\Gamma (1/2)]^{r}}{\Gamma (\frac{r+1}{2})}(kh)^{\varepsilon
+(r-2)/2}\zeta _{\rho }^{S}(y-x).  \label{eq:25}
\end{equation}
The assertion of the lemma follows now from (\ref{eq:14}),
(\ref{eq:18}), (\ref{eq:19}) and (\ref{eq:25}).
\bigskip

\noindent \textbf{Lemma 10. }\textit{Let
}$A(s,t,x,y),B(s,t,x,y),C(s,t,x,y)$\textit{%
  be some functions with absolute value less than
}$C(t-s)^{-1/2}\zeta^{\kappa} _{\sqrt{t-s}}(y-x)$\textit{ for a
constant }$C$ \textit{ and an integer } $\kappa \geq S^\prime d$.
\textit{ Then}
\[
\sum_{r=0}^{\infty }A\otimes
_{h}(B+C)^{(r)}(ih,jh,x,y)-\sum_{r=0}^{\infty }A\otimes
_{h}B^{(r)}(ih,jh,x,y)
\]
\[
=\sum_{r=1}^{\infty }\left[ A\otimes _{h}\Phi \right] \otimes
_{h}\left[ C\otimes _{h}\Phi \right] ^{(r)}(ih,jh,x,y),
\]
\textit{where }$\Phi =\sum_{r=0}^{\infty }B^{(r)}.$ \bigskip

\noindent \textsc{Proof of Lemma 10. }Under the conditions of the
lemma all series are absolutely convergent. The assertion of this
lemma is a consequence of the linearity of the operation $\otimes
_{h}$ and of the possibility to permutate the terms in absolutely
convergent series.\bigskip

\section{\label{sec:mainres} Proof of Theorem 1.}

We now come to the proof of Theorem 1. Main tools for the proof have
been given in Subsections \ref{sec:Pardiff}, \ref{sec:Parchain},
\ref{sec:bonedg}  and \ref{sec:bonker}. From Lemmas 1 and 2 we get
that
\[
p(0,T,x,y)=\sum_{r=0}^{n}\tilde{p}\otimes
H^{(r)}(0,T,x,y)+o(h^{2}T)\phi _{C,%
\sqrt{T}}(y-x).
\]
With Lemma 3 this gives
\begin{equation}
p(0,T,x,y)-p_{h}(0,T,x,y)=T_{1}+...+T_{7}+o(h^{2}T)\phi
_{C,\sqrt{T}}(y-x), \label{eq:26}
\end{equation}
where
\begin{eqnarray*}
T_{1}&=&\sum_{r=0}^{n}\tilde{p}\otimes
H^{(r)}(0,T,x,y)-\sum_{r=0}^{n}\tilde{p}%
\otimes _{h}H^{(r)}(0,T,x,y),
\\
T_{2}&=&\sum_{r=0}^{n}\tilde{p}\otimes
_{h}H^{(r)}(0,T,x,y)-\sum_{r=0}^{n}%
\tilde{p}\otimes _{h}(H+M_{h}^{\prime \prime }+\sqrt{h}%
N_{1})^{(r)}(0,T,x,y),
\\
T_{3}&=&\sum_{r=0}^{n}\tilde{p}\otimes _{h}(H+M_{h}^{\prime \prime
}+\sqrt{h}%
N_{1})^{(r)}(0,T,x,y)-\sum_{r=0}^{n}\tilde{p}\otimes
_{h}(H+M_{h}+\sqrt{h}%
N_{1})^{(r)}(0,T,x,y),
\\
T_{4}&=&\sum_{r=0}^{n}\tilde{p}\otimes _{h}(H+M_{h}+\sqrt{h}%
N_{1})^{(r)}(0,T,x,y)-\sum_{r=0}^{n}\tilde{p}\otimes
_{h}(K_{h}+M_{h})^{(r)}(0,T,x,y),
\\
T_{5}&=&\sum_{r=0}^{n}\tilde{p}\otimes
_{h}(K_{h}+M_{h})^{(r)}(0,T,x,y)-\sum_{r=0}^{n}\tilde{p}_{h}\otimes
_{h}(K_{h}+M_{h})^{(r)}(0,T,x,y),
\\
T_{6}&=&\sum_{r=0}^{n}\tilde{p}_{h}\otimes
_{h}(K_{h}+M_{h})^{(r)}(0,T,x,y)-\sum_{r=0}^{n}\tilde{p}_{h}\otimes
_{h}(K_{h}+M_{h}+R_{h})^{(r)}(0,T,x,y),
\\
T_{7}&=&\sum_{r=0}^{n}\tilde{p}_{h}\otimes
_{h}(K_{h}+M_{h}+R_{h})^{(r)}(0,T,x,y)-\sum_{r=0}^{n}\tilde{p}_{h}\otimes
_{h}H_{h}^{(r)}(0,T,x,y).
\end{eqnarray*}
Here we put $N_{1}(s,t,x,y)=(L-\widetilde{L})\widetilde{\pi
}_{1}(s,t,x,y).$

We now discuss the asymptotic behaviour of the terms
$T_{1},...,T_{7}$.\bigskip

\noindent \textit{Asymptotic treatment of the term $T_{1}$.}

We start from the recurrence relations for $r=1,2,3,...$%
\[
\left(\widetilde{p}\otimes
H^{\left(r\right)}\right)\left(0,jh,x,y\right)-\left(\widetilde{p}%
\otimes_{h}H^{\left(r\right)}\right)\left(0,jh,x,y\right)\qquad\;\;\;\;
\]

\[
=\left[ \left( \widetilde{p}\otimes H^{\left( r-1\right) }\right)
\otimes
H-\left( \widetilde{p}\otimes H^{\left( r-1\right) }\right) \otimes _{h}H%
\right] \left( 0,jh,x,y\right) \qquad
\]
\begin{equation}
\,\,\,\,+\left[ \left( \widetilde{p}\otimes H^{\left( r-1\right)
}\right) \ -\left( \widetilde{p}\otimes _{h}H^{\left( r-1\right)
}\right) \right] \otimes _{h}H\left( 0,jh,x,y\right) .\qquad
\label{MR-006}
\end{equation}
By summing up the identities in (\ref{MR-006}) from $r=1$ to $\infty
$ and by using the linearity of the operations $\otimes $ and
$\otimes _{h}$ we get
\[
(p-p^{d})\left( 0,jh,x,y\right) =\left( p\otimes H-p\otimes
_{h}H\right) \left( 0,jh,x,y\right)
\]
\begin{equation}
+(p-p^{d})\otimes _{h}H\left( 0,jh,x,y\right), \label{MR-007}
\end{equation}
where we put
\begin{equation}
p^{d}(ih,i^{\prime }h,x,y)=\sum_{r=0}^{\infty }(\widetilde{p}\otimes
_{h}H^{(r)})(ih,i^{\prime }h,x,y).  \label{MR-008}
\end{equation}
By iterative application of (\ref{MR-007}) we obtain
\[
(p-p^{d})\left( 0,jh,x,y\right) =\left( p\otimes H-p\otimes
_{h}H\right) \left( 0,jh,x,y\right)
\]
\begin{equation}
+\left( p\otimes H-p\otimes _{h}H\right) \otimes _{h}\Phi \left(
0,jh,x,y\right) ,\qquad \,\,\,\,\qquad \;\;  \label{MR-009}
\end{equation}
where $\Phi (ih,i^{\prime }h,z,z^{\prime })=H(ih,i^{\prime
}h,z,z^{\prime })+H\otimes _{h}H(ih,i^{\prime }h,z,z^{\prime
})+...=\sum_{r=1}^{\infty }H^{(r)}(ih,i^{\prime }h,z,z^{\prime }).$

By application of a Taylor expansion we get
\[
\left(p\otimes
H-p\otimes_{h}H\right)(0,jh,x,z)\qquad\qquad\qquad\qquad\qquad\qquad\;\;
\]
\[
=\sum_{i=0}^{j-1}\int_{ih}^{(i+1)h}du\int_{R^{d}}\left[\lambda\left(u%
\right)-\lambda\left(ih\right)\right]dv\qquad\qquad\qquad\qquad\;\;\,\,
\]
\[
=\sum_{i=0}^{j-1}\int_{ih}^{(i+1)h}(u-ih)du\int_{R^{d}}\lambda^{%
\prime}(ih)dv\qquad\qquad\qquad\qquad\qquad
\]
\begin{equation}
+\sum_{i=0}^{j-1}\int_{ih}^{(i+1)h}(u-ih)^{2}\int_{0}^{1}(1-%
\delta)\int_{R^{d}}\lambda^{\prime\prime}(s)\mid_{s=s_{i}}dvd\delta
du, \label{MR-010}
\end{equation}
where $\lambda\left(u\right)=p(0,u,x,v)H(u,jh,v,z),$ $s_{i}=s_{i}(u,i,%
\delta,h)=ih+\delta(u-ih).$

Note that
\[
\int_{R^{d}}\lambda^{\prime}(ih)dv=\int_{R^{d}}\frac{\partial}{\partial s}%
p(0,s,x,v)\mid_{s=ih}H(ih,jh,v,z)dv\qquad\;\;\;\;\;
\]
\[
+\int_{R^{d}}p(0,ih,x,v)\frac{\partial}{\partial s}H(s,jh,v,z)\mid_{s=ih}dv=%
\int_{R^{d}}L^{t}p(0,ih,x,v)\;\;
\]
\[
\times(L-\widetilde{L})\widetilde{p}(ih,jh,v,z)dv-\int_{R^{d}}p(0,ih,x,v)[(L-%
\widetilde{L})\widetilde{L}\widetilde{p}(ih,jh,v,z)
\]
\[
-H_{1}(ih,jh,v,z)]dv=\int_{R^{d}}p(0,ih,x,v)H_{1}(ih,jh,v,z)dv\qquad\qquad\;%
\;
\]
\begin{equation}
+\int_{R^{d}}p(0,ih,x,v)(L^{2}-2L\widetilde{L}+\widetilde{L}^{2})\widetilde{p%
}(ih,jh,v,z)dv,\qquad\qquad\qquad\;\;  \label{MR-011}
\end{equation}
where $H_{1}(s,t,v,z)$ is defined below in (\ref{MR-015}). We get
from (\ref {MR-011})
\[
\sum_{i=0}^{j-1}\int_{ih}^{(i+1)h}(u-ih)du\int_{R^{d}}\lambda^{\prime}(ih)dv=%
\frac{h}{2}(p\otimes_{h}H_{1})(0,jh,x,z)
\]
\begin{equation}
+\frac{h}{2}(p\otimes_{h}A_{0})(0,jh,x,z),\qquad\qquad\qquad\qquad\qquad%
\qquad\qquad\;\;\;  \label{MR-012}
\end{equation}
where $A_{0}(s,jh,v,z)=(L^{2}-2L\widetilde{L}+\widetilde{L}^{2})\widetilde{p}%
(s,jh,v,z).$ The direct calculation shows that
\[
A_{0}(s,jh,v,z)=\frac{1}{4}\sum_{p,q,r,l=1}^{d}(\sigma_{pq}(s,v)-%
\sigma_{pq}(s,z))(\sigma_{rl}(s,v)-\sigma_{rl}(s,z))\qquad\;
\]
\[
\times\frac{\partial^{4}\widetilde{p}(s,jh,v,z)}{\partial
v_{p}\partial
v_{q}\partial v_{r}\partial v_{l}}+\sum_{p,q,r=1}^{d}(\sigma_{pq}(s,v)-%
\sigma_{pq}(s,z))(m_{r}(s,v)-m_{r}(s,z))\qquad\;
\]
\begin{equation}
\times\frac{\partial^{3}\widetilde{p}(s,jh,v,z)}{\partial
v_{p}\partial
v_{q}\partial v_{r}}+\frac{1}{2}\sum_{p,q,r,l=1}^{d}\sigma_{pq}(s,v)\frac{%
\partial\sigma_{rl}(s,v)}{\partial v_{p}}\frac{\partial^{3}\widetilde{p}%
(s,jh,v,z)}{\partial v_{q}\partial v_{r}\partial v_{l}}+(\leq2),
\label{MR-012a}
\end{equation}
where we denote by $(\leq2)$ the sum of terms containing the derivatives of $%
\widetilde{p}(s,jh,v,z)$ of the order less or equal than 2. Note
that for a
constant $C$ $<\infty$ and any $0<\varepsilon<\frac{1}{2}$%
\[
\left|\frac{h}{2}(p\otimes_{h}H_{1})(0,jh,x,z)\right|\leq Ch\phi_{C,\sqrt{jh}%
}\left(z-x\right),\qquad\qquad\;\;\;
\]
\begin{equation}
\left|\frac{h}{2}(p\otimes_{h}A_{0})(0,jh,x,z)\right|\leq
C(\varepsilon)h^{1/2}j^{-(1/2-\varepsilon)}\phi_{C,\sqrt{jh}%
}\left(z-x\right).  \label{MR-012b}
\end{equation}
First inequality (\ref{MR-012b}) follows from (B1) and the well know
estimates for the diffusion density $p$ and for the kernel $H_{1}$ .
The second inequality (\ref{MR-012b}) follows from (B1),
(\ref{MR-012a}) and the following estimate
\begin{eqnarray} \nonumber &&
\frac{h}{2}\sum_{i=0}^{j-1}h\left| \int_{R^{d}}p(0,ih,x,v)\frac{\partial ^{3}%
\widetilde{p}(ih,jh,v,z)}{\partial v_{q}\partial v_{r}\partial v_{l}}%
dv\right| \\ \nonumber &&
\leq \frac{h^{3}}{2}\left| \frac{\partial ^{3}\widetilde{p}%
(0,jh,x,z)}{\partial v_{q}\partial v_{r}\partial v_{l}}\right|
+\frac{h}{2}\sum_{i=1}^{j-1}h\left| \int_{R^{d}}\frac{\partial p(0,ih,x,v)}{%
\partial v_{q}}\frac{\partial ^{2}\widetilde{p}(ih,jh,v,z)}{\partial
v_{r}\partial v_{l}}dv\right| \\ && \leq \frac{h^{3}}{2}\left| \frac{\partial ^{3}%
\widetilde{p}(0,jh,x,z)}{\partial v_{q}\partial v_{r}\partial
v_{l}}\right| +Ch^{1/2}j^{-(1/2-\varepsilon )}B(\frac{1}{2},\varepsilon )\phi _{C,\sqrt{jh}%
}\left( z-x\right) . \label{MR-012c}
\end{eqnarray}
Now we shall estimate the second summand in the right hand side of
(\ref {MR-010}). Clearly
\[
\lambda^{\prime\prime}(s)=\frac{\partial^{2}}{\partial s^{2}}%
p(0,s,x,v)H(s,jh,v,z)+2\frac{\partial}{\partial s}p(0,s,x,v)
\]
\begin{equation}
\times\frac{\partial}{\partial s}H(s,jh,v,z)+p(0,s,x,v)\frac{\partial^{2}}{%
\partial s^{2}}H(s,jh,v,z).\qquad\;\;  \label{MR-013}
\end{equation}
Using forward and backward Kolmogorov equations we get from
(\ref{MR-013}) after long but simple calculations
\[
\sum_{i=0}^{j-1}\int_{ih}^{(i+1)h}(u-ih)^{2}\int_{0}^{1}(1-%
\delta)\int_{R^{d}}\lambda^{\prime\prime}(s)\mid_{s=s_{i}}dvd\delta
du\qquad\qquad\qquad\qquad\qquad\;\;
\]
\begin{equation}
=\sum_{i=0}^{j-1}\int_{ih}^{(i+1)h}(u-ih)^{2}\int_{0}^{1}(1-%
\delta)\sum_{k=1}^{4}\int_{R^{d}}p(0,s,x,v)A_{k}(s,jh,v,z)\mid_{s=s_{i}}dvd%
\delta du,  \label{MR-014}
\end{equation}
where
\[
A_{1}(s,jh,v,z)=(L^{3}-3L^{2}\widetilde{L}+3L\widetilde{L}^{2}-\widetilde{L}%
^{3})\widetilde{p}(s,jh,v,z),\;\;\;\;\,
\]
\[
A_{2}=(L_{1}H+2LH_{1})(s,jh,v,z),\qquad\qquad\qquad\qquad\qquad\,\,
\]
\[
A_{3}(s,jh,v,z)=[(L-\widetilde{L})\widetilde{L}_{1}+2(L_{1}-\widetilde{L}%
_{1})\widetilde{L}]\widetilde{p}(s,jh,v,z),
\]
\begin{equation}
A_{4}(s,jh,v,z)=H_{2}(s,jh,v,z).\qquad\qquad\qquad\qquad\qquad\;\;\;
\label{MR-014a}
\end{equation}
and
\[
H_{l}(s,t,v,z)=(L_{l}-\widetilde{L}_{l})\widetilde{p}(s,t,v,z)\;\;\;\qquad%
\qquad\qquad\qquad\qquad\;
\]
\[
=\frac{1}{2}\sum_{i,j=1}^{d}\left(\frac{\partial^{l}\sigma_{ij}(s,v)}{%
\partial s^{l}}-\frac{\partial^{l}\sigma_{ij}(s,z)}{\partial s^{l}}\right)%
\frac{\partial^{2}\widetilde{p}(s,t,v,z)}{\partial v_{i}\partial v_{j}}%
\qquad\qquad
\]
\begin{equation}
+\sum_{i=1}^{d}\left(\frac{\partial^{l}m_{i}(s,v)}{\partial s^{l}}-\frac{%
\partial^{l}m_{i}(s,z)}{\partial s^{l}}\right)\frac{\partial\widetilde{p}%
(s,t,v,z)}{\partial v_{i}},\;\;\; l=1,2.\;  \label{MR-015}
\end{equation}
Using integration by parts and the definition (\ref{MR-014a}) of $%
A_{2},A_{3} $ and $A_{4}$ it is easy to get that for any
$0<\varepsilon<1/2$
and for $k=2,3,4$%
\[
\left|\sum_{i=0}^{j-1}\int_{ih}^{(i+1)h}(u-ih)^{2}\int_{0}^{1}(1-%
\delta)\int_{R^{d}}p(0,s,x,v)A_{k}(s,jh,v,z)\mid_{s=s_{i}}dvd\delta
du\right|
\]
\begin{equation}
\leq C(\varepsilon)h^{3/2-\varepsilon}\phi_{C,\sqrt{jh}}\left(z-x\right).%
\qquad\qquad\qquad\qquad\qquad\qquad\qquad\qquad\qquad\qquad\;
\label{MR-015a}
\end{equation}
For $k=1$ we shall prove the following estimate for any $0<\varepsilon<\frac{%
1}{2}$%
\[
\left|\sum_{i=0}^{j-1}\int_{ih}^{(i+1)h}(u-ih)^{2}\int_{0}^{1}(1-%
\delta)\int_{R^{d}}p(0,s,x,v)A_{1}(s,jh,v,z)\mid_{s=s_{i}}dvd\delta
du\right|
\]
\begin{equation}
\leq C(\varepsilon)hj^{-(1/2-\varepsilon)}\phi_{C,\sqrt{jh}%
}\left(z-x\right).\qquad\qquad\qquad\qquad\qquad\qquad\qquad\qquad\qquad%
\qquad\;  \label{MR-015b}
\end{equation}
Note that the function $A_{1}(s,jh,v,z)$ can be written as the
following sum
\[
A_{1}(s,jh,v,z)=\frac{1}{8}\sum_{i,j,p,q,l,r=1}^{d}(\sigma_{ij}(s,v)-%
\sigma_{ij}(s,z))(\sigma_{pq}(s,v)-\sigma_{pq}(s,z))(\sigma_{lr}(s,v)\qquad
\]
\[
-\sigma_{lr}(s,z))\frac{\partial^{6}\widetilde{p}(s,jh,v,z)}{\partial
v_{i}\partial v_{j}\partial v_{p}\partial v_{q}\partial v_{l}\partial v_{r}}+%
\frac{3}{4}\sum_{i,j,p,q,l=1}^{d}(\sigma_{ij}(s,v)-\sigma_{ij}(s,z))(%
\sigma_{pq}(s,v)\qquad\;\;\;
\]
\[
-\sigma_{pq}(s,z))(m_{l}(s,v)-m_{l}(s,z))\frac{\partial^{5}\widetilde{p}%
(s,jh,v,z)}{\partial v_{i}\partial v_{j}\partial v_{p}\partial
v_{q}\partial
v_{l}}+\frac{3}{4}\sum_{i,j,p,q,l,r=1}^{d}\sigma_{ij}(s,v)\frac{%
\partial\sigma_{pq}(s,v)}{\partial v_{i}}
\]
\begin{equation}
(\sigma_{lr}(s,v)-\sigma_{lr}(s,z))\frac{\partial^{5}\widetilde{p}%
(s,jh,v,z)}{\partial v_{j}\partial v_{p}\partial v_{q}\partial
v_{l}\partial
v_{r}}+(\leq4),\;\;\;\qquad\qquad\qquad\qquad\qquad\qquad\,
\label{MR-016}
\end{equation}
where we denote by $(\leq4)$ the sum of terms containing the derivatives of $%
\widetilde{p}(s,jh,v,z)$ of the order less or equal than 4. By (B1)
and (\ref {MR-016}) it is clear that the estimate for the left hand
side of (\ref {MR-015a}) for $k=1$ will be the same up to a constant
as for the following sum for fixed $p,q,r,l$

\[
\left|\sum_{i=0}^{j-1}\int_{ih}^{(i+1)h}(u-ih)^{2}\int_{0}^{1}(1-%
\delta)\int_{R^{d}}p(0,s,x,v)\frac{\partial^{4}\widetilde{p}(s,jh,v,z)}{%
\partial v_{p}\partial v_{q}\partial v_{l}\partial v_{r}}\mid_{s=s_{i}}dvd%
\delta du\right|
\]
After integration by parts w.r.t. $v_{p}$ and with the substitution $%
hw=(u-ih)$ in each integral we obtain
\[
\left|\sum_{i=0}^{j-1}\int_{ih}^{(i+1)h}(u-ih)^{2}\int_{0}^{1}(1-%
\delta)\int_{R^{d}}p(0,s,x,v)\frac{\partial^{4}\widetilde{p}(s,jh,v,z)}{%
\partial v_{p}\partial v_{q}\partial v_{l}\partial v_{r}}\mid_{s=s_{i}}dvd%
\delta du\right|\qquad\;\;\;\;\;
\]
\[
=\left|\sum_{i=0}^{j-1}\int_{ih}^{(i+1)h}(u-ih)^{2}%
\int_{0}^{1}(1-\delta)\int_{R^{d}}\frac{\partial p(0,s,x,v)}{\partial v_{p}}%
\frac{\partial^{3}\widetilde{p}(s,jh,v,z)}{\partial v_{q}\partial
v_{l}\partial v_{r}}\mid_{s=s_{i}}dvd\delta du\right|\;\;\;\;
\]
\[
\leq Ch^{2}\phi_{C,\sqrt{jh}}\left(z-x\right)\int_{0}^{1}w^{2}%
\int_{0}^{1}(1-\delta)\sum_{i=0}^{j-1}h\frac{1}{\sqrt{ih+\delta hw}}\frac{1}{%
[(j-i)h-\delta hw)]^{3/2}}d\delta dw\qquad\;\;
\]
\[
\leq Ch^{3/2-\varepsilon}\phi_{C,\sqrt{jh}}\left(z-x\right)\int_{0}^{1}w^{2}%
\int_{0}^{1}(1-\delta)^{1/2-\varepsilon}\sum_{i=0}^{j-1}h\frac{1}{\sqrt{%
ih+\delta hw}}\frac{1}{[(j-\delta w)h-ih)]^{1-\varepsilon}}d\delta
dw
\]
\[
\leq Ch^{3/2-\varepsilon}\phi_{C,\sqrt{jh}}\left(z-x\right)%
\int_{0}^{1}w^{2}dw\int_{0}^{1}(1-\delta)^{1/2-\varepsilon}d\delta%
\int_{0}^{(j-1)h}\frac{dt}{\sqrt{t}[(j-1)h-t]^{1-\varepsilon}}%
\qquad\qquad\;\;\;
\]
\begin{equation}
\leq Chj^{-(1/2-\varepsilon)}B(\frac{1}{2},\varepsilon)\phi_{C,\sqrt{jh}%
}\left(z-x\right),\qquad\qquad\qquad\qquad\qquad\qquad\qquad\qquad\qquad%
\qquad\qquad\;\;\;\;  \label{MR-017}
\end{equation}
where $B(p,q)$ is a Beta function and
$\phi_{C,\rho}\left(z-x\right)$ is defined in Lemma 2. As we
mentioned above (\ref{MR-015b}) follows now from (\ref{MR-017}).
By (B2), (\ref{MR-010}), (\ref{MR-012}), (\ref {MR-012a}), (\ref
{MR-012c}), (\ref{MR-015a}) and (\ref{MR-015b}) we obtain for any $%
0<\varepsilon<\frac{1}{2}$ and $j=1,2,...n$%
\begin{equation}
\left| (p\otimes H-p\otimes _{h}H)(0,T,x,z)\right| \leq
C(\varepsilon )h^{1/2}n^{-(1/2-\varepsilon )}\phi
_{C,\sqrt{T}}(z-x). \label{MR-018}
\end{equation}
We use now the following estimate for
$\Phi(ih,i^{\prime}h,z,z^{\prime})$ that was proved in Konakov and
Mammen (2002) (formula (5.7) on page 284)
\begin{equation}
\left|\Phi(ih,i^{\prime}h,z,z^{\prime})\right|\leq C\frac{1}{\sqrt{%
i^{\prime}h-ih}}\phi_{C,\sqrt{i^{\prime}h-ih}}\left(z^{\prime}-z\right).
\label{MR-019}
\end{equation}
From (B2), (\ref{MR-010}), (\ref{MR-012}), (\ref{MR-017}),
(\ref{MR-018}) and (\ref{MR-019}) we obtain
 the following
representation
\[
(p-p^{d})\left( 0,T,x,y\right) =\frac{h}{2}(p\otimes _{h}H_{1})(0,T,x,y)+%
\frac{h}{2}(p\otimes _{h}A_{0})(0,T,x,y)
\]
\[
+\frac{h}{2}(p\otimes _{h}H_{1}\otimes _{h}\Phi )\left( 0,T,x,y\right) +%
\frac{h}{2}(p\otimes _{h}A_{0}\otimes _{h}\Phi )(0,T,x,y)\qquad \;
\]
\begin{equation}
+R(0,T,x,y),\qquad \qquad \qquad \qquad \qquad \qquad \qquad \qquad
\qquad \qquad \qquad \qquad \;  \label{MR-021}
\end{equation}
where for any $0<\varepsilon <1/2$%
\[
\left| R(0,T,x,y)\right| \leq C(\varepsilon )(h^{3/2-\varepsilon
}+hn^{-(1/2-\varepsilon )})\phi _{C,\sqrt{T}}\left( y-x\right)
\]
\[
=\phi _{C,\sqrt{T}}\left( y-x\right) o(h^{1+\delta}).\qquad \qquad
\qquad \qquad \qquad \qquad \qquad \;\;\,
\]
This representation implies that
\begin{eqnarray} \nonumber
T_{1}&=& \frac{h}{2}[p\otimes
_{h}(L^{2}-2L\widetilde{L}+\widetilde{L}^{2})%
\widetilde{p}\otimes _{h}\Phi )](0,T,x,y) \\ &&
+\frac{h}{2}[p\otimes _{h}(L^{\prime }-\widetilde{L}^{\prime
})\widetilde{p}%
\otimes _{h}\Phi ](0,T,x,y) +R_{T}(0,T,x,y),  \label{eq:28}
\end{eqnarray}
where for any $0<\varepsilon <1/2$%
\begin{equation}
\left| R_{T}(0,T,x,y)\right| \leq C(\varepsilon
)hn^{-1/2+\varepsilon }\phi _{C,\sqrt{T}}(y-x)\leq C(\varepsilon
)h^{1+\delta }\phi _{C,\sqrt{T}}(y-x)  \label{eq:29}
\end{equation}
for $\delta>0$ small enough and where $\Phi
(s,t,x,y)=\sum_{r=0}^{\infty }H^{(r)}(s,t,x,y)$. Here the summand
$H^{(0)}(s,t,x,y)$ is introduced to shorten the notation. By
definition we suppose that $g\otimes
_{h}H^{(0)}(s,t,x,y)=g(s,t,x,y)$ for a
function $%
g$. Note, that in the homogenous case $\sigma _{ij}(s,x)=\sigma
_{ij}(x),m_{i}(s,x)=m_{i}(x)$ and thus the second summand in
(\ref{eq:28}) is equal to $0.$
\bigskip

\noindent  \textit{Asymptotic treatment of the term $T_{2}$. }We
will show that
\begin{eqnarray} \nonumber
&&\left| T_{2}-3\sum_{r=0}^{\infty }\tilde{p}\otimes
_{h}H^{(r)}(0,T,x,y)+\sum_{r=0}^{\infty }\tilde{p}\otimes
_{h}(H+M_{h,1}+%
\sqrt{h}N_{1})^{(r)}(0,T,x,y)\right.
\\ && \qquad \nonumber
\left. +\sum_{r=0}^{\infty }\tilde{p}\otimes
_{h}(H+M_{h,2})^{(r)}(0,T,x,y)+\sum_{r=0}^{\infty
}\tilde{p}\otimes _{h}(H+M_{h,3}^{\prime \prime
})^{(r)}(0,T,x,y)\right|
\\ &&  \leq Chn^{-\delta }\zeta
_{\sqrt{T}}(y-x)  \label{eq:T201}
\end{eqnarray}
with some positive $\delta >0.$ Note that it is enough to consider
the case $r\geq 2$ because for $r=1,2$ the left hand side of
(\ref{eq:T201}) is equal to zero. Note that (\ref{eq:T201})
immediately follows from the following bounds for $r=2,3,...$
\begin{eqnarray} \nonumber &&
\left| \tilde{p}\otimes _{h}(H+M_{h}^{\prime \prime }+\sqrt{h}%
N_{1})^{(r)}(0,T,x,y)\right. \\ \nonumber && \qquad
-\tilde{p}\otimes
_{h}(H+M_{h,1}+M_{h,2}+\sqrt{h}N_{1})^{(r)}(0,T,x,y)
\\ \nonumber && \qquad
\left. -[\tilde{p}\otimes _{h}(H+M_{h,3}^{\prime \prime
})^{(r)}-\tilde{p}%
\otimes _{h}H^{(r)}](0,T,x,y)\right|
\\  &&
\leq C(\varepsilon )h^{3/2-2\varepsilon }\frac{C^{r}}{\Gamma
(\frac{r-1}{2})}%
T^{3\varepsilon +\frac{r-3}{2}}\zeta _{\sqrt{kh}}(v-x),
\label{eq:T202}
\end{eqnarray}
and
\begin{eqnarray}  \nonumber &&
\left| \tilde{p}\otimes
_{h}(H+M_{h,1}+M_{h,2}+\sqrt{h}N_{1})^{(r)}(0,T,x,y)%
\right.
\\ \nonumber && \qquad
-\tilde{p}\otimes _{h}(H+M_{h,1}+\sqrt{h}N_{1})^{(r)}(0,T,x,y)
\\ \nonumber && \qquad
\left. -[\tilde{p}\otimes _{h}(H+M_{h,2})^{(r)}-\tilde{p}\otimes
_{h}H^{(r)}](0,T,x,y)\right| \\  && \leq C(\varepsilon
)h^{3/2-2\varepsilon }\frac{C^{r}}{\Gamma
(\frac{r-1}{2})}%
T^{3\varepsilon +\frac{r-4}{2}}\zeta _{\sqrt{T}}(v-x)
\label{eq:T203}
\end{eqnarray}
for all sufficiently small $\varepsilon >0$ with a constant
$C(\varepsilon)$ that fulfills $\lim_{\varepsilon \longrightarrow
0}C(\varepsilon )=+\infty $. First we prove the bound
(\ref{eq:T202}). Denote the expression under the sign of the
absolute value in (\ref{eq:T202}) by $\Gamma _{r}.$ Note that
$\Gamma _{0}=\Gamma _{1}=0.$ For $r\geq 2$ we make use of the
following recurrence formula \begin{eqnarray}\nonumber \Gamma
_{r}&=&\Gamma _{r-1}\otimes _{h}H+\left[ \tilde{p}\otimes
_{h}(H+M_{h}^{\prime \prime }+\sqrt{h}N_{1})^{(r-1)}\right.
\\ \nonumber &&
\left. -\tilde{p}\otimes
_{h}(H+M_{h,1}+M_{h,2}+\sqrt{h}N_{1})^{(r-1)}\right] \otimes
_{h}(M_{h}^{\prime \prime }+\sqrt{h}N_{1})
\\ \nonumber &&
+\left[ \tilde{p}\otimes
_{h}(H+M_{h,1}+M_{h,2}+\sqrt{h}N_{1})^{(r-1)}-%
\tilde{p}\otimes _{h}(H+M_{h,3}^{\prime \prime })^{(r-1)}\right]
\otimes _{h}M_{h,3}^{\prime \prime }\\&=&I+II+III. \label{eq:T203aa}
\end{eqnarray}
We start with bounding $II$ . First we will give an estimate for
\begin{equation}
\left| \tilde{p}\otimes _{h}(H+M_{h}^{\prime \prime
}+\sqrt{h}N_{1})^{(r-1)}-%
\tilde{p}\otimes
_{h}(H+M_{h,1}+M_{h,2}+\sqrt{h}N_{1})^{(r-1)}\right| .
\label{eq:T203a}
\end{equation}
For $r=2$ we have $\tilde{p}\otimes _{h}M_{h,3}^{\prime \prime
}(0,kh,x,y). $ It follows from (\ref{eq:012h}) that it is enough to
estimate
\begin{equation}
J_{1}=h^{2}\sum_{i=0}^{k-2}h\int \widetilde{p}%
(0,ih,x,v)(f(ih,v)-f(ih,y))D_{v}^{\nu
}\widetilde{p}_{h}((i+1)h,kh,v,y)dv \label{eq:T204}
\end{equation}
for $\left| \nu \right| =4$ and
\[
J_{2}=h^{3/2}\sum_{i=0}^{k-2}h\int
\widetilde{p}(0,ih,x,v)(f(ih,v)-f(ih,y))%
\int q(ih,v,\theta )\theta ^{\nu }\int_{0}^{1}(1-\delta )^{2}
\]
\begin{equation}
\times D_{v}^{\nu
+e_{l}+e_{q}}\widetilde{p}_{h}((i+1)h,kh,v+\delta
\widetilde{h}(\theta ),y)d\delta d\theta dv  \label{eq:T205}
\end{equation}
for $\left| \nu \right| =3.$ Here $f(t,x)$ is a function with
$D_{x}^{\nu }f(t,x),\left| \nu \right| =0,1,2,3$ bounded uniformly
in $(t,x).$ An estimate for $J_{1}$ follows from (\ref{eq:15}). This
gives
\begin{equation}
\left| J_{1}\right| \leq Ch^{2-\varepsilon }(kh)^{\varepsilon
-1/2}B(\frac{1%
}{2},\varepsilon )\zeta _{\sqrt{kh}}^{S}(y-x).  \label{eq:T206}
\end{equation}
The estimate for $J_{2}$ can be obtained analogously to the estimate
of $I_{2}$ (see (\ref{eq:17})). By integrating by parts we get
\begin{eqnarray*}
J_{2}&=&h^{3/2}\sum_{i=0}^{k-2}h\int_{0}^{1}(1-\delta )^{2}d\delta
\int d\theta \cdot \theta ^{\nu }\int
D_{v}^{e_{l}+e_{q}}[\widetilde{p}(0,ih,x,v)
\\
&&\times (f(ih,v)-f(ih,y))q(ih,v,\theta )]D_{v}^{\nu }\widetilde{p}%
_{h}((i+1)h,kh,v+\delta \widetilde{h}(\theta ),y)dv.
\end{eqnarray*}
The derivative
\[
D_{v}^{e_{l}+e_{q}}[\widetilde{p}(0,ih,x,v)(f(ih,v)-f(ih,y))q(ih,v,\theta
)]
\]
is a sum of 9 summands. By using integration by parts once more for
all summands which contain $D_{\nu }^{\mu }\widetilde{p}(0,ih,x,v)$
with $\left| \mu \right| <2$ , we obtain
\begin{eqnarray}\nonumber
\left| J_{2}\right| &\leq& Ch^{3/2-2\varepsilon }\zeta _{\sqrt{kh}%
}^{S-3}(y-x)\int \psi (\theta )\left\| \theta \right\|
^{3}(h^{(S-3)/2}\left\| \theta \right\| ^{S-3}+1)d\theta
\\ &&
\sum_{i=1}^{k-2}h\frac{1}{(ih)^{1-\varepsilon }}\times \frac{1}{%
(kh-ih)^{1-\varepsilon }}\leq Ch^{3/2-2\varepsilon }B(\varepsilon
,\varepsilon )(kh)^{2\varepsilon -1}\zeta _{\sqrt{kh}}^{S-3}(y-x)
\label{eq:T207}
\end{eqnarray}
for any $\varepsilon \in (0,1/4).$ It follows from (\ref{eq:T206})
and (\ref
{eq:T207}) that for $r=2$ (\ref{eq:T203a}) does not exceed $%
Ch^{3/2-2\varepsilon }B(\varepsilon ,\varepsilon )(kh)^{\varepsilon
-1}\zeta _{\sqrt{kh}}^{S-3}(y-x).$ For $r\geq 3$ we use the
recurrence relation
\begin{eqnarray} \nonumber
&&\tilde{p}\otimes _{h}(H+M_{h}^{\prime \prime
}+\sqrt{h}N_{1})^{(r-1)}-\tilde{%
p}\otimes _{h}(H+M_{h,1}+M_{h,2}+\sqrt{h}N_{1})^{(r-1)}
\\ \nonumber &&
=\left[ \tilde{p}\otimes _{h}(H+M_{h}^{\prime \prime }+\sqrt{h}%
N_{1})^{(r-2)}-\tilde{p}\otimes
_{h}(H+M_{h,1}+M_{h,2}+\sqrt{h}N_{1})^{(r-2)}%
\right]
\\ \nonumber && \qquad
\otimes _{h}(H+M_{h}^{\prime \prime
}+\sqrt{h}N_{1})+[\tilde{p}\otimes
_{h}(H+M_{h,1}+M_{h,2}+\sqrt{h}N_{1})^{(r-2)}]\otimes
_{h}M_{h,3}^{\prime \prime }
\\ &&
=I^{\prime }+II^{\prime }.  \label{eq:T208}
\end{eqnarray}
From (\ref{eq:T208}) we obtain for $r=3$
\begin{eqnarray*}
\left| I^{\prime }\right| &\leq& Ch^{3/2-2\varepsilon }B(\varepsilon
,\varepsilon )\zeta _{\sqrt{kh}}^{S-3}(y-x)\sum_{i=1}^{k-2}h(ih)^{%
\varepsilon -1}(kh-ih)^{-1/2}
\\
&\leq& Ch^{3/2-2\varepsilon }B(\varepsilon ,\varepsilon )B(\frac{1}{2}%
,\varepsilon )(kh)^{\varepsilon -1/2}\zeta
_{\sqrt{kh}}^{S-3}(y-x).
\end{eqnarray*}
To estimate $II^{\prime }$ we use the following estimates
\begin{eqnarray}
&&\left|
D_{v}^{a}D_{x}^{b}(H+M_{h,1}+M_{h,2}+\sqrt{h}N_{1})(jh,kh,x,v)\right|
\leq C\rho ^{-1-\left| a\right| -\left| b\right| }\zeta _{\rho
}(v-x), \label{eq:T209}
\\ &&
\left|
D_{x}^{b}(H+M_{h,1}+M_{h,2}+\sqrt{h}N_{1})(jh,kh,x,x+v)\right| \leq
C\rho ^{-1}\zeta _{\rho }(v-x) . \label{eq:T210}
\end{eqnarray}
To prove (\ref{eq:T209}) one can use the following estimates for the
summands in $M_{h,1},M_{h,2}$ and $\sqrt{h}N_{1}$
\begin{eqnarray*}
&& \left |h^{1/2}D_{v}^{a}D_{x}^{b}\left[ D_{x}^{\nu }\widetilde{p}%
_{h}((j+1)h,kh,x,v)(f(jh,x)-f(jh,v))\right]\right | \leq C\rho
^{-1-\left| a\right| -\left| b\right| }\zeta _{\rho }(v-x),\left|
\nu \right| =3, \label{eq:T211}
\\ && \left |
hD_{v}^{a}D_{x}^{b}\left[ D_{x}^{\nu }\widetilde{p}%
_{h}((j+1)h,kh,x,v)(f(jh,x)-f(jh,v))\right]\right | \leq C\rho
^{-1-\left| a\right| -\left| b\right| }\zeta _{\rho }(v-x),\left|
\nu \right| =4, \label{eq:T212}
\\ &&
\left | h^{1/2}D_{v}^{a}D_{x}^{b}\left[ D_{x}^{\nu +e_{p}+e_{q}}\widetilde{p}%
_{h}((j+1)h,kh,x,v)\rho ^{2}(f(jh,x)-f(jh,v))\right]\right | \leq
C\rho ^{-1-\left| a\right| -\left| b\right| }\zeta _{\rho
}(v-x),\left| \nu \right| =3, \label{eq:T213}
\end{eqnarray*}
for a function $f(t,x)$ with $\left| a\right| +\left| b\right| $
derivatives w.r.t.\ $x$ that are uniformly bounded w.r.t.\ $t$.
These estimates are direct consequences of Lemma 7. To prove
(\ref{eq:T210}) one can use the following estimates for the summands
in $M_{h,1},M_{h,2}$
and $\sqrt{h}N_{1}$%
\begin{eqnarray*}
&& \left |h^{1/2}D_{x}^{b}\left[ D_{x}^{\nu
}\widetilde{p}_{h}((j+1)h,kh,x,y)\mid
_{y=x+v}(f(jh,x)-f(jh,x+v))\right] \right |\leq C\rho ^{-1}\zeta
_{\rho }(v-x),\left| \nu \right| =3,  \label{eq:T214}
\\ && \left |
hD_{x}^{b}\left[ D_{x}^{\nu }\widetilde{p}_{h}((j+1)h,kh,x,y)\mid
_{y=x+v}(f(jh,x)-f(jh,x+v))\right] \right |\leq C\rho ^{-1}\zeta
_{\rho }(v-x),\left| \nu \right| =4,  \label{eq:T215}
\\ && \left |
h^{1/2}D_{x}^{b}\left[ D_{x}^{\nu +e_{p}+e_{q}}\widetilde{p}%
_{h}((j+1)h,kh,x,v)\mid _{y=x+v}\rho
^{2}(f(jh,x)-f(jh,x+v))\right]\right | \leq C\rho ^{-1}\zeta
_{\rho }(v-x),\left| \nu \right| =3. \label{eq:T216}
\end{eqnarray*}
Again, these estimates follow from the estimates obtained in the
proof of Lemma 7. Note that with $z(V_{j,k}^{-1/2}(y)\mu
_{j,k}(y),x,y)=V_{j,k}^{-1/2}(y)(y-x-\mu _{j,k}(y))$ it holds
$$
\left| \frac{\partial z}{\partial y}\right| \leq \frac{C}{\rho
},\left| \frac{\partial z}{\partial x}\right| \leq \frac{C}{\rho }
$$
and that with $z(V_{j,k}^{-1/2}(x+v),\mu
_{j,k}(x+v),x,x+v)=V_{j,k}^{-1/2}(x+v)(v-\mu _{j,k}(x+v))$ it holds
$$
\left| \frac{\partial z}{\partial x}\right| \leq \left|
\frac{\partial
z}{%
\partial V_{j,k}^{-1/2}}\right| \left| \frac{\partial V_{j,k}^{-1/2}}{%
\partial x}\right| +\left| \frac{\partial z}{\partial \mu
_{j,k}}\right| \left| \frac{\partial \mu _{j,k}}{\partial x}\right|
\leq C(\left\| v\right\| +1). $$ Now with the inequalities
(\ref{eq:T209}),(\ref{eq:T210}) we can proceed like in the proof of
Theorem 2.3 in Konakov and Mammen (2002). This gives the following
estimate for $r=3,4,... $
\begin{eqnarray} \nonumber
&& \nonumber \left| D_{v}^{a}D_{x}^{b}[\widetilde{p}\otimes
_{h}(H+M_{h,1}+M_{h,2}+\sqrt{h%
}N_{1})^{(r-2)}](jh,kh,x,v)\right| \\ && \leq
C^{r}B(1,\frac{1}{2})\times ...\times
B(\frac{r-1}{2},\frac{1}{2})\rho ^{r-2-\left| a\right| -\left|
b\right| }\zeta _{\rho }(v-x).\label{eq:T217}
\end{eqnarray}
Now we denote $\widetilde{p}_{1,r}=\widetilde{p}\otimes
_{h}(H+M_{h,1}+M_{h,2}+\sqrt{h}N_{1})^{(r)},\widetilde{p}_{0}=\widetilde{p}$.
To estimate $\widetilde{p}_{1,r-2}\otimes _{h}M_{h,3}^{\prime \prime
}$ it is enough to make the same calculations with integration by
parts as it was done above for $J_{1}$ and $J_{2}.$ This gives
\begin{eqnarray}\label{eq:T217a}
\left| II^{\prime }\right| &\leq& \left|
\widetilde{p}_{1,r-2}\otimes _{h}M_{h,3}^{\prime \prime
}(0,kh,x,y)\right|
\\
&\leq& C^{r}h^{3/2-2\varepsilon
}B(1,\frac{1}{2})B(1+\frac{1}{2},\frac{1}{2}%
)\times ...\times
B(1+\frac{r-3}{2},\frac{1}{2})B(\frac{r-2}{2},\varepsilon
)(kh)^{\varepsilon +\frac{r-4}{2}}\zeta _{\sqrt{kh}}(v-x) \nonumber
\end{eqnarray}
and by induction
\begin{equation}
\left| I^{\prime }\right| \leq C^{r}h^{3/2-2\varepsilon
}B(\varepsilon
,%
\frac{1}{2})B(\varepsilon +\frac{1}{2},\frac{1}{2})\times
....\times B(\varepsilon +\frac{r-3}{2},\frac{1}{2})B(\varepsilon
,\varepsilon )(kh)^{\varepsilon +\frac{r-4}{2}}\zeta
_{\sqrt{kh}}(v-x), \label{eq:T217b}
\end{equation}
$r=3,4,...$. Comparing (\ref{eq:T217a}) and (\ref{eq:T217b}) we
obtain \
that for $r\geq 3$%
\begin{eqnarray} \nonumber
&&\left| \tilde{p}\otimes _{h}(H+M_{h}^{\prime \prime
}+\sqrt{h}N_{1})^{(r-1)}-%
\tilde{p}\otimes
_{h}(H+M_{h,1}+M_{h,2}+\sqrt{h}N_{1})^{(r-1)}\right|
\\ &&
\leq C^{r}h^{3/2-2\varepsilon }B(\varepsilon
,\frac{1}{2})B(\varepsilon
+\frac{1}{%
2},\frac{1}{2})\times ...\times B(\varepsilon
+\frac{r-3}{2},\frac{1}{2}%
)B(\varepsilon ,\varepsilon )(kh)^{\varepsilon
+\frac{r-4}{2}}\zeta
_{\sqrt{%
kh}}(v-x).  \label{eq:T217c}
\end{eqnarray}
From (\ref{eq:T217c}) we get the following estimate for $II$%
\begin{equation}
\left| II\right| \leq C^{r}h^{3/2-2\varepsilon }B(\varepsilon
,\varepsilon )B(\varepsilon ,\frac{1}{2})B(\varepsilon
+\frac{1}{2},\frac{1}{2})\times ...\times B(\varepsilon
+\frac{r-2}{2},\frac{1}{2})T^{\varepsilon
+\frac{r-3%
}{2}}\zeta _{\sqrt{T}}(v-x).  \label{eq:T218}
\end{equation}
To estimate $III$ note that the following inequalities  that are
similar to
(\ref{eq:T209}), (%
\ref{eq:T210}),(\ref{eq:T217}) hold for $H+M_{h,3}^{\prime \prime
},$
\begin{eqnarray}
\nonumber && \left| D_{v}^{a}D_{x}^{b}(H+M_{h,3}^{\prime \prime
})(jh,kh,x,v)\right| \leq C\rho ^{-1-\left| a\right| -\left|
b\right| }\zeta _{\rho }(v-x) , \\ \nonumber && \left|
D_{x}^{b}(H+M_{h,3}^{\prime \prime })(jh,kh,x,x+v)\right| \leq C\rho
^{-1}\zeta _{\rho }(v-x),\\ \nonumber &&  \left|
D_{v}^{a}D_{x}^{b}[\widetilde{p}\otimes _{h}(H+M_{h,3}^{\prime
\prime })^{(r)}](jh,kh,x,v)\right|
\\  && \qquad
\leq C^{r}B(1,\frac{1}{2})\times ...\times
B(\frac{r+1}{2},\frac{1}{2})\rho ^{r-\left| a\right| -\left|
b\right| }\zeta _{\rho }(v-x). \label{eq:T221}
\end{eqnarray}
To prove the last three inequalities it is enough to get the
corresponding estimates for summands in $M_{h,3}^{\prime \prime }$
$\
$(see (%
\ref{eq:012h}))$.$ These estimates can be proved by the same
arguments as used in the proofs of (\ref{eq:T209}), (\ref{eq:T210}),
and (\ref{eq:T217}). To estimate $III$ we have now to estimate
$\widetilde{p}_{1,r}\otimes _{h}M_{h,3}^{\prime \prime }$ and
$\widetilde{p}_{2,r}\otimes _{h}M_{h,3}^{\prime \prime }$ where
\[
\widetilde{p}_{2,r}=\widetilde{p}\otimes _{h}(H+M_{h,3}^{\prime
\prime })^{(r)}.
\]
Using integration by parts and inequality (\ref{eq:T221}), we obtain
for $\widetilde{p}_{2,r-1}\otimes _{h}M_{h,3}^{\prime \prime }$ the
same estimate as for $\widetilde{p}_{1,r}\otimes _{h}M_{h,3}^{\prime
\prime
}$%
\begin{eqnarray*} &&
\left| \widetilde{p}_{2,r}\otimes _{h}M_{h,3}^{\prime \prime
}(0,kh,x,y)\right|
\\ &&
\leq C^{r}h^{3/2-2\varepsilon }B(1,\frac{1}{2})\times ...\times
B(\frac{r+1}{%
2},\frac{1}{2})B(\frac{r}{2},\varepsilon )(kh)^{\varepsilon
+\frac{r-2}{2}%
}\zeta _{\sqrt{kh}}(y-x)
\end{eqnarray*}
for $i=1,2$. Hence for $r=2,3,...$
\begin{eqnarray} \nonumber
\left| III\right| &\leq& \left| \widetilde{p}_{1,r-1}\otimes
_{h}M_{h,3}^{\prime \prime }(0,kh,x,y)\right| +\left| \widetilde{p}%
_{2,r-1}\otimes _{h}M_{h,3}^{\prime \prime }(0,kh,x,y)\right|
\\
&\leq& C^{r}h^{3/2-2\varepsilon }B(1,\frac{1}{2})\times ...\times
B(\frac{r}{2}%
,\frac{1}{2})B(\frac{r-1}{2},\varepsilon )(kh)^{\varepsilon
+\frac{r-3}{2}%
}\zeta _{\sqrt{kh}}(y-x).  \label{eq:T223}
\end{eqnarray}
From (\ref{eq:T203aa}), (\ref{eq:T218}) and (\ref{eq:T223}) we get
for $r=2,3,...$
\[
\left| \Gamma _{r}(0,kh,x,y)\right| \leq C^{r}h^{3/2-2\varepsilon
}B(\varepsilon ,\varepsilon )B(\varepsilon ,\frac{1}{2})\times
....\times B(\varepsilon
+\frac{r-2}{2},\frac{1}{2})(kh)^{\varepsilon
+\frac{r-3}{2}%
}\zeta _{\sqrt{kh}}(v-x).
\]
In particular,
\begin{equation}
\left| \Gamma _{r}(0,T,x,y)\right| \leq h^{3/2-2\varepsilon
}\frac{\Gamma ^{3}(\varepsilon )}{\Gamma (2\varepsilon
)}\frac{C^{r}}{\Gamma
(\varepsilon +%
\frac{r-1}{2})}T^{\varepsilon +\frac{r-3}{2}}\zeta
_{\sqrt{kh}}(v-x),r=2,3... \label{eq:T224}
\end{equation}
for any $\varepsilon \in (0,1/4).$ Now we estimate the left hand
side of (\ref{eq:T203}). Denote the expression under the sign of the
absolute value in (\ref{eq:T203}) by $\digamma _{r}.$ Note that
$\digamma _{0}=\digamma _{1}=0$. For $r\geq 2$ we make use of the
following recurrence formula
\begin{eqnarray*}
\digamma _{r}&=&\digamma _{r-1}\otimes _{h}H+\left[
\widetilde{p}\otimes
_{h}(H+M_{h,1}+M_{h,2}+\sqrt{h}N_{1})^{(r-1)}\right.
\\ && \qquad
\left. -\widetilde{p}\otimes
_{h}(H+M_{h,1}+\sqrt{h}N_{1})^{(r-1)}\right] \otimes
_{h}(M_{h,1}+M_{h,2}+\sqrt{h}N_{1})
\\ && \qquad
+\left[ \widetilde{p}\otimes _{h}(H+M_{h,1}+\sqrt{h}N_{1})^{(r-1)}-%
\widetilde{p}\otimes _{h}(H+M_{h,1})^{(r-1)}\right] \otimes
_{h}M_{h,2}
\\ &=& I+II+III.
\end{eqnarray*}
We start again from the estimation of
\[
A_{r-1}=\widetilde{p}\otimes
_{h}(H+M_{h,1}+M_{h,2}+\sqrt{h}N_{1})^{(r-1)}-%
\widetilde{p}\otimes _{h}(H+M_{h,1}+\sqrt{h}N_{1})^{(r-1)}.
\]
For $r=2$ we have $A_{1}=(\widetilde{p}\otimes
_{h}M_{h,2})(0,kh,x,y).$ It is enough to estimate
\[
J_{3}=h\sum_{i=0}^{k-2}h\int
\widetilde{p}(0,ih,x,v)(f(ih,v)-f(ih,y))D_{v}^{%
\nu }\widetilde{p}_{h}((i+1)h,kh,v,y)dv
\]
for $\left| \nu \right| =4.$ Analogously to (\ref{eq:15}) we obtain
that
\[
\left| J_{3}\right| \leq Ch^{1-\varepsilon }(kh)^{-1/2+\varepsilon
}B(\frac{1%
}{2},\varepsilon )\zeta _{\sqrt{kh}}^{S}(y-x).
\]
and, hence,
\begin{equation}
\left| A_{1}\right| \leq Ch^{1-\varepsilon }(kh)^{-1/2+\varepsilon
}B(\frac{1%
}{2},\varepsilon )\zeta _{\sqrt{kh}}^{S}(y-x).  \label{eq:T224a}
\end{equation}
For $r\geq 3$ we use the recurrence relation
\begin{eqnarray} \nonumber
A_{r-1}&=&A_{r-2}\otimes _{h}(H+M_{h,1}+M_{h,2}+\sqrt{h}N_{1})
\\ \nonumber && \qquad +\left[ \widetilde{p}\otimes
_{h}(H+M_{h,1}+\sqrt{h}N_{1})^{(r-2)}\right] \otimes _{h}M_{h,2}\\
&=& I^{\prime }+II^{\prime }.  \label{eq:T225}
\end{eqnarray}
From (\ref{eq:T224a}) and (\ref{eq:T225}) we obtain for $r=3$%
\begin{eqnarray} \nonumber
\left| I^{\prime }\right| &\leq& Ch^{1-\varepsilon
}B(\frac{1}{2},\varepsilon )\zeta
_{\sqrt{kh}}^{S}(y-x)\sum_{i=0}^{k-2}h(ih)^{\varepsilon
-1/2}(kh-ih)^{-1/2}
\\
&\leq& Ch^{1-\varepsilon }B(\frac{1}{2},\varepsilon
)B(\frac{1}{2},\varepsilon +\frac{1}{2})(kh)^{\varepsilon }\zeta
_{\sqrt{kh}}^{S}(y-x). \label{eq:T225a}
\end{eqnarray}
To estimate $II^{\prime }$ we use the following inequality for
$r=3,4,... $
\begin{eqnarray} \nonumber &&
\left| D_{v}^{a}D_{x}^{b}[\widetilde{p}\otimes _{h}(H+M_{h,1}+\sqrt{h}%
N_{1})^{(r-2)}](jh,kh,x,v)\right|
\\ &&
\leq C^{r}B(1,\frac{1}{2})\times ...\times
B(\frac{r-1}{2},\frac{1}{2})\rho ^{r-2-\left| a\right| -\left|
b\right| }\zeta _{\rho }(v-x) .\label{eq:T226}
\end{eqnarray}
This inequality follows from (\ref{eq:T217}). We have
\begin{equation}
\left| II^{\prime }\right| \leq Ch^{1-\varepsilon }B(1,\varepsilon
)(kh)^{\varepsilon }\zeta _{\sqrt{kh}}^{S}(y-x).  \label{eq:T226a}
\end{equation}
Comparing (\ref{eq:T225a}) and (\ref{eq:T226a}) we obtain that
$\left| A_{2}\right| \leq $$C^{2}h^{1-\varepsilon
}B(\frac{1}{2},\varepsilon
)B(%
\frac{1}{2},\varepsilon +\frac{1}{2})(kh)^{\varepsilon }\zeta
_{\sqrt{kh}%
}^{S}(y-x).$ By induction we easily get for $r=2,3...$
\begin{equation}
\left| A_{r-1}(0,kh,x,y)\right| \leq C^{r}h^{1-\varepsilon
}B(\frac{1}{2}%
,\varepsilon )\times ...\times B(\frac{1}{2},\varepsilon
+\frac{r-2}{2}%
)(kh)^{\varepsilon +\frac{r-3}{2}}\zeta _{\sqrt{kh}}^{S}(y-x).
\label{eq:T227}
\end{equation}
To estimate
\[
A_{r-1}\otimes _{h}(M_{h,1}+M_{h,2}+\sqrt{h}N_{1})
\]
it is enough to estimate
$$
J_{4}=h^{1/2}\sum_{i=0}^{k-2}h\int
A_{r-1}(0,ih,x,v)(f(ih,v)-f(ih,y))D_{v}^{\nu }\widetilde{p}%
_{h}((i+1)h,kh,v,y)dv $$
for $\left| \nu \right| =3$,%
$$
J_{5}=h\sum_{i=0}^{k-2}h\int
A_{r-1}(0,ih,x,v)(f(ih,v)-f(ih,y))D_{v}^{\nu }%
\widetilde{p}_{h}((i+1)h,kh,v,y)dv $$ for $\left| \nu \right| =4$,
and
$$
J_{6}=h^{1/2}\sum_{i=0}^{k-2}h\int
A_{r-1}(0,ih,x,v)(f(ih,v)-f(ih,y))(kh-ih)D_{v}^{\nu
+e_{p}+e_{q}}\widetilde{p%
}(ih,kh,v,y)dv
$$
for $\left| \nu \right| =3.$ It follows from (\ref{eq:T227}) that
\[
\left| J_{4}\right| \leq C^{r}h^{3/2-2\varepsilon
}B(\frac{1}{2},\varepsilon )\times ...\times
B(\frac{1}{2},\varepsilon +\frac{r-2}{2})B(\varepsilon
,\varepsilon +\frac{r-1}{2})(kh)^{2\varepsilon
+\frac{r-3}{2}}\zeta
_{\sqrt{%
kh}}^{S}(y-x).
\]
Clearly, the same estimate holds for $J_{5}$ and $J_{6}$ . Thus we
obtain
$$
\left| II\right| \leq C^{r}h^{3/2-2\varepsilon
}B(\frac{1}{2},\varepsilon )\times ...\times
B(\frac{1}{2},\varepsilon +\frac{r-2}{2})B(\varepsilon
,\varepsilon +\frac{r-1}{2})(kh)^{2\varepsilon
+\frac{r-3}{2}}\zeta
_{\sqrt{%
kh}}^{S}(y-x). $$ Now we give an estimate for $III$. We write
\[
B_{r-1}=\widetilde{p}\otimes _{h}(H+M_{h,1}+\sqrt{h}N_{1})^{(r-1)}-%
\widetilde{p}\otimes _{h}(H+M_{h,1})^{(r-1)}.
\]
Using the recurrence equation
\[
B_{r-1}=B_{r-2}\otimes
_{h}(H+M_{h,1}+\sqrt{h}N_{1})+\widetilde{p}\otimes
_{h}(H+M_{h,1})^{(r-2)}\otimes _{h}\sqrt{h}N_{1},B_{0}=0
\]
we obtain that
\begin{equation}
III=\sum_{l=0}^{r-2}\widetilde{p}_{3,l}\otimes
_{h}\sqrt{h}N_{1}\otimes
_{h}(H+M_{h,1}+\sqrt{h}N_{1})^{(r-l-2)}\otimes
_{h}M_{h,2}(0,T,x,y), \label{eq:T232}
\end{equation}
where $\widetilde{p}_{3,l}=\widetilde{p}\otimes
_{h}(H+M_{h,1})^{(l)}.$ To estimate $III$ it is enough to estimate a
typical term in the last sum. Thus we have to estimate
\begin {eqnarray*}
&&h^{3/2}\sum_{k=0}^{n-2}h\int \Big\{ \sum_{j=0}^{k-1}h\Big[ \int
\sum_{i=0}^{j-1}h\int
\widetilde{p}_{3,l}(0,ih,x,w)(jh-ih)D_{w}^{\mu
+e_{n}+e_{m}}\widetilde{p}(ih,jh,w,z)
\\ && \qquad
 \times (g(ih,w)-g(ih,z))dw\Big] (H+M_{h,1}+\sqrt{h}%
N_{1})^{(r-l-2)}(ih,kh,z,v)dz\Big\} (f(kh,v)-f(kh,y))
\\ && \qquad
\times D_{v}^{\nu }\widetilde{p}_{h}((k+1)h,T,v,y)dv.
\end{eqnarray*}
To estimate this term we apply two times an integration by parts in
the internal integral $\int ...dw$ and then we make two times an
integration by parts in $\int ...dv.$ We also use the following
estimates
\begin{eqnarray*}
&&\left| D_{w}^{a}D_{x}^{b}\widetilde{p}_{3,l}(0,ih,x,w)\right|
\\ && \qquad
\leq C^{l}B(1,\frac{1}{2})\times ...\times
B(\frac{l+1}{2},\frac{1}{2})(ih)^{%
\frac{l-\left| a\right| -\left| b\right| }{2}}\zeta
_{\sqrt{ih}}(w-x), \\ && \left|
D_{v}^{a}D_{z}^{b}(H+M_{h,1}+\sqrt{h}N_{1})^{(r-l-2)}(ih,kh,z,v)%
\right|
\\ && \qquad
\leq C^{r-l-2}B(\frac{1}{2},\frac{1}{2})\times ...\times
B(\frac{1}{2},\frac{%
r-l-3}{2})(kh-ih)^{\frac{r-l-4-\left| a\right| -\left| b\right|
}{2}}\zeta _{%
\sqrt{kh-ih}}(v-z)
\end{eqnarray*}
for $0\leq l\leq r-3$ with  $B(\frac{1}{2},0)=1$. This gives the
following estimate for
any $0\leq l\leq r-3,r\geq 2$%
\begin{eqnarray} \nonumber
&&\left| \widetilde{p}_{3,l}\otimes _{h}\sqrt{h}N_{1}\otimes
_{h}(H+M_{h,1}+%
\sqrt{h}N_{1})^{(r-l-2)}\otimes _{h}M_{h,2}(0,T,x,y)\right|
\\ && \qquad
\leq C^{r}h^{3/2-3\varepsilon }\frac{\Gamma (\varepsilon )}{\Gamma
(3\varepsilon +\frac{r-1}{2})}T^{3\varepsilon +\frac{r-3}{2}}\zeta
_{\sqrt{T}%
}(y-x).  \label{eq:T236}
\end{eqnarray}
For $l=r-2$ we have to estimate
\[
\widetilde{p}_{3,r-2}\otimes _{h}\sqrt{h}N_{1}{}_{h}\otimes
_{h}M_{h,2}(0,T,x,y).
\]
This is a finite sum of terms corresponding to the different
summands
in $%
N_{1,h}$ and $M_{h,2}$. To estimate a typical term
\begin{eqnarray*} &&
h^{3/2}\sum_{k=0}^{n-2}h\int \Big\{ \sum_{j=0}^{k-1}h\int
\widetilde{p}%
_{3,r-2}(0,jh,x,w)(kh-jh)D_{w}^{\mu +e_{n}+e_{m}}\widetilde{p}%
(jh,kh,w,v)
\\ && \qquad
 \times (g(jh,w)-g(jh,v))dw\Big\}
(f(kh,v)-f(kh,y))D_{v}^{\nu
}%
\widetilde{p}_{h}((k+1)h,T,v,y)dv
\end{eqnarray*}
again we apply integration by parts and after direct calculations we
obtain the following estimate for $r=2,3,... $
\begin{eqnarray}&& \nonumber
\left| \widetilde{p}_{3,r-2}\otimes _{h}\sqrt{h}N_{1}{}_{h}\otimes
_{h}M_{h,2}(0,T,x,y)\right|
\\ && \qquad
\leq C^{r}h^{3/2-3\varepsilon }\frac{\Gamma (\varepsilon
+\frac{r-2}{2})}{%
\Gamma (3\varepsilon +\frac{r-2}{2})}\Gamma ^{2}(\varepsilon
)\frac{1}{%
\Gamma (\frac{r}{2})}T^{3\varepsilon +\frac{r-4}{2}}\zeta _{\sqrt{T}%
}(y-x). \label{eq:T238}
\end{eqnarray}
The inequalities (\ref{eq:T202}) and (\ref{eq:T203}) follow now from
(\ref {eq:T224}), (\ref{eq:T236}) and (\ref{eq:T238}).
\bigskip

\noindent \textit{Asymptotic treatment of the term $T_{3}$. }We will
show that
\[
\left| T_{3}-\left[ \sum_{r=0}^{\infty }\tilde{p}\otimes
_{h}(H+A)^{(r)}(0,T,x,y)-\sum_{r=0}^{\infty }\tilde{p}\otimes
_{h}H^{(r)}(0,T,x,y)\right] \right|
\]
\begin{equation}
\leq Chn^{-\delta }\zeta _{\sqrt{T}}(y-x),  \label{T301}
\end{equation}
where $A=M_{h}^{\prime \prime }-M_{h}=-\frac{h}{2}(L_{\star }^{2}-2L%
\widetilde{L}+\widetilde{L}^{2})\lambda (x)$. Write
\begin{eqnarray*}
C_{r}&=&\tilde{p}\otimes _{h}(H+M_{h}^{\prime \prime }+\sqrt{h}%
N_{1})^{(r)}(0,T,x,y)
\\ &&
-\tilde{p}\otimes _{h}(H+M_{h}+\sqrt{h}N_{1})^{(r)}(0,T,x,y)
\\ &&
-[\tilde{p}\otimes _{h}(H+A)^{(r)}-\tilde{p}\otimes
_{h}H^{(r)}](0,T,x,y).
\end{eqnarray*}
Similarly as in (\ref{eq:T203aa}) we have the following recurrence
relation
\begin{eqnarray} \nonumber
C_{r}&=&C_{r-1}\otimes _{h}H+\left[ \tilde{p}\otimes
_{h}(H+M_{h}^{\prime \prime }+\sqrt{h}N_{1})^{(r-1)}\right.
\\ && \nonumber
\left. -\tilde{p}\otimes
_{h}(H+M_{h}+\sqrt{h}N_{1})^{(r-1)}\right] \otimes
_{h}(M_{h}^{\prime \prime }+\sqrt{h}N_{1})
\\ && \nonumber
+\left[ \tilde{p}\otimes _{h}(H+M_{h}+\sqrt{h}N_{1})^{(r-1)}-\tilde{p}%
\otimes _{h}(H+A)^{(r-1)}\right] \otimes _{h}A \\ &=&I+II+III.
\label{T302}
\end{eqnarray}
With the notation
\[
D_{r-1}=\tilde{p}\otimes _{h}(H+M_{h}+\sqrt{h}N_{1})^{(r-1)}-\tilde{p}%
\otimes _{h}(H+A)^{(r-1)}
\]
we get
\[
D_{r-1}=D_{r-2}\otimes
_{h}(H+M_{h}+\sqrt{h}N_{1})+\tilde{p}_{h}\otimes
_{h}(H+A)^{(r-2)}\otimes _{h}(M_{h}-A+\sqrt{h}N_{1}).
\]
Iterative application gives
\begin{eqnarray*} \nonumber
III&=&D_{r-1}\otimes _{h}A \\&=&
\sum_{l=0}^{r-2}\widetilde{p}_{4,l}\otimes _{h}(M_{h}-A+\sqrt{h}%
N_{1})\otimes _{h}(H+M_{h}+\sqrt{h}N_{1})^{(r-l-2)}\otimes
_{h}A(0,T,x,y),
\end{eqnarray*}
where $\widetilde{p}_{4,l}=\widetilde{p}\otimes _{h}(H+A)^{(l)}.$
This sum can be estimated  in exactly the same way as the sum in
(\ref {eq:T232}). This gives for $r=2,3,...$
\begin{equation}
\left| III\right| \leq C(\varepsilon )h^{3/2-2\varepsilon
}\frac{C^{r}}{%
\Gamma (\frac{r-1}{2})}T^{3\varepsilon +\frac{r-4}{2}}\zeta _{\sqrt{T}%
}(v-x).  \label{T304}
\end{equation}
To estimate $II$ we write
\[
E_{r-1}=\tilde{p}\otimes _{h}(H+M_{h}^{\prime \prime }+\sqrt{h}%
N_{1})^{(r-1)}-\tilde{p}\otimes
_{h}(H+M_{h}+\sqrt{h}N_{1})^{(r-1)}.
\]
For $r=2$ we have $E_{1}=\tilde{p}\otimes _{h}A$ and analogously
to (%
\ref{eq:T206})
$$
\left| E_{1}\right| \leq Ch^{1-\varepsilon }(kh)^{\varepsilon
-1/2}B(\frac{1%
}{2},\varepsilon )\zeta _{\sqrt{kh}}^{S}(y-x). $$ For $r\geq 3$
 similarly as in (\ref{eq:T208}) we use the recurrence relation
\begin{eqnarray*}
 E_{r-1}&=&E_{r-2}\otimes _{h}(H+M_{h}^{\prime \prime
}+\sqrt{h}N_{1})+[\tilde{p%
}\otimes _{h}(H+M_{h}+\sqrt{h}N_{1})^{(r-2)}]\otimes _{h}A
\\ &=& I^{\prime }+II^{\prime }.
\end{eqnarray*}
The terms $I^{\prime }$ and $II^{\prime }$ have a similar structure
as the corresponding terms
in (%
\ref{eq:T225}) and they can be estimated similarly. This gives the
following estimates for $r=2,3,...$
\begin{eqnarray*}
\left| E_{r-1}\right| &\leq& C^{r}h^{1-\varepsilon
}B(\frac{1}{2},\varepsilon )\times ...\times
B(\frac{1}{2},\varepsilon +\frac{r-2}{2})(kh)^{\varepsilon
+\frac{r-3}{2}}\zeta _{\sqrt{kh}}^{S}(y-x),
\\
\left| II\right| &=&\left| E_{r-1}\otimes _{h}(M_{h}^{\prime \prime
}+\sqrt{h}%
N_{1})(0,T,x,y)\right|
\\
&\leq& C(\varepsilon )h^{3/2-2\varepsilon }\frac{C^{r}}{\Gamma
(\frac{r-1}{2})}%
T^{3\varepsilon +\frac{r-4}{2}}\zeta _{\sqrt{T}}(v-x).
\end{eqnarray*}
The claim (\ref{T301}) follows from (\ref{T302}), (\ref{T304}) and
the last two inequalities.
\bigskip

\noindent \textit{Asymptotic treatment of the term $T_{4}$. }We will
show that
\begin{eqnarray} \label{eqadd:3a}
T_{4}&=&\sum_{r=1}^{\infty }\tilde{p}\otimes _{h}H^{(r)}(0,T,x,y)
\\ \nonumber && \qquad
-\sum_{r=1}^{\infty }\tilde{p}\otimes
_{h}[H+hN_{2}]^{(r)}(0,T,x,y)+R_{h}^{\ast }(x,y),
\end{eqnarray}
with $N_{2}(s,t,x,y)=(L-\widetilde{L})\widetilde{\pi
}_{2}(s,t,x,y),\left| R_{h}^{\ast }(x,y)\right| \leq Chn^{-\delta
}\zeta _{\sqrt{T}}^{S}(y-x)$  for $\delta >0$ small enough  and with
a constant $C$  depending on $\delta $. For the proof of
(\ref{eqadd:3a}) it suffices to show that for $\delta $ small enough
\begin{eqnarray} \label{eqadd:3b} &&
\left| \sum_{r=1}^{n}\tilde{p}\otimes _{h}(H+M_{h}+\sqrt{h}%
N_{1}+hN_{2})^{(r)}(0,T,x,y)\right.
\\ \nonumber && \qquad
\left. -\sum_{r=1}^{n}\tilde{p}\otimes
_{h}(K_{h}+M_{h})^{(r)}(0,T,x,y)\right|
\\ \nonumber && \leq \left[ \sum_{k=1}^{n}\frac{C^{k}%
}{\Gamma (\frac{k}{2})}\right] hn^{-\delta }\zeta
_{\sqrt{T}}^{S}(y-x),
\\ \label{eqadd:3c} &&
\left| \sum_{r=1}^{n}\tilde{p}\otimes _{h}(H+M_{h}+\sqrt{h}%
N_{1})^{(r)}(0,T,x,y)\right.
\\ \nonumber && \qquad
-\sum_{r=1}^{n}\tilde{p}\otimes _{h}(H+M_{h}+\sqrt{h}%
N_{1}+hN_{2})^{(r)}(0,T,x,y)
\\ \nonumber && \qquad
\left. -\left[ \sum_{r=1}^{n}\tilde{p}\otimes
_{h}H^{(r)}(0,T,x,y)-\sum_{r=1}^{n}\tilde{p}\otimes
_{h}[H+hN_{2}]^{(r)}(0,T,x,y)\right] \right|
\\ \nonumber &&
\leq \left[ \sum_{k=1}^{n}\frac{C^{k}}{\Gamma (\frac{k}{2})}\right]
Chn^{-\delta }\zeta _{\sqrt{T}}^{S}(y-x).
\end{eqnarray}
Denote $D_{3,0}\equiv 0$ and
\begin{eqnarray*}
D_{3,m}(0,jh,x,y)&=&\sum_{r=1}^{m}\tilde{p}\otimes
_{h}(K_{h}+M_{h})^{(r)}(0,jh,x,y)
\\ \nonumber &&
-\sum_{r=1}^{m}\tilde{p}\otimes _{h}(H+M_{h}+\sqrt{h}%
N_{1}+hN_{2})^{(r)}(0,jh,x,y).
\end{eqnarray*}
Then (\ref{eqadd:3b}) can be rewritten as
\[
\left| D_{3,n}(0,T,x,y)\right| \leq Chn^{-\delta }\zeta _{\sqrt{T}%
}^{S}(y-x).
\]
We now make iterative use of
\begin{equation}
D_{3,m}=D_{3,m-1}\otimes _{h}(H+M_{h}+\sqrt{h}N_{1}+hN_{2})+g_{m-1},
\label{eqadd:3d}
\end{equation}
for $m=1,2,...,$ \ where
\begin{eqnarray*}
g_{m}(0,jh,x,y)&=&-\left[ \sum_{r=0}^{m}\tilde{p}\otimes
_{h}(K_{h}+M_{h})^{(r)}\right] \otimes _{h}(H-K_{h}+\sqrt{h}%
N_{1}+hN_{2})(0,jh,x,y)
\\
&=&S_{h,m}\otimes _{h}(L-\widetilde{L})d_{h}(0,jh,x,y)
\end{eqnarray*}
with
\begin{eqnarray*}
g_{0}(0,jh,x,y)&=&-\tilde{p}\otimes _{h}(H-K_{h}+\sqrt{h}%
N_{1}+hN_{2})(0,jh,x,y),
\\
d_{h}&=&\widetilde{p}_{h}-\widetilde{p}-\sqrt{h}\widetilde{\pi }_{1}-h%
\widetilde{\pi }_{2},
\\
S_{h,m}(0,ih,x,y)&=&\sum_{r=0}^{m}\tilde{p}\otimes
_{h}(K_{h}+M_{h})^{(r)}(0,ih,x,y).
\end{eqnarray*}
Iterative application of (\ref{eqadd:3d}) gives
\[
D_{3,n}(0,T,x,y)=\sum_{r=0}^{n-1}g_{r}\otimes _{h}(H+M_{h}+\sqrt{h}%
N_{1}+hN_{2})^{(n-r-1)}(0,T,x,y).
\]
To prove (\ref{eqadd:3b}) we will show that
\begin{eqnarray} \label{eqadd:3e}
&&\left| g_{r}\otimes _{h}(H+M_{h}+\sqrt{h}N_{1}+hN_{2})^{(n-r-1)}(0,T,x,y)%
\right|
\\ \nonumber && \qquad
\leq \frac{C^{n-r}}{\Gamma (\frac{n-r}{2})}hn^{-\delta }\zeta _{\sqrt{T}%
}^{S}(y-x).
\end{eqnarray}
For this purpose we decompose the left handside of (\ref{eqadd:3e})
into four terms
\begin{eqnarray*}
a_{r,1}&=&\sum_{0\leq i\leq n/2}h\int g_{r}(0,ih,x,u)(H+M_{h}+\sqrt{h}%
N_{1}+hN_{2})^{(n-r-1)}(ih,T,u,y)du,
\\
a_{r,2}&=&\sum_{n/2<i\leq n}h^{2}\sum_{0\leq k\leq i/2}\int \int
S_{h,r}(0,kh,x,v)(L-\widetilde{L})d_{h}(kh,ih,v,u)
\\ &&
\times (H+M_{h}+\sqrt{h}N_{1}+hN_{2})^{(n-r-1)}(ih,T,u,y)dvdu,
\\
a_{r,3}&=&\sum_{n/2<i\leq n}h^{2}\sum_{i/2<k\leq i-n^{\delta
^{\prime }}}\int \int
(L^{T}-\widetilde{L}^{T})S_{h,r}(0,kh,x,v)d_{h}(kh,ih,v,u)
\\ &&
\times (H+M_{h}+\sqrt{h}N_{1}+hN_{2})^{(n-r-1)}(ih,T,u,y)dvdu,
\\
a_{r,4}&=&\sum_{n/2<i\leq n}h^{2}\sum_{i-n^{\delta ^{\prime }}<k\leq
i-1}\int \int
(L^{T}-\widetilde{L}^{T})S_{h,r}(0,kh,x,v)d_{h}(kh,ih,v,u)
\\ &&
\times (H+M_{h}+\sqrt{h}N_{1}+hN_{2})^{(n-r-1)}(ih,T,u,y)dvdu.
\end{eqnarray*}
Here $L^{T}$ and $\widetilde{L}^{T}$ denote the adjoint operators of
$L$ and $\widetilde{L}$, and $\delta ^{\prime }$ satisfies
inequalities $2\varkappa <\delta ^{\prime }<\frac{3}{5}(1-\varkappa
),$ where $\varkappa $ is defined in (B2). For the proof of
(\ref{eqadd:3e}) it suffices to show for $l=1,2,3,4$%
\begin{eqnarray} \label{eqadd:3f}
\left| a_{r,l}\right| &\leq& hn^{-\delta
}C^{n-r}B(1,\frac{1}{2})\times ...\times
B(\frac{n-r-1}{2},\frac{1}{2})\zeta _{\sqrt{T}}^{S}(y-x)
\\ \nonumber
&\leq & \left[ \frac{C^{n-r}}{\Gamma (\frac{n-r}{2})}\right]
hn^{-\delta }\zeta _{\sqrt{T}}^{S}(y-x)
\end{eqnarray}
for some $\delta >0.$

\textit{Proof \ of \ (\ref{eqadd:3f}) for }$l=2.$ Note that \ $k\leq
i/2,i>n/2$ \ imply \ $ih-kh\geq \frac{T}{4}$. The claim follows from
the inequalities
\begin{eqnarray} \label{eqadd:3g}
&& \max \{\left| K_{h}(ih,jh,x,y)\right| ,\left|
M_{h}(ih,jh,x,y)\right| ,\left| \sqrt{h}N_{1}(ih,jh,x,y)\right|
\\ \nonumber && \qquad
\left| hN_{2}(ih,jh,x,y)\right| ,\left| H(ih,jh,x,y)\right| \}
\\ \nonumber &&
\leq C\rho ^{-1}\zeta _{\rho }(y-x)\mbox{ with }\rho ^{2}=jh-ih
\mbox { for }0\leq i<j\leq n,
\end{eqnarray}
\begin{eqnarray} \label{eqadd:3h} &&
\left| S_{h,m}(0,kh,x,v)\right| \leq C\zeta _{\sqrt{kh}}^{S-2}(v-x),
\\ \label{eqadd:03a} &&
\left| (L-\widetilde{L})d_{h}(kh,ih,v,u)\right| \leq
Ch^{3/2}(ih-kh)^{-2}\zeta _{\sqrt{ih-kh}}^{S-8}(u-v)
\\ \nonumber &&
=O(hn^{-1/2+3/2\varkappa })\zeta _{\sqrt{ih-kh}}^{S-8}(u-v),
\\ \label{eqadd:03b} &&
\left| (H+M_{h}+\sqrt{h}N_{1}+hN_{2})^{(n-r-1)}(ih,T,u,y)\right|
\\ \nonumber &&
\leq C^{n-r}\rho ^{n-r-3}B(\frac{1}{2},\frac{1}{2})\times ...\times B(\frac{%
n-r-2}{2},\frac{1}{2})\zeta _{\sqrt{T-ih}}^{S-2}(y-u)
\\ \nonumber &&
\leq \left[ \frac{C^{n-r}}{\Gamma (\frac{n-r-1}{2})}\right]
(T-ih)^{-1/2}\zeta _{\sqrt{T-ih}}^{S-2}(y-u)
\end{eqnarray}
for $n-r-3=-1,0,1,...,n-3$ with $\rho ^{2}=T-ih$. We put
$B(\frac{1}{2},0)=1)$.
Inequality (\ref{eqadd:3g}) follows from the definitions of the functions $%
K_{h},...,H$. Inequalities (\ref{eqadd:3h}) and (\ref{eqadd:03b})
can be proved by the same method as used in the proof of Theorem 2.3
in Konakov and Mammen (2002) (pp. 282 - 284). Inequality
(\ref{eqadd:03a}) follows from the inequality $ih-kh\geq
\frac{T}{4}$, Lemma 5 and the arguments used in the proof of Lemma
7.

\textit{Proof \ of \ (\ref{eqadd:3f}) for }$l=3.$ Note that
$n/2<i,k>i/2$ \ imply $\ kh>\frac{T}{4}.$ We use the following
inequalities
\begin{eqnarray} \label{eqadd:03ba} &&
\left| d_{h}(kh,ih,v,u)\right| \leq Ch^{3/2}(ih-kh)^{-3/2}\zeta _{\sqrt{ih-kh%
}}^{S-6}(u-v),
\\ \nonumber &&
\left| (L^{T}-\widetilde{L}^{T})S_{h,r}(0,kh,x,v)\right| \leq CT^{-1}\zeta _{%
\sqrt{kh}}^{S-2}(v-x),
\\ \label{eqadd:03c} &&
\left| h\sum_{i/2<k\leq i-n^{\delta ^{\prime }}}\int (L^{T}-\widetilde{L}%
^{T})S_{h,r}(0,kh,x,v)d_{h}(kh,ih,v,u)dv\right|
\\ \nonumber &&
\leq Ch^{3/2}T^{-1}\sum_{i/2<k\leq i-n^{\delta ^{\prime }}}h\frac{1}{%
(ih-kh)^{3/2}}\zeta _{\sqrt{ih}}(u-x)
\\ \nonumber &&
\leq Ch^{3/2}T^{-1}\int_{ih/2}^{ih-n^{\delta ^{\prime }}h}\frac{du}{%
(ih-u)^{3/2}}\zeta _{\sqrt{ih}}(u-x)\leq
Ch^{3/2}T^{-3/2}n^{(1-\delta ^{\prime })/2}\zeta _{\sqrt{ih}}(u-x)
\\ \nonumber &&
\leq Chn^{-\delta ^{\prime \prime }}\zeta _{\sqrt{ih}}(u-x),
\end{eqnarray}
where $\delta ^{\prime \prime }=\delta ^{\prime }/2-\varkappa >0.$
Claim (\ref {eqadd:3f}) for $l=3$ now follows from (\ref{eqadd:03c})
and (\ref{eqadd:03b}).

\textit{Proof of  (\ref{eqadd:3f}) for }$l=4$. For $i-n^{\delta
^{\prime }}<k\leq i-1,n/2<i$ we have $ih>T/2,$ $kh>T/3,$
$(i-k)<n^{\delta ^{\prime }}$ for sufficiently large $n$. The
integral
\[
\int (L^{T}-\widetilde{L}^{T})S_{h,r}(0,kh,x,v)\widetilde{p}%
_{h}(kh,ih,v,u)dv
\]
is a finite sum of integrals. We show how to estimate a typical term
of this sum. The other terms can be estimated analogously. We
consider for
fixed $j,l$%
\begin{eqnarray} \label{eqadd:03d} &&
\int \frac{\partial ^{2}S_{h,r}(0,kh,x,v)}{\partial v_{j}\partial v_{l}}%
(\sigma _{jl}(kh,v)-\sigma _{jl}(kh,u))h^{-d/2}
\\ \nonumber &&
\qquad \times
q^{(i-k)}[kh,u,h^{-1/2}(u-v-h\sum_{l=k}^{i-1}m(lh,u))]dv
\\ \nonumber &&
=\int \frac{\partial ^{2}S_{h,r}(0,kh,x,v)}{\partial v_{j}\partial v_{l}}%
\mid _{v=u^{\ast }-\sqrt{h}w}
\\ \nonumber && \qquad
\times \lbrack \sigma _{jl}(kh,u^{\ast }-\sqrt{h}w)-\sigma
_{jl}(kh,u)]q^{(i-k)}(kh,u,w)dw,
\end{eqnarray}
where $u^{\ast }=u-h\sum_{l=k}^{i-1}m(lh,u).$ Now using a Tailor
expansion we obtain that the right hand side of (\ref{eqadd:03d}) is
equal to
\begin{eqnarray*} &&
\int \left[ \frac{\partial ^{2}S_{h,r}(0,kh,x,u^{\ast })}{\partial
v_{j}\partial v_{l}}-\sqrt{h}\sum_{\left| \nu \right|
=1}\frac{w^{\nu }}{\nu !}\int_{0}^{1}D_{v}^{\nu }\frac{\partial
^{2}S_{h,r}(0,kh,x,u^{\ast }-\delta \sqrt{h}w)}{\partial
v_{j}\partial v_{l}}d\delta \right]
\\ &&
\times \left[ -\sqrt{h}\sum_{\left| \nu \right| =1}\frac{[w+\sqrt{h}%
\sum_{l=k}^{i-1}m(lh,u)]^{\nu }}{\nu !}D_{u}^{\nu }\sigma
_{jl}(kh,u)\right.
\\ && \qquad +2h\sum_{\left| \nu \right| =2}\frac{[w+\sqrt{h}%
\sum_{l=k}^{i-1}m(lh,u)]^{\nu }}{\nu !}
\\ &&
\left. \qquad \times \int_{0}^{1}D_{u}^{\nu }\sigma _{jl}(kh,u-\delta \sqrt{h}%
w-\delta h\sum_{l=k}^{i-1}m(lh,u))d\delta \right]
q^{(i-k)}(kh,u,w)dw.
\end{eqnarray*}
Note that
\begin{eqnarray}\nonumber &&
-\sqrt{h}\int \frac{\partial ^{2}S_{h,r}(0,kh,x,u^{\ast })}{\partial
v_{j}\partial v_{l}}(w_{p}+\sqrt{h}%
\sum_{l=k}^{i-1}m_{p}(lh,u))q^{(i-k)}(kh,u,w)dw
\\ \nonumber &&
=-h\frac{\partial ^{2}S_{h,r}(0,kh,x,u^{\ast })}{\partial v_{j}\partial v_{l}%
}\sum_{l=k}^{i-1}m_{p}(lh,u),
\\  \label{eqadd:03e} &&
 h\int \frac{\partial ^{2}S_{h,r}(0,kh,x,u^{\ast
})}{\partial v_{j}\partial
v_{l}}(w_{p}+\sqrt{h}\sum_{l=k}^{i-1}m_{p}(lh,u))(w_{q}+\sqrt{h}%
\sum_{l=k}^{i-1}m_{q}(lh,u))
\\ \nonumber && \qquad
\times \left\{ \int_{0}^{1}D_{u}^{\nu }\sigma _{jl}(kh,u)d\delta
+\int_{0}^{1}\left[ D_{u}^{\nu }\sigma _{jl}(kh,u-\delta \sqrt{h}w
\right . \right . \\ \nonumber && \qquad \qquad \left. \left .
-\delta h\sum_{l=k}^{i-1}m(lh,u))-D_{u}^{\nu }\sigma
_{jl}(kh,u)\right] d\delta \right\} q^{(i-k)}(kh,u,w)dw
\\ \nonumber &&
=h\frac{\partial ^{2}S_{h,r}(0,kh,x,u^{\ast })}{%
\partial v_{j}\partial v_{l}}D_{u}^{\nu }\sigma _{jl}(kh,u)\int
w_{p}w_{q}q^{(i-k)}(kh,u,w)dw
\\ \nonumber && \qquad
+h^{2}\frac{\partial ^{2}S_{h,r}(0,kh,x,u^{\ast })}{\partial
v_{j}\partial
v_{l}}\sum_{l=k}^{i-1}m_{p}(lh,u)\sum_{l=k}^{i-1}m_{q}(lh,u)+R,
\end{eqnarray}
where by (A3$^{\prime }$) we have for $%
j_{0}<(i-k)<n^{\delta ^{\prime }},w^{\prime }=(i-k)^{-1/2}w$ \
\textbf{\ }
\begin{eqnarray} \nonumber
\left| R\right| &\leq& Ch^{3/2}\left| \frac{\partial
^{2}S_{h,r}(0,kh,x,u^{\ast })}{\partial v_{j}\partial v_{l}}\right|
\int \left( n^{\delta ^{\prime }/2}\left\| w^{\prime }\right\|
+O(T^{1/2}n^{-1/2+\delta ^{\prime }})\right) ^{3}\psi (w^{\prime
})dw^{\prime }
\\ \label{eqadd:03f}
&\leq& Ch\zeta _{\sqrt{ih}}(u-x)(h^{S}n^{S\delta ^{\prime
}}+1)T^{-3/2}n^{-1/2+3\delta ^{\prime }/2}\int \left\| w^{\prime
}\right\| ^{3}\psi (w^{\prime })dw^{\prime } \\ \nonumber &\leq&
Chn^{-1/2+\varkappa /2+3\delta ^{\prime }/2}\zeta
_{\sqrt{ih}}(u-x)\leq Chn^{-1/2(1-3\delta ^{\prime }-\varkappa
)}\zeta _{\sqrt{ih}}(u-x),
\end{eqnarray}
We obtain analogously
\begin{eqnarray*}
&& h\int
w_{p}(w_{q}+\sqrt{h}\sum_{l=k}^{i-1}m_{q}(lh,u))D_{u}^{e_{q}}\sigma
_{jl}(kh,u)\int_{0}^{1}\frac{\partial ^{3}S_{h,r}(0,kh,x,u^{\ast
}-\delta \sqrt{h}w)}{\partial v_{p}\partial v_{j}\partial
v_{l}}d\delta
\\ && \qquad
\times q^{(i-k)}(kh,u,w)dw
\\ &&
=h\frac{\partial \sigma _{jl}(kh,u)}{\partial u_{q}%
}\frac{\partial ^{3}S_{h,r}(0,kh,x,u^{\ast })}{\partial
v_{p}\partial v_{j}\partial v_{l}}\int
w_{p}w_{q}q^{(i-k)}(kh,u,w)dw+R,
\end{eqnarray*}
where
\[
\left| R\right| \leq Chn^{-1/2(1-3\delta ^{\prime }-3\varkappa )}\zeta _{%
\sqrt{ih}}(u-x),1-3\delta ^{\prime }-3\varkappa >0
\]
and, for $1-3\delta ^{\prime }-2\varkappa >0$
\begin{eqnarray} \label{eqadd:03g} &&
\left| h^{3/2}\int w_{p}\int_{0}^{1}\frac{\partial
^{3}S_{h,r}(0,kh,x,u^{\ast }-\delta \sqrt{h}w)}{\partial
v_{p}\partial v_{j}\partial v_{l}}d\delta \right .\\ \nonumber
&&\left. \qquad \int_{0}^{1}\frac{\partial ^{2}\sigma
_{jl}(kh,u-\delta \sqrt{h}w-\delta
h\sum_{l=k}^{i-1}m_{q}(lh,u))}{\partial u_{r}\partial u_{s}}d\delta
\right.
\\ \nonumber && \qquad \qquad
\left. \times (w_{r}+\sqrt{h}\sum_{l=k}^{i-1}m_{r}(lh,u))(w_{s}+\sqrt{h}%
\sum_{l=k}^{i-1}m_{s}(lh,u))q^{(i-k)}(kh,u,w)dw\right|
\\ \nonumber &&
\leq Chn^{-1/2(1-3\delta ^{\prime }-2\varkappa )}\zeta _{\sqrt{ih}%
}(u-x).
\end{eqnarray}
For $1\leq i-k\leq j_{0}$ \ the same estimates remain true \ because
the following bound holds \
\begin{equation}
\int \left\| w\right\| ^{S}q^{(j)}(t,x,w)dw\leq C(j_{0}).
\label{eqadd:03h}
\end{equation}
The same estimates hold for $\widetilde{p}(kh,ih,v,u)$ \ with \
$\phi ^{(i-k)}(kh,u,w)$ \ instead of $q^{(i-k)}(kh,u,w),$ where
$\phi (kh,u,w)$ \ is a gaussian density with the mean \ $0$ \ and \
with \ the covariance matrix equal to \ $\sigma (kh,u).$ The first
two moments of $q^{(i-k)}$ \ and $\phi ^{(i-k)}$ \ coinside so after
substraction we obtain uniformly for \ $i-n^{\delta ^{\prime
}}<k\leq i-1$
\begin{eqnarray} \label{eqadd:03i} &&
\left| \sum_{n/2<i\leq n}h^{2}\sum_{i-n^{\delta ^{\prime }}<k\leq
i-1}\int \int (L^{T}-\widetilde{L}^{T})S_{h,r}(0,kh,x,v)\right. \\
\nonumber &&
\qquad \left. \times(\widetilde{p}_{h}(kh,ih,v,u)-%
\widetilde{p}(kh,ih,v,u))dv\right.
\\ \nonumber && \qquad
\left. \times
(H+M_{h}+\sqrt{h}N_{1}+hN_{2})^{(n-r-1)}(ih,T,u,y)du\right|
\\ \nonumber &&
\leq \left[ \frac{C^{n-r}}{\Gamma (\frac{n-r-1}{2})}\right]
hT^{3/2}n^{-3/2(1-\varkappa -5\delta ^{\prime }/3)}\zeta
_{\sqrt{ih}}(u-x).
\end{eqnarray}
To estimate the other terms in $\ d_{h}(kh,ih,v,u)$ we need bounds
for the following expressions
\begin{eqnarray*} &&
h\sum_{i-n^{\delta ^{\prime }}<k\leq i-1}\int (L^{T}-\widetilde{L}%
^{T})S_{h,r}(0,kh,x,v)\sqrt{h}(ih-kh)\\
&& \qquad \times D_{v}^{\nu }\widetilde{p}(kh,ih,v,u)dv \mbox{ for }
\left| \nu \right| =3,
\\
&&
h\sum_{i-n^{\delta ^{\prime }}<k\leq i-1}\int (L^{T}-\widetilde{L}%
^{T})S_{h,r}(0,kh,x,v)h(ih-kh)\\
&& \qquad \times D_{v}^{\nu }\widetilde{p}(kh,ih,v,u)dv \mbox{ for }
\left| \nu \right| =4,
\\
&&
h\sum_{i-n^{\delta ^{\prime }}<k\leq i-1}\int (L^{T}-\widetilde{L}%
^{T})S_{h,r}(0,kh,x,v)h(ih-kh)^{2}\\
&& \qquad \times D_{v}^{\nu }\widetilde{p}(kh,ih,v,u)dv \mbox{ for }
\left| \nu \right| =6.
\end{eqnarray*}
We have
\begin{eqnarray} \label{eqadd:03j} &&
\left| h\sum_{i-n^{\delta ^{\prime }}<k\leq i-1}\int (L^{T}-\widetilde{L}%
^{T})S_{h,r}(0,kh,x,v)\sqrt{h}(ih-kh)\right.
\\ \nonumber && \qquad \left. D_{v}^{\nu }\widetilde{p}%
(kh,ih,v,u)dv\right|
\\ \nonumber &&
=\left| h\sum_{i-n^{\delta ^{\prime }}<k\leq i-1}\int D^{e_{p}+e_{q}}(L^{T}-%
\widetilde{L}^{T})S_{h,r}(0,kh,x,v)\sqrt{h}(ih-kh)\right.
\\ \nonumber && \qquad \left. D_{v}^{\nu -e_{p}-e_{q}}%
\widetilde{p}(kh,ih,v,u)dv\right| \\\nonumber  && \leq
CT^{-2}n^{\delta ^{\prime }}h^{3/2}\sum_{i-n^{\delta ^{\prime
}}<k\leq i-1}\frac{h}{\sqrt{ih-kh}}\zeta _{\sqrt{ih}}(u-x)
\\ \nonumber &&
\leq Chn^{-(1-\varkappa -3\delta ^{\prime }/2)}\zeta
_{\sqrt{ih}}(u-x). \nonumber
\end{eqnarray}
Clearly, the same estimate (\ref{eqadd:03j}) holds for \ $\left| \nu
\right| =4$ \ and $\left| \nu \right| =6.$ \ Now (\ref{eqadd:3f})
for $l=4$ \ follows from this remark and (\ref{eqadd:03i}) and
(\ref{eqadd:03j}).

\textit{Proof \ of \ (\ref{eqadd:3f}) for }$l=1.$ Note that for this case $%
T-ih\geq T/2.$%
\begin{eqnarray}  \label{eqadd:03ja} &&
a_{r,1}=\sum_{0\leq i\leq n/2}h^{2}\sum_{0\leq k\leq i-1}\int \int (L^{T}-%
\widetilde{L}^{T})S_{h,r}(0,kh,x,v)\\ \nonumber && \qquad \times
d_{h}(kh,ih,v,u)\Psi _{h,r}(ih,T,u,y)dvdu
\\ \nonumber &&
=\sum_{0\leq k\leq n/2-1}h\int (L^{T}-\widetilde{L}^{T})S_{h,r}(0,kh,x,v)%
\\ \nonumber && \qquad \times \left\{ \sum_{k+1\leq i\leq k+n^{\delta ^{\prime }}}h\int
d_{h}(kh,ih,v,u)\Psi _{h,r}(ih,T,u,y)du\right.
\\ \nonumber && \qquad
\left. +\sum_{k+n^{\delta ^{\prime }}<i\leq n/2}h\int
d_{h}(kh,ih,v,u)\Psi _{h,r}(ih,T,u,y)du\right\} dv,
\end{eqnarray}
where we denote
\[
\Psi
_{h,r}(ih,T,u,y)=(H+M_{h}+\sqrt{h}N_{1}+hN_{2})^{(n-r-1)}(ih,T,u,y).
\]
We consider
\begin{eqnarray*} &&
\sum_{k+1\leq i\leq k+n^{\delta ^{\prime }}}h\int
h^{-d/2}q^{(i-k)}(kh,u,h^{-1/2}[u-v-h\sum_{l=k}^{i-1}m(lh,u)])\Psi
_{h,r}(ih,T,u,y)du
\\&&
=\sum_{k+1\leq i\leq k+n^{\delta ^{\prime }}}h\int \left\{ q^{(i-k)}(kh,v,w)+%
\sqrt{h}\sum_{\left| \nu \right|
=1}(w+\sqrt{h}\sum_{l=k}^{i-1}m(lh,u))^{\nu }D_{v}^{\nu
}q^{(i-k)}(kh,v,w)\right.
\\&& \qquad
+h\sum_{\left| \nu \right|
=2}\frac{(w+\sqrt{h}\sum_{l=k}^{i-1}m(lh,u))^{\nu }}{\nu
!}D_{v}^{\nu }q^{(i-k)}(kh,v,w)
\\&& \qquad
+3h^{3/2}\sum_{\left| \nu \right| =3}%
\frac{(w+\sqrt{h}\sum_{l=k}^{i-1}m(lh,u))^{\nu }}{\nu !}
\\&& \qquad
\left. \times \int_{0}^{1}(1-\delta )^{2}D_{v}^{\nu
}q^{(i-k)}(kh,v+\delta h^{1/2}w+\delta
h\sum_{l=k}^{i-1}m(lh,u),w)d\delta \right\}
\\&& \qquad
\times \left\{ \Psi _{h,r}(ih,T,v,y)+\sqrt{h}\sum_{\left| \nu \right| =1}(w+%
\sqrt{h}\sum_{l=k}^{i-1}m(lh,u))^{\nu }D_{v}^{\nu }\Psi
_{h,r}(ih,T,v,y)\right.
\\&& \qquad
+h\sum_{\left| \nu \right|
=2}\frac{(w+\sqrt{h}\sum_{l=k}^{i-1}m(lh,u))^{\nu }}{\nu
!}D_{v}^{\nu }\Psi _{h,r}(ih,T,v,y)
\\&& \qquad +3h^{3/2}\sum_{\left| \nu
\right| =3}\frac{(w+\sqrt{h}\sum_{l=k}^{i-1}m(lh,u))^{\nu }}{\nu !}
\\&& \qquad
\left. \times \int_{0}^{1}(1-\delta )^{2}D_{v}^{\nu }\Psi
_{h,r}(ih,T,v+\delta h^{1/2}w+\delta
h\sum_{l=k}^{i-1}m(lh,u),y)d\delta \right\} dw
\end{eqnarray*}
This integral is a sum of $4\times 4=16$ integrals. We estimate only
two of \ them. Other integrals can be estimated by similar methods.
First, we estimate
\[
\sum_{k+1\leq i\leq k+n^{\delta ^{\prime }}}h\int
q^{(i-k)}(kh,v,w)\Psi _{h,r}(ih,T,v,y)dw=\sum_{k+1\leq i\leq
k+n^{\delta ^{\prime }}}h\Psi _{h,r}(ih,T,v,y)dw.
\]
Note that we get the same term when we replace $q^{(i-k)}(kh,v,w)$ by $%
\phi ^{(i-k)}(kh,v,w)$. After the replacement this term disappears.
Second, we estimate
\begin{eqnarray*} &&
\sum_{k+1\leq i\leq k+n^{\delta ^{\prime }}}h\int q^{(i-k)}(kh,v,w)\sqrt{h}%
\sum_{\left| \nu \right| =1}(w+\sqrt{h}\sum_{l=k}^{i-1}m(lh,u))^{\nu
}D_{v}^{\nu }\Psi _{h,r}(ih,T,v,y)dw
\\ &&
=h^{3/2}\sum_{j=1}^{d}\sum_{k+1\leq i\leq k+n^{\delta ^{\prime
}}}D_{v}^{e_{j}}\Psi _{h,r}(ih,T,v,y)\int q^{(i-k)}(kh,v,w)[w_{j}+\sqrt{h}%
\sum_{l=k}^{i-1}m_{j}(lh,v)
\\ && \qquad
+O(hn^{\delta ^{\prime }}\left\| w\right\| +h^{3/2}n^{2\delta
^{\prime }})]dw \\ && =h^{2}\sum_{j=1}^{d}\sum_{k+1\leq i\leq
k+n^{\delta ^{\prime }}}D_{v}^{e_{j}}\Psi
_{h,r}(ih,T,v,y)\sum_{l=k}^{i-1}m_{j}(lh,v)
\\ && \qquad
+O\left( h^{2}n^{2\delta ^{\prime }}\sum_{j=1}^{d}\sum_{k+1\leq
i\leq k+n^{\delta ^{\prime }}}h\left| D_{v}^{e_{j}}\Psi
_{h,r}(ih,T,v,y)\right| \right) \\ && \qquad +O\left(
h^{3/2}n^{\delta ^{\prime }}\sum_{j=1}^{d}\sum_{k+1\leq i\leq
k+n^{\delta ^{\prime }}}h \left| D_{v}^{e_{j}}\Psi
_{h,r}(ih,T,v,y)\right| \int q^{(i-k)}(kh,v,w)\left\| w\right\|
dw\right)
\\ &&
=h^{2}\sum_{j=1}^{d}\sum_{k+1\leq i\leq k+n^{\delta ^{\prime
}}}D_{v}^{e_{j}}\Psi _{h,r}(ih,T,v,y)\sum_{l=k}^{i-1}m_{j}(lh,v)+R,
\end{eqnarray*}
where
\[
\left| R\right| \leq \frac{C^{n-r}}{\Gamma (\frac{n-r-1}{2})}%
T^{1/2}hn^{-3/2+2\delta ^{\prime }}\zeta _{\sqrt{T-kh}}(y-v).
\]
The first term in the right hand side of this equation will be the
same if we replace $q^{(i-k)}(kh,v,w)$ \ \ by $\phi
^{(i-k)}(kh,v,w)$. After the replacement this term disappears. For a
proof of this equation we consider the function $u(w)$ that is
defined as an implicit function and we used the following change of
variables
\[
h^{1/2}w=u-v-h\sum_{l=k}^{i-1}m(lh,u)
\]
to obtain
\[
\sqrt{h}\sum_{l=k}^{i-1}m(lh,u(w))=\sqrt{h}\sum_{l=k}^{i-1}m(lh,v)+O\left(
h(i-k)\left\| w\right\| +h^{3/2}(i-k)^{2}\right)
\]
because of $\ (i-k)\leq n^{\delta ^{\prime }}$. By similar methods
we get
\begin{eqnarray} \label{eqadd:03m} &&
\left| \sum_{k+1\leq i\leq k+n^{\delta ^{\prime }}}h\int [\sqrt{h}\widetilde{%
\pi }_{1}(kh,ih,v,u)+h\widetilde{\pi }_{2}(kh,ih,v,u)]\Psi
_{h,r}(ih,T,u,y)du\right|
\\ \nonumber && \qquad
\leq \frac{C^{n-r}}{\Gamma (\frac{n-r-1}{2})}hn^{-3/2+2\delta
^{\prime }+\varkappa /2}\zeta _{\sqrt{T-kh}}(y-v).
\end{eqnarray}
It remains to estimate
\[
\sum_{k+n^{\delta ^{\prime }}<i\leq n/2}h\int d_{h}(kh,ih,v,u)\Psi
_{h,r}(ih,T,u,y)du.
\]
From (\ref{eqadd:03b}) and (\ref{eqadd:03ba}) we obtain
\begin{eqnarray} \label{eqadd:03p}
&& \left| \sum_{k+n^{\delta ^{\prime }}<i\leq n/2}h\int
d_{h}(kh,ih,v,u)\Psi
_{h,r}(ih,T,u,y)du\right| \\ \nonumber && \qquad \leq \left[ \frac{C^{n-r}}{\Gamma (\frac{n-r-1}{2})%
}\right] T^{-1/2}h^{3/2}
\int_{kh+n^{\delta ^{\prime }}h}^{T/2}\frac{du}{(u-kh)^{3/2}}\zeta _{%
\sqrt{T-kh}}(y-v)
\\ \nonumber && \qquad
\leq \left[ \frac{C^{n-r}}{\Gamma (\frac{n-r-1}{2})}\right]
T^{-1/2}hn^{-\delta ^{\prime }/2}\zeta _{\sqrt{T-kh}}(y-v)
\\ \nonumber && \qquad
\leq \left[ \frac{C^{n-r}}{\Gamma (\frac{n-r-1}{2})}\right]
hn^{-1/2(\delta ^{\prime }-\varkappa )}\zeta _{\sqrt{T-kh}}(y-v).
\end{eqnarray}
Now we substitute the estimate (\ref{eqadd:03p}) into
(\ref{eqadd:03ja}). This
gives the following estimate for any $0<\varepsilon <\varkappa $%
\begin{eqnarray} \label{eqadd:03q}&&
\left| \sum_{k=1}^{n/2-1}h\int
(L^{T}-\widetilde{L}^{T})S_{h,r}(0,kh,x,v) \right .
\\ \nonumber && \qquad \qquad\left .
\sum_{k+n^{\delta ^{\prime }}<i\leq n/2}h\int d_{h}(kh,ih,v,u)\Psi
_{h,r}(ih,T,u,y)du\right|
\\ \nonumber && \qquad
\leq \left[ \frac{C^{n-r}}{\Gamma (\frac{n-r-1}{2})}\right]
hn^{-1/2(\delta ^{\prime }-\varkappa )}h^{-\varepsilon
}\sum_{k=1}^{n/2}h(kh)^{\varepsilon -1}\zeta _{\sqrt{T}}(y-x)
\\ \nonumber && \qquad
\leq C(\varepsilon )\left[ \frac{C^{n-r}}{\Gamma
(\frac{n-r-1}{2})}\right] hn^{-1/2(\delta ^{\prime }-\varkappa
-\varepsilon )}\zeta _{\sqrt{T}}(y-x).
\end{eqnarray}
For $k=0$ we get with $S_{h,r}(0,0,x,v)=\delta (x-v)$ where
$\delta (\cdot )$ is the Dirac function that
\begin{eqnarray*} &&
\left| \sum_{1\leq i\leq i/2}h^{2}\int \int S_{h,r}(0,0,x,v)(L-\widetilde{L}%
)d_{h}(0,ih,v,u)\Psi _{h,r}(ih,T,u,y)du\right|
\\ && \qquad
\leq C(\varepsilon )\left[ \frac{C^{n-r}}{\Gamma
(\frac{n-r-1}{2})}\right] hn^{-(1/2-\varepsilon )}\zeta
_{\sqrt{T}}(y-x).
\end{eqnarray*}
This completes the proof (\ref{eqadd:3f}) for $l=1$. The estimate
(\ref {eqadd:3c}) may be proved by the same arguments as were used
in the treatment of $T_3$.
\bigskip

\noindent  \textit{Asymptotic treatment of the term $T_{5}.$ }We
will show that,
\begin{eqnarray}\nonumber
T_{5}&=&-\sqrt{h}\sum_{r=0}^{\infty }\widetilde{\pi }_{1}\otimes
_{h}(H+M_{h,1}+\sqrt{h}N_{1})^{(r)}(0,T,x,y)
\\ &&
-h\sum_{r=0}^{\infty }\widetilde{\pi }_{2}\otimes
_{h}H^{(r)}(0,T,x,y)+R_{h}(x,y),  \label{eq:T500}
\end{eqnarray}
where $\left| R_{h}(x,y)\right| \leq Chn^{-\gamma }\zeta _{\sqrt{T}%
}^{S-2}(y-x)$ for some $\gamma >0.$ Note that with $\
S_{h}(s,t,x,y)=\sum_{r=1}^{n}(K_{h}+M_{h})^{(r)}(s,t,x,y)$ the term
$T_{5}$ can be rewritten as
\[
T_{5}=(\widetilde{p}-\widetilde{p}_{h})(0,T,x,y)+(\widetilde{p}-\widetilde{p}%
_{h})\otimes _{h}S_{h}(0,T,x,y).
\]
We start by showing that for $\varkappa <\delta <\frac{1-\varkappa
}{4}$ \
uniformly for $x,y\in R$%
\begin{equation}
\left| h\sum_{1\leq j\leq n^{\delta }}\int
(\widetilde{p}_{h}-\widetilde{p}%
)(0,jh,x,u)S_{h}(jh,T,u,y)du\right| \leq O(hn^{-1/2(1-\varkappa
-4\delta )})\zeta _{\sqrt{T}}^{S-2}(y-x)  \label{eq:29a}
\end{equation}
for $\delta $ small enough. For the proof of (\ref{eq:29a}) we
will show
that uniformly for $1\leq j\leq n^{\delta }$ and for $x,y\in R^{d}$%
\begin{eqnarray} \nonumber &&
\int \widetilde{p}_{h}(0,jh,x,u)S_{h}(jh,T,u,y)du=S_{h}(jh,T,x,y)
\\ && \qquad
+O[h^{1/2}T^{-1/2}n^{-1/2+\delta }+h^{1/2}T^{-1}+n^{\delta
/2}h^{1/2}]\zeta _{\sqrt{T}}^{S-2}(y-x),  \label{eq:29b}
\\ && \nonumber
\int \widetilde{p}(0,jh,x,u)S_{h}(jh,T,u,y)du=S_{h}(jh,T,x,y)
\\ && \qquad
+O[h^{1/2}T^{-1/2}n^{-1/2+\delta }+h^{1/2}T^{-1}+n^{\delta
/2}h^{1/2}]\zeta _{\sqrt{T}}^{S-2}(y-x).  \label{eq:29c}
\end{eqnarray}
Claim (\ref{eq:29a}) immediately follows from
(\ref{eq:29b})-(\ref{eq:29c}%
). For the proof we will make use of the fact that for all $1\leq
j\leq n^{\delta }$ and for all $x,y\in R^{d}$ and $\left| \nu
\right|
=1$%
\begin{equation}
\left| D_{x}^{\nu }S_{h}(jh,T,x,y)\right| \leq C(T-jh)^{-1}\zeta
_{\sqrt{T-jh%
}}^{S-2}(y-x).  \label{eq:29d}
\end{equation}
Claim (\ref{eq:29d}) can be shown with the same arguments as in the
proof of (5.7) in Konakov and Mammen (2002). Note that the function
$\Phi $ in that paper has a similar structure as $S_{h}$. For $1\leq
j\leq n^{\delta }$ the bound (\ref{eq:29d}) immediately implies for
a constant
$C^{\prime }$%
\begin{equation}
\left| D_{x}^{\nu }S_{h}(jh,T,x,y)\right| \leq C^{\prime
}T^{-1}\zeta
_{%
\sqrt{T}}^{S-2}(y-x).  \label{eq:29e}
\end{equation}
We have
$\widetilde{p}_{h}(0,jh,x,u)=h^{-d/2}q^{(j)}[0,u,h^{-1/2}(u-x-h%
\sum_{i=0}^{j-1}m(ih,u))].$ Denote the determinant of the Jacobian
matrix of $u-h\sum_{i=0}^{j-1}m(ih,u)$ by $\Delta _{h}.$ From the
condition (A3) and (\ref{eq:29e}) we get that for $1\leq j\leq
n^{\delta }$
\begin{eqnarray*}
&&\int \widetilde{p}_{h}(0,jh,x,u)S_{h}(jh,T,u,y)du
\\ &&
=\int
h^{-d/2}q^{(j)}[0,u,h^{-1/2}(u-x-h\sum_{i=0}^{j-1}m(ih,u))]S_{h}(jh,T,u,y)du
\\ &&
=\int q^{(j)}(0,x+h^{1/2}w+h\sum_{i=0}^{j-1}m(ih,u(w)),w)\left|
\Delta _{h}^{-1}\right|
S_{h}(jh,T,x+h^{1/2}w+h\sum_{i=0}^{j-1}m(ih,u(w)),y)dw
\\ &&
=\int [q^{(j)}(0,x,w)+O(j^{-d/2}h^{1/2})(\left\| w\right\| +1)\psi
(j^{-1/2}w)][1+O(jh)][S_{h}(jh,T,x,y)
\\ && \qquad
+O(h^{1/2}T^{-1})\zeta _{\sqrt{T}}^{S-2}(y-x)(1+h^{(S-2)/2}\left\|
w\right\| ^{S-2})(\left\| w\right\| +1)]dw
\\ &&
=S_{h}(jh,T,x,y)+O[h^{1/2}T^{-1/2}n^{-1/2+\delta
}+h^{1/2}T^{-1}+h^{1/2}n^{\delta /2}]\zeta _{\sqrt{T}}^{S-2}(y-x)
\end{eqnarray*}
with $u=u(w)$ in $\sum_{i=0}^{j-1}m(ih,u)$ defined by the Inverse
Function Theorem from the substitution $w=$ $h^{-1/2}(u-x-h%
\sum_{i=0}^{j-1}m(ih,u))$. This proves (\ref{eq:29b}). Claim (\ref
{eq:29c}) follows by similar arguments. From (\ref{eq:29a}) we get
that for $\delta <\frac{1-\varkappa }{4}$ (with $\varkappa $ defined
as in (B2))
\[
T_{5}=(\widetilde{p}-\widetilde{p}_{h})(0,T,x,y)+h\sum_{n^{\delta
}<j<n}\int
(\widetilde{p}-\widetilde{p}_{h})(0,jh,x,u)S_{h}(jT,u,y)du+R_{h}(x,y)
\]
with $\left| R_{h}(x,y)\right| \leq O(hn^{-1/2(1-\varkappa
-4\delta )})\zeta
_{\sqrt{T}}^{S-2}(y-x)$ $.$ We now make use of the expansion of $%
\widetilde{p}_{h}-\widetilde{p}$ given in Lemma 5. We have with
$\rho =(jh)^{1/2}\geq h^{1/2}n^{\delta /2}$
\begin{equation}
\left| h\sum_{j=n^{\delta }}^{n}h^{3/2}\rho ^{-3}\int \zeta _{\rho
}^{S}(u-x)S_{h}(jh,T,u,y)du\right| \leq Ch^{2}T^{-\delta ^{\prime
}}n^{-\delta ^{\prime \prime }}\sum_{j=n^{\delta }}^{n}\rho
^{-2+2\delta ^{\prime }}\int \left| \zeta _{\rho
}^{S}(u-x)S_{h}(jh,T,u,y)\right| du, \label{eq:29g}
\end{equation}
where $\delta ^{\prime }<\frac{1}{2}\delta (1-\delta )^{-1},2$
$\delta ^{\prime \prime }=\delta +2\delta \delta ^{\prime
}-2\delta ^{\prime }.$ Now we get that
\begin{equation}
h\sum_{j=n^{\delta }}^{n}\rho ^{-2+2\delta ^{\prime }}\int \left|
\zeta _{\rho }^{S}(u-x)S_{h}(jh,T,u,y)\right| du\leq CB(\delta
^{\prime },1/2)T^{\delta ^{\prime }-1/2}\zeta
_{\sqrt{T}}^{S-2}(y-x) \label{eq:30}
\end{equation}
for a constant $C.$ This shows that for $\delta ^{\prime }>0$ small
enough
\begin{eqnarray*}
T_{5}&=&-[\sqrt{h}\widetilde{\pi }_{1}+h\widetilde{\pi
}_{2}](0,T,x,y)
\\ &&
-h\sum_{n^{\delta }<j<n}\int [\sqrt{h}\widetilde{\pi
}_{1}+h\widetilde{\pi }%
_{2}](0,jh,x,u)S_{h}(jh,T,u,y)du+R_{h}^{\prime }(x,y)
\end{eqnarray*}
with $\left| R_{h}^{\prime }(x,y)\right| \leq O(hn^{-(\delta
^{\prime \prime }-\varkappa /2)})\zeta _{\sqrt{T}}^{S-2}(y-x)$ \
with a constant in $O(\cdot )$ depending on $\delta ^{\prime }.$
It follows from (\ref{eq:29a}), (\ref {eq:29g}) and (\ref{eq:30})
that
\begin{equation}
T_{5}=-\sum_{r=0}^{\infty }[\sqrt{h}\widetilde{\pi
}_{1}+h\widetilde{\pi }%
_{2}]\otimes _{h}(K_{h}+M_{h})^{(r)}(0,T,x,y)+R_{h}^{\prime \prime
}(x,y), \label{eq:31}
\end{equation}
where $\left| R_{h}^{\prime \prime }(x,y)\right| \leq
O(hn^{-(\delta ^{\prime \prime }-\varkappa /2)})\zeta
_{\sqrt{T}}^{S-2}(y-x).$ Now we apply Lemma 10 with
$A=\sqrt{h}\widetilde{\pi }_{1},$
$B=H+M_{h,1}+\sqrt{%
h}N_{1},C=(K_{h}-H-\sqrt{h}N_{1})+(M_{h}-M_{h,1})$ to
\begin{equation}
-\sum_{r=0}^{\infty }\sqrt{h}\widetilde{\pi }_{1}\otimes
_{h}(K_{h}+M_{h})^{(r)}(0,T,x,y)+\sum_{r=0}^{\infty
}\sqrt{h}\widetilde{\pi }%
_{1}\otimes _{h}(H+M_{h,1}+\sqrt{h}N_{1})^{(r)}(0,T,x,y)
\label{eq:T501}
\end{equation}
and with $A=h\widetilde{\pi }_{2},B=H,C=(K_{h}-H)+M_{h}$ to
\begin{equation}
-\sum_{r=0}^{\infty }h\widetilde{\pi }_{2}\otimes
_{h}(K_{h}+M_{h})^{(r)}(0,T,x,y)+\sum_{r=0}^{\infty }h\widetilde{\pi }%
_{2}\otimes _{h}H^{(r)}(0,T,x,y).  \label{eq:T502}
\end{equation}
The estimate (\ref{eq:T500}) follows from (\ref{eq:30}),
(\ref{eq:T501}), (%
\ref{eq:T502}), Lemma 10 and Lemma 5 . \bigskip

\noindent \textit{Asymptotic treatment of the term $T_{6}$. }By
application of Lemma 9 we get that
\[
\left| T_{6}\right| \leq C(\varepsilon )hn^{-1/2+\varepsilon
}\zeta
_{\sqrt{T%
}}^{S}(y-x).
\]
\bigskip

\noindent \textit{Asymptotic treatment of the term $T_{7}$. } From
the recurrence
relation for $r=2,3,...$%
\begin{eqnarray*}&&
\widetilde{p}_{h}\otimes
_{h}(K_{h}+M_{h}+R_{h})^{(r)}(0,T,x,y)-\widetilde{p}%
_{h}\otimes _{h}H_{h}{}^{(r)}(0,T,x,y)
\\ &&
=\left[ \widetilde{p}_{h}\otimes
_{h}(K_{h}+M_{h}+R_{h})^{(r-1)}-\widetilde{p%
}_{h}\otimes _{h}H_{h}{}^{(r-1)}\right] \otimes _{h}H_{h}(0,T,x,y)
\\ && \qquad
+[\widetilde{p}_{h}\otimes _{h}(K_{h}+M_{h}+R_{h})^{(r-1)}\otimes
_{h}(K_{h}+M_{h}+R_{h}-H_{h})](0,T,x,y)
\end{eqnarray*}
and from Lemma 8 with $r=1$ we get by similar arguments as in the
proof of Lemma 9 that
\[
\left| T_{7}\right| \leq Ch^{3/2}T^{-1/2}\zeta _{\sqrt{T}%
}^{S}(y-x)=Chn^{-1/2}\zeta _{\sqrt{T}}^{S}(y-x).
\]\bigskip

\noindent \textit{Plugging in the asymptotic expansions of }
$T_{1},...,T_{7}$. We now plug the asymptotic expansions of
$T_{1},...,T_{7}$ into (\ref {eq:26}). Using Lemma 10, Theorem 2.1
in Konakov and Mammen (2002) we get \begin{eqnarray} && \nonumber
p_{h}(0,T,x,y)-p(0,T,x,y)
\\ && \nonumber
=\sqrt{h}\left[ \widetilde{\pi }_{1}+p^{d}\otimes _{h}\Re
_{1}\right] \otimes _{h}\Phi (0,T,x,y)
\\ && \nonumber \qquad
+h\Big\{ \left[ \widetilde{\pi }_{2}+\widetilde{\pi }_{1}\otimes
_{h}\Phi \otimes _{h}\Re _{1}+p^{d}\otimes _{h}\Re
_{2}+p^{d}\otimes _{h}\Re
_{3}%
\right] \otimes _{h}\Phi (0,T,x,y)
\\ && \nonumber \qquad
+p^{d}\otimes _{h}\left( \Re _{1}\otimes _{h}\Phi \right)
^{(2)}(0,T,x,y)
\\ \nonumber &&  \qquad
 +\frac{1}{2}p\otimes _{h}(L_{\star
}^{2}-L^{2})p^{d}(0,T,x,y)-\frac{1}{%
2}p\otimes _{h}(L^{\prime }-\widetilde{L}^{\prime
})p^{d}(0,T,x,y)\Big\}\\ &&  \qquad +O(h^{1+\delta }\zeta
_{\sqrt{T}}(y-x)), \label{eq:P01}
\end{eqnarray}
where
\begin{eqnarray*}
p^d(ih,i^\prime h,x,y) &=& \sum_{r=0}^\infty \widetilde p
\otimes_h
H^{(r)}(ih,i^\prime h,x,y),\\
\Re
_{1}(s,t,x,y)&=&N_{1}(s,t,x,y)+M_{1}(s,t,x,y)-\widetilde{M}_{1}(s,t,x,y),
\\
\Re _{2}(s,t,x,y)&=&N_{2}(s,t,x,y)+\Pi _{1}(s,t,x,y)-\widetilde{\Pi }%
_{1}(s,t,x,y),
\\
\Re _{3}(s,t,x,y)&=&\sum_{\left| \nu \right| =4}\frac{\chi _{\nu
}(s,x)-\chi _{\nu }(s,y)}{\nu !}D_{x}^{\nu }\widetilde{p}(s,t,x,y),
\\
M_{1}(s,t,x,y)&=&\sum_{\left| \nu \right| =3}\frac{\chi _{\nu
}(s,x)}{\nu
!}%
D_{x}^{\nu
}\widetilde{p}(s,t,x,y),
\\ \widetilde{M}_{1}(s,t,x,y)&=&\sum_{\left|
\nu \right| =3}\frac{\chi _{\nu }(s,y)}{\nu !}D_{x}^{\nu
}\widetilde{p}%
(s,t,x,y),
\\
\Pi _{1}(s,t,x,y)&=&\sum_{\left| \nu \right| =3}\frac{\chi _{\nu
}(s,x)}{\nu !}%
D_{x}^{\nu }\widetilde{\pi }_{1}(s,t,x,y),\\
 \widetilde{\Pi } _{1}(s,t,x,y)&=&\sum_{\left| \nu \right|
=3}\frac{\chi _{\nu }(s,y)}{\nu
!}%
D_{x}^{\nu }\widetilde{\pi }_{1}(s,t,x,y).
\end{eqnarray*}
Note that for the homogenous case and $T=[0,1]$ (\ref{eq:P01})
coincides with formula (53) on page 623 in Konakov and Mammen
(2005).
\bigskip

\noindent \textit{Asymptotic replacement of }$p^{d}$ by $p$. It
follows from (\ref{MR-008}), (\ref{MR-017}) and (\ref{MR-018}) that
\begin{equation}
\left| (p^{d}-p)(ih,jh,x,z)\right| \leq C(\varepsilon
)h^{1-\varepsilon }(jh-ih)^{\varepsilon -1/2}\phi
_{\sqrt{(j-i)h}}(z-x)  \label{eq:P02}
\end{equation}
for any $0<\varepsilon <1/2$. Using (\ref{eq:P02}) and making an
integration by parts we can replace $p^{d}$ by $p$ in \
(\ref{eq:P01}%
). For example the operator $L_{\star }^{2}-L^{2}$ is an operator
of order three. Applying integration by parts we get for $\left|
\nu
\right| =3$%
\begin{eqnarray*}
\left| \sum_{i=1}^{n-1}h\int D_{z}^{\nu
}p(0,ih,x,z)(p^{d}-p)(ih,T,z,y)dz\right| &\leq& C(\varepsilon
)h^{1-\varepsilon }\sum_{i=1}^{n-1}h\frac{1}{(ih)^{3/2}}\frac{1}{%
(T-ih)^{1/2-\varepsilon }}\phi _{\sqrt{T}}(y-x) \\
&\leq& C(\varepsilon )h^{1/2-2\varepsilon }T^{2\varepsilon
-1/2}B(\varepsilon ,\varepsilon +\frac{1}{2})\phi _{\sqrt{T}}(y-x).
\end{eqnarray*}
By (B2) we have $0<\varkappa <1-4\varepsilon .$ This implies
\begin{eqnarray*}
\left| \frac{h}{2}p\otimes _{h}(L_{\star
}^{2}-L^{2})(p^{d}-p)(0,T,x,y)\right| &\leq& C(\varepsilon
)hT^{1/2}n^{-(1/2-2\varepsilon -\varkappa /2)}\phi
_{\sqrt{T}}(y-x)
\\
&\leq& C(\varepsilon )h^{1+\delta }\phi _{\sqrt{T}}(y-x)
\end{eqnarray*}
for some $0<\delta <1/2.$ The other terms in (\ref{eq:P01})
containing $%
p^{d}$ can be estimated analogously. Thus we get the following
representation
\begin{eqnarray*}
&& p_{h}(0,T,x,y)-p(0,T,x,y)
\\ &&
=\sqrt{h}\left[ \widetilde{\pi }_{1}+p\otimes _{h}\Re _{1}\right]
\otimes _{h}\Phi (0,T,x,y)
\\ && \qquad
+h\Big\{ \left[ \widetilde{\pi }_{2}+\widetilde{\pi }_{1}\otimes
_{h}\Phi \otimes _{h}\Re _{1}+p\otimes _{h}\Re _{2}+p\otimes _{h}\Re
_{3}\right] \otimes _{h}\Phi (0,T,x,y)
\\ && \qquad
+p\otimes _{h}\left( \Re _{1}\otimes _{h}\Phi \right)
^{(2)}(0,T,x,y)
\\ && \qquad
 +\frac{1}{2}p\otimes _{h}(L_{\star
}^{2}-L^{2})p(0,T,x,y)-\frac{1}{2}%
p\otimes _{h}(L^{\prime }-\widetilde{L}^{\prime })p(0,T,x,y)\Big\}
+O(h^{1+\delta }\zeta _{\sqrt{T}}(y-x)).
\end{eqnarray*}
In the further analysis we make use of the following binary
operation $\otimes _{h}^{\prime }$. This operator generalizes the
binary operation
$%
\otimes $ introduced in Konakov and Mammen (2005). For $s\in
[0,t-h]$ and $t\in \{h,2h,...,T\}$ the operation $\otimes
_{h}^{\prime }$ is defined as follows
\[
f\otimes _{h}^{\prime }g(s,t,x,y)=\sum_{s\leq jh\leq t-h}h\int
f(s,jh,x,z)g(jh,t,z,y)dz.
\]
Note that for $s\in \{0,h,2h,...,T\}$ the two operations $\otimes
_{h}^{\prime }$ and $\otimes _{h}$ coincide.
\bigskip

\noindent \textit{Asymptotic replacement of }$(p\otimes _{h}\Re
_{i})\otimes _{h}\Phi (0,T,x,y)$ \textit{by} $p\otimes (\Re
_{i}\otimes _{h}^{\prime }\Phi )(0,T,x,y)=(p\otimes \Re
_{i})\otimes _{h}\Phi (0,T,x,y),i=1,2,3,$ \textit{}$[p\otimes
_{h}\left( \Re _{1}\otimes _{h}\Phi \right) ]\otimes _{h}\left(
\Re _{1}\otimes _{h}\Phi \right) (0,T,x,y)$ \textit{by }$p\otimes
\lbrack (\Re _{i}\otimes _{h}^{\prime }\Phi )\otimes _{h}^{\prime
}(\Re _{i}\otimes _{h}^{\prime }\Phi )](0,T,x,y), $ $p\otimes
_{h}(L_{\star }^{2}-L^{2})p(0,T,x,y)$ \textit{by }$p\otimes
(L_{\star }^{2}-L^{2})p(0,T,x,y)$ \textit{and} $p\otimes
_{h}(L^{\prime }-\widetilde{L}^{\prime })p(0,T,x,y)$ \textit{by }
$%
p\otimes (L^{\prime }-\widetilde{L}^{\prime })p(0,T,x,y).$

These replacements follow from the definitions of $\Re
_{i},i=1,2,3,$ and can be proved by the same method as in the
treatment of $T_1$. There an estimate for the replacement error of
$p\otimes _{h}H$ by $p\otimes H$ is given. Linearity of the
operation $\otimes _{h}$ implies that it is enough to consider the
functions $p\otimes _{h}\Im $ where $\Im (u,t,z,v)$ is a function
that has one of the following forms:
\begin{eqnarray*}
&&\frac{\chi _{\nu }(u,z)-\chi _{\nu }(u,v)}{\nu !}D_{x}^{\nu
}\widetilde{p}%
(u,t,z,v) \mbox{ with }\left| \nu \right| =3,4,\\ && \frac{\chi
_{\nu }(u,z)-\chi _{\nu
}(u,v)}{%
\nu !}D_{x}^{\nu }\widetilde{\pi }_{1}(u,t,z,v) \mbox{ with } \left|
\nu \right| =3,
\\ &&
(L-\widetilde{L})\widetilde{\pi
}_{1}(u,t,z,v) \mbox{ or } (L-\widetilde{L})\widetilde{%
\pi }_{2}(u,t,z,v). \end{eqnarray*}
We consider the case $\Im (u,t,z,v)=(L-\widetilde{L})\widetilde{\pi }%
_{1}(u,t,z,v)$. The other cases can be treated similarly. It is
enough to consider a typical term of
$(L-\widetilde{L})\widetilde{\pi }_{1}(u,t,z,v)$. We will give
bounds for
\begin{eqnarray} \nonumber
&& \int_{0}^{jh}du\int p(0,u,x,z)\left( \int_{u}^{jh}\chi _{\nu
}(w,v)dw\right) D_{z}^{\nu
}(L-\widetilde{L})\widetilde{p}(u,jh,z,v)dz
\\ && \qquad \nonumber
-\sum_{i=0}^{j-1}h\int p(0,ih,x,z)\left( \int_{ih}^{jh}\chi _{\nu
}(w,v)dw\right) D_{z}^{\nu
}(L-\widetilde{L})\widetilde{p}(ih,jh,z,v)dz
\\ &&= \nonumber
\sum_{i=0}^{j-1}\int_{ih}^{(i+1)h}du\int [\lambda (u)-\lambda
(ih)]dz\\ &&= \nonumber
\sum_{i=0}^{j-1}\int_{ih}^{(i+1)h}(u-ih)du\int \lambda ^{\prime
}(ih)dz
\\ && \qquad
+\sum_{i=0}^{j-1}\int_{ih}^{(i+1)h}(u-ih)^{2}du\int_{0}^{1}(1-%
\delta )\int \lambda ^{\prime \prime }(s)\mid _{s=s_{i}}dzd\delta
du \label{eq:P04aa}
\end{eqnarray}
where $\lambda (u)=p(0,u,x,z)\left( \int_{u}^{jh}\chi _{\nu
}(w,v)dw\right) D_{z}^{\nu }H(u,jh,z,v),s_{i}=ih+\delta (u-ih).$ As
in the treatment of $T_1$, we obtain that
\begin{eqnarray} \nonumber
&& \sum_{i=0}^{j-1}\int_{ih}^{(i+1)h}(u-ih)du\int \lambda ^{\prime
}(ih)dz
\\ \nonumber
&& =\frac{h}{2}\sum_{i=0}^{j-1}h\int_{ih}^{jh}\chi _{\nu
}(s,v)ds\int p(0,ih,x,z)D_{z}^{\nu }A_{0}(ih,jh,z,v)dz \\
\nonumber
&&
 \qquad +\frac{h}{2}\sum_{i=0}^{j-1}h\int_{ih}^{jh}\chi _{\nu }(s,v)ds\int
p(0,ih,x,z)D_{z}^{\nu }H_{1}(ih,jh,z,v)dz
 \\
 &&
 \qquad
-\frac{h}{2}\sum_{i=0}^{j-1}h\chi _{\nu }(ih,v)\int
p(0,ih,x,z)D_{z}^{\nu }H(ih,jh,z,v)dz,  \label{eq:P04a}
\end{eqnarray}
where
\begin{eqnarray*}
A_{0}(s,jh,z,v)&=&(L^{2}-2L\widetilde{L}+\widetilde{L}^{2})\widetilde{p}%
(s,jh,z,v),
\\
H_{l}(s,t,z,v)&=&\frac{1}{2}\sum_{i,j=1}^{d}\left( \frac{\partial
^{l}\sigma _{ij}(s,z)}{\partial s^{l}}-\frac{\partial ^{l}\sigma
_{ij}(s,v)}{\partial s^{l}}\right) \frac{\partial
^{2}\widetilde{p}(s,t,z,v)}{\partial z_{i}\partial z_{j}}
\\ && \qquad
\sum_{i=1}^{d}\left( \frac{\partial ^{l}m_{i}(s,z)}{\partial
s^{l}}-\frac{%
\partial ^{l}m_{i}(s,v)}{\partial s^{l}}\right) \frac{\partial
\widetilde{p}%
(s,t,z,v)}{\partial z_{i}} \end{eqnarray*} for $l=0,1,2$ with
$H_{0}\equiv H$. The differential operator $A_{0}$ was introduced
before equation (\ref{MR-012a}). It is a fourth order differential
operator. From the structure of this operator and from
(\ref{eq:P04a}) it is clear that it is enough to estimate
\begin{equation}
I\triangleq \frac{h}{2}\sum_{j=0}^{n-1}h\int
\sum_{i=0}^{j-1}h\int_{ih}^{jh}\chi _{\nu }(s,v)ds\int
p(0,ih,x,z)D_{z}^{\nu +\mu }\widetilde{p}(ih,jh,z,v)dz\Phi
(jh,T,v,y)dv  \label{eq:P04b}
\end{equation}
for $\left| \nu \right| =3,\left| \mu \right| =3.$ To estimate
(\ref{eq:P04b}%
) we consider three possible cases: a) $jh>T/2,ih\leq jh/2$
$\Longrightarrow jh-ih>T/4$ b) $jh>T/2,ih>jh/2$ $\Longrightarrow $
$ih>T/4$ c)
$jh<T/2$ $%
\Longrightarrow $ $T-jh>T/2.$ In the case a) we apply integration by
parts. This transfers two derivatives to $p(0,ih,x,z)$. This gives
\begin{eqnarray*}
&&\left| \frac{h}{2}\sum_{i=0}^{j-1}h\int_{ih}^{jh}\chi _{\nu
}(s,v)ds\int D_{z}^{e_{k}+e_{l}}p(0,ih,x,z)D_{z}^{\nu +\mu
-e_{k}-e_{l}}\widetilde{p}%
(ih,jh,z,v)dz\right|
\\ &&
\leq Ch^{1-2\varepsilon }\int_{0}^{jh}\frac{1}{u^{1-\varepsilon
}}\frac{du}{%
(jh-u)^{1-\varepsilon }}\phi _{\sqrt{jh}}(v-x)
\\ &&
\leq C(\varepsilon )h^{1-2\varepsilon }(jh)^{2\varepsilon -1}\phi
_{\sqrt{jh}%
}(v-x)
\end{eqnarray*}
and
\begin{eqnarray} \nonumber
\left| I\right| &\leq& C(\varepsilon )h^{1-2\varepsilon
}\int_{0}^{T}\frac{du}{%
u^{1-2\varepsilon }(T-u)^{1/2}}\phi _{\sqrt{T}}(y-x) \\
\nonumber &\leq& C(\varepsilon )h^{3/4}T^{1/2}n^{-(1/4-2\varepsilon
-\varkappa /2)}\phi _{\sqrt{T}}(y-x)\\ &\leq& C(\varepsilon
)T^{1/2-\delta }h^{3/4+\delta }\phi
_{%
\sqrt{T}}(y-x),  \label{eq:P04c}
\end{eqnarray}
where $\delta =(1/4-2\varepsilon -\varkappa /2)>0$ if $\varkappa
<1/2-4\varepsilon ,$ $0<\varepsilon <0,05$ (see the condition (B2)).
In the case b) we apply integration by parts and transfer four
derivatives to $p(0,ih,x,z).$ This gives the same estimate as in
(\ref{eq:P04c}). At last, in the case c) we make an integration by
parts and transfer three derivatives to $\Phi (jh,T,v,y)$ and one
derivative to $p(0,ih,x,z).$ This
gives the same estimate as in (\ref{eq:P04c}). To pass from $D_{z}^{\mu }%
\widetilde{p}(ih,jh,z,v)$ to $D_{v}^{\mu }\widetilde{p}(ih,jh,z,v)$
we use the following estimate
\[
\left| D_{z}^{\mu }\widetilde{p}(ih,jh,z,v)+D_{v}^{\mu }\widetilde{p}%
(ih,jh,z,v)\right| \leq C\phi _{\sqrt{jh-ih}}(v-z).
\]
Clearly, the same estimate (\ref{eq:P04c}) holds true for the
other summands in the right hand side of (\ref{eq:P04a}). This
gives
\begin{eqnarray*}
&& \frac{h}{2}\left| \sum_{j=0}^{n-1}h\int
\sum_{i=0}^{j-1}h\int_{ih}^{jh}\chi _{\nu }(s,v)ds\int
p(0,ih,x,z)D_{z}^{\nu }H_{1}(ih,jh,z,v)dz\Phi (jh,T,v,y)dv\right|
\\ &&
\leq C(\varepsilon )T^{1/2-\delta }h^{3/4+\delta }\phi
_{\sqrt{T}}(y-x),
\\ &&
\frac{h}{2}\sum_{i=0}^{j-1}h\int \chi _{\nu
}(ih,v)\sum_{i=0}^{j-1}h\int p(0,ih,x,z)D_{z}^{\nu
}H(ih,jh,z,v)dz\Phi (jh,T,v,y)dv
\\&&
\leq C(\varepsilon )T^{1/2-\delta }h^{3/4+\delta }\phi
_{\sqrt{T}}(y-x).
\end{eqnarray*}
We now estimate the second summand in the right hand side of (\ref
{eq:P04aa}). Similarly as in (\ref{MR-013}) we obtain
\begin{eqnarray} \nonumber &&
\sum_{i=0}^{j-1}\int_{ih}^{(i+1)h}(u-ih)^{2}du\int_{0}^{1}(1-%
\delta )\int \lambda ^{\prime \prime }(s)\mid _{s=s_{i}}dzd\delta
du\\ \nonumber &&=
\sum_{i=0}^{j-1}\int_{ih}^{(i+1)h}(u-ih)^{2}\int_{0}^{1}(1-\delta
)\sum_{k=1}^{4}\int_{s}^{jh}\chi _{\nu }(\tau ,v)d\tau \int
p(0,s,x,z)D_{z}^{\nu }A_{k}(s,jh,z,v)\mid _{s=s_{i}}dzd\delta du
\\ \nonumber && \qquad
-\sum_{i=0}^{j-1}\int_{ih}^{(i+1)h}(u-ih)^{2}\int_{0}^{1}(1-\delta
)\chi _{\nu }(s,v)\int {}p(0,s,x,z)D_{z}^{\nu }A_{0}(s,jh,z,v)\mid
_{s=s_{i}}dzd\delta du
\\ \nonumber && \qquad
-\sum_{i=0}^{j-1}\int_{ih}^{(i+1)h}(u-ih)^{2}\int_{0}^{1}(1-\delta
)\chi _{\nu }(s,v)\int p(0,s,x,z)D_{z}^{\nu }H_{1}(s,jh,z,v)\mid
_{s=s_{i}}dzd\delta du
\\  && \qquad
-\sum_{i=0}^{j-1}\int_{ih}^{(i+1)h}(u-ih)^{2}\int_{0}^{1}(1-\delta
)\frac{%
\partial \chi _{\nu }(s,v)}{\partial s}\int p(0,s,x,z)D_{z}^{\nu
}H(s,jh,z,v)\mid _{s=s_{i}}dzd\delta du,  \label{eq:P04d}
\end{eqnarray}
where the operators $A_{i},i=1,2,3,4,$ are defined as follows:
\begin{eqnarray*}
A_{1}(s,jh,z,v)&=&(L^{3}-3L^{2}\widetilde{L}+3L\widetilde{L}^{2}-\widetilde{L}%
^{3})\widetilde{p}(s,jh,z,v), \\
A_{2}(s,jh,z,v)&=&(L_{1}H+2LH_{1})(s,jh,z,v),
\\
A_{3}(s,jh,z,v)&=&[(L-\widetilde{L})\widetilde{L}_{1}+2(L_{1}-\widetilde{L}%
_{1})\widetilde{L}]\widetilde{p}(s,jh,z,v),
\\
A_{4}(s,jh,z,v)&=&H_{2}(s,jh,z,v).
\end{eqnarray*}
The operator $A_{1}$ was introduced in (\ref{MR-014}).
For this operator it is enough to estimate for fixed $p,q,r,l$%
\[
\sum_{i=0}^{j-1}\int_{ih}^{(i+1)h}(u-ih)^{2}\int_{0}^{1}(1-\delta
)\sum_{k=1}^{4}\int_{s}^{jh}\chi _{\nu }(\tau ,v)d\tau \int
p(0,s,x,z)D_{z}^{\nu }\left( \frac{\partial
^{4}\widetilde{p}(s,jh,z,v)}{%
\partial z_{p}\partial z_{q}\partial z_{l}\partial z_{r}}\right) \mid
_{s=s_{i}}dzd\delta du.
\]
As in (\ref{MR-016}) we obtain that this term does not exceed
\begin{equation}
C(\varepsilon )h^{3/2-\varepsilon }(jh)^{2\varepsilon -1}\phi
_{\sqrt{jh}%
}(v-x).  \label{eq:P04f}
\end{equation}
It follows from the explicit form of these operators that the same
estimate (\ref{eq:P04f}) holds for $A_{2},A_{3}$ and $A_{4}.$ The
other three terms in the right hand side of (\ref{eq:P04d}) do not
contain the factor $\int_{s}^{jh}\chi _{\nu }(\tau ,v)d\tau $ and
they can be estimated separately. Clearly, it is enough to estimate
the term containing $A_{0}.$ The remaining two summands are less
singular. From the explicit form of $A_{0}$ (compare also
(\ref{MR-012})) we obtain that it is enough to estimate for fixed
$q,l,r$
$$
\sum_{i=0}^{j-1}\int_{ih}^{(i+1)h}(u-ih)^{2}\int_{0}^{1}(1-\delta
)\chi
_{\nu }(s,v)\int {}p(0,s,x,z)D_{z}^{\nu }\left( \frac{\partial ^{3}%
\widetilde{p}(s,jh,z,v)}{\partial z_{q}\partial z_{l}\partial
z_{r}}\right) (s,jh,z,v)\mid _{s=s_{i}}dzd\delta du.
$$
Analogously to (\ref{MR-012a}) we get that this term does not exceed
\begin{equation}
C(\varepsilon )h^{1-2\varepsilon }(jh)^{2\varepsilon -1}\phi
_{\sqrt{jh}%
}(v-x).  \label{eq:P04h}
\end{equation}
Now from (\ref{eq:P04aa}), (\ref{eq:P04c}), (\ref{eq:P04d}),
(\ref{eq:P04f}) and (\ref{eq:P04h}) we obtain that
\[
\left| \lbrack p\otimes _{h}(L-\widetilde{L})\widetilde{\pi
}_{1}]\otimes
_{h}\Phi (0,T,x,y)-p\otimes \lbrack (L-\widetilde{L})\widetilde{\pi }%
_{1}\otimes _{h}^{\prime }\Phi ](0,T,x,y)\right|
\]
\begin{equation}
\leq Ch^{3/4+\delta }\phi _{\sqrt{T}}(y-x)  \label{eq:P04i}
\end{equation}
for some $\delta >0.$ The other replacements can be shown
analogously. Thus we come to the following representation
\begin{eqnarray} \nonumber &&
p_{h}(0,T,x,y)-p(0,T,x,y)
\\ \nonumber &&
=\sqrt{h}\left[ \widetilde{\pi }_{1}\otimes _{h}^{\prime }\Phi
(0,T,x,y)+p\otimes (\Re _{1}\otimes _{h}^{\prime }\Phi
)(0,T,x,y)\right]
\\ \nonumber && \qquad
+h\left[ \widetilde{\pi }_{2}\otimes _{h}^{\prime }\Phi
(0,T,x,y)+p\otimes (\Re _{2}\otimes _{h}^{\prime }\Phi
)(0,T,x,y)+p\otimes _{h}(\Re _{3}\otimes _{h}^{\prime }\Phi
)(0,T,x,y)\right]
\\ \nonumber && \qquad
+h\left[ \widetilde{\pi }_{1}\otimes _{h}^{\prime }\Phi +p\otimes
\left( \Re _{1}\otimes _{h}^{\prime }\Phi \right) \right] \otimes
_{h}^{\prime }(\Re _{1}\otimes _{h}^{\prime }\Phi )(0,T,x,y)
\\  && \qquad
+\frac{h}{2}p\otimes (L_{\star
}^{2}-L^{2})p(0,T,x,y)-\frac{h}{2}p\otimes (L^{\prime
}-\widetilde{L}^{\prime })p(0,T,x,y)+O(h^{1+\delta }\zeta
_{\sqrt{%
T}}(y-x)).  \label{eq:P05}
\end{eqnarray}
We now further simplify our expansion of $p_{h}-p.$ We start by
showing the following expansion
\begin{eqnarray} \nonumber &&
p_{h}(0,T,x,y)-p(0,T,x,y)
\\ \nonumber &&
=\sqrt{h}(p\otimes \mathcal{F}_{1}[p_{\Delta }])(0,T,x,y)+h\left(
p\otimes \mathcal{F}_{2}[p_{\Delta }]\right) (0,T,x,y)
\\ \nonumber && \qquad
+h\left( p\otimes \mathcal{F}_{1}[p\otimes \mathcal{F}_{1}[p_{\Delta
}]]\right) (0,T,x,y)
\\ && \qquad
+\frac{h}{2}p\otimes (L_{\star
}^{2}-L^{2})p(0,T,x,y)-\frac{h}{2}p\otimes (L^{\prime
}-\widetilde{L}^{\prime })p(0,T,x,y)+O(h^{1+\delta }\zeta
_{\sqrt{%
T}}(y-x)),  \label{eq:P06}
\end{eqnarray}
where for $s\in \lbrack 0,t-h],t\in \{h,2h,...,T\}$%
\begin{eqnarray*}
p_{\Delta }(s,t,z,y)&=&(\widetilde{p}\otimes _{h}^{\prime }\Phi
)(s,t,z,y)
\\
&=&\widetilde{p}(s,t,z,y)+\sum_{s\leq jh\leq t-h}h\int \widetilde{p}%
(s,jh,z,v)\Phi _{1}(jh,t,v,y)dv.
\end{eqnarray*}
Here $\Phi _{1}=H+H\otimes _{h}^{\prime }H+H\otimes _{h}^{\prime
}H\otimes _{h}^{\prime }H+...$. We now treat the term $\ p\otimes
\widetilde{%
L}\widetilde{\pi }_{1}(s,t,x,y).$%
\begin{eqnarray} \nonumber
p\otimes \widetilde{L}\widetilde{\pi
}_{1}(s,t,x,y)&=&\int_{s}^{t}d\tau \int p(s,\tau ,x,v)(t-\tau
)\sum_{\left| \nu \right| =3}\frac{\overline{\chi
}%
_{\nu }(\tau ,t,y)}{\nu !}D_{v}^{\nu
}(\widetilde{L}_{v}\widetilde{p}(\tau ,t,v,y))dv \\
 \nonumber &=&-\sum_{\left| \nu
\right| =3}\frac{1}{\nu !}\int dv\left[ \int_{s}^{t}p(s,\tau
,x,v)\left( \int_{\tau }^{t}\chi _{\nu }(u,y)du\right)
\frac{\partial }{\partial \tau }D_{v}^{\nu }\widetilde{p}(\tau
,t,v,y)d\tau %
\right]
\\
 \nonumber &=&
 -\sum_{\left| \nu \right| =3}\frac{1}{\nu !}\int
dv\int_{s}^{\frac{s+t}{2}%
}...-\sum_{\left| \nu \right| =3}\frac{1}{\nu !}\int
dv\int_{\frac{s+t}{2}%
}^{t}...
\\
  &=&I+II.  \label{eq:P07}
\end{eqnarray}
By integrating by parts w.r.t.\ the time variable we obtain for $I$.
\begin{eqnarray} \nonumber
I&=&-\sum_{\left| \nu \right| =3}\frac{1}{\nu !}\int dv\left[
p(s,\tau
,x,v)\left( \int_{\tau }^{t}\chi _{\nu }(u,y)du\right) D_{v}^{\nu }%
\widetilde{p}(\tau ,t,v,y)\mid _{\tau =s}^{\tau =(s+t)/2}\right.
\\ \nonumber && \qquad
\left. -\int_{s}^{\frac{s+t}{2}}D_{v}^{\nu }\widetilde{p}(\tau
,t,v,y)\left( \frac{\partial p(s,\tau ,x,v)}{\partial \tau
}\int_{\tau }^{t}\chi _{\nu }(u,y)du-p(s,\tau ,x,v)\chi _{\nu
}(\tau ,y)\right) d\tau \right]
\\ \nonumber &=&-\sum_{\left| \nu \right| =3}\frac{1}{\nu !}\int dv\left[
p(s,\frac{s+t}{2}%
,x,v)\left( \int_{\frac{s+t}{2}}^{t}\chi _{\nu }(u,y)du\right)
D_{v}^{\nu }%
\widetilde{p}(\frac{s+t}{2},t,v,y)\right.
\\ \nonumber && \qquad
\left. +\sum_{\left| \nu \right| =3}\frac{1}{\nu !}\left(
\int_{s}^{t}\chi _{\nu }(u,y)du\right) D_{v}^{\nu
}\widetilde{p}(s,t,x,y)\right]\\ \nonumber && \qquad  +\sum_{\left|
\nu \right| =3}\frac{1}{\nu !}\int_{s}^{\frac{s+t}{2}}d\tau \left(
\int_{\tau }^{t}\chi _{\nu }(u,y)du\right)
 \int L^{T}p(s,\tau ,x,v)D_{v}^{\nu }\widetilde{p}(\tau
,t,v,y)dv\\ && \qquad -\sum_{\left| \nu \right| =3}\frac{1}{\nu
!}\int_{s}^{\frac{s+t}{2}%
}\chi _{\nu }(\tau ,y)d\tau \int p(s,\tau ,x,v)D_{v}^{\nu
}\widetilde{p}%
(\tau ,t,v,y)dv.  \label{eq:P08}
\end{eqnarray}
For the second term we get \begin{eqnarray} \nonumber
II&=&\sum_{\left| \nu \right| =3}\frac{1}{\nu !}\int p(s,\frac{s+t}{2}%
,x,v)\left( \int_{\frac{s+t}{2}}^{t}\chi _{\nu }(u,y)du\right)
D_{v}^{\nu }%
\widetilde{p}(\frac{s+t}{2},t,v,y)dv
\\ && \qquad \nonumber
+\sum_{\left| \nu \right| =3}\frac{1}{\nu
!}\int_{\frac{s+t}{2}}^{t}d\tau \left( \int_{\tau }^{t}\chi _{\nu
}(u,y)du\right) \int L^{T}p(s,\tau ,x,v)D_{v}^{\nu
}\widetilde{p}(\tau ,t,v,y)dv
\\ && \qquad
-\sum_{\left| \nu \right| =3}\frac{1}{\nu
!}\int_{\frac{s+t}{2}}^{t}\chi _{\nu }(\tau ,y)d\tau \int p(s,\tau
,x,v)D_{v}^{\nu }\widetilde{p}(\tau ,t,v,y)dv.  \label{eq:P09}
\end{eqnarray}
From (\ref{eq:P07})- (\ref{eq:P09}) we have
$$
p\otimes \widetilde{L}\widetilde{\pi }_{1}(s,t,x,y)=\widetilde{\pi }%
_{1}(s,t,x,y)+p\otimes L\widetilde{\pi }_{1}(s,t,x,y)-p\otimes
\widetilde{M}%
_{1}(s,t,x,y). $$ This shows that
\begin{eqnarray}
\nonumber &&\widetilde{\pi }_{1}(s,t,x,y)+p\otimes \Re _{1}(s,t,x,y)
=\widetilde{\pi }_{1}(s,t,x,y)\\ \nonumber && \qquad +p\otimes
L\widetilde{\pi
}_{1}(s,t,x,y)-p%
\otimes \widetilde{L}\widetilde{\pi }_{1}(s,t,x,y)+p\otimes
M_{1}(s,t,x,y) -p\otimes \widetilde{M}_{1}(s,t,x,y)\\  && =p\otimes
M_{1}(s,t,x,y). \label{eq:P11}
\end{eqnarray}
It follows from (\ref{eq:P11}) and the definitions of the operations
$%
\otimes $ and $\otimes _{h}^{\prime }$ that
\begin{eqnarray} \nonumber
&& \sqrt{h}\left[ \widetilde{\pi }_{1}\otimes _{h}^{\prime }\Phi
(s,t,x,y)+(p\otimes \Re _{1})\otimes _{h}^{\prime }\Phi
(s,t,x,y)\right] \\
\nonumber && =\sqrt{h}(\widetilde{\pi }_{1}+p\otimes \Re
_{1})\otimes _{h}^{\prime }\Phi (s,t,x,y)\\
\nonumber &&=\sqrt{h}(p\otimes M_{1})\otimes _{h}^{\prime }\Phi
(s,t,x,y)
\\
\nonumber && =\sqrt{h}\sum_{0\leq jh\leq t-h}h\int (p\otimes
M_{1})(s,jh,x,z)\Phi (jh,t,z,y)dz
\\
\nonumber && =\sqrt{h}\sum_{0\leq jh\leq t-h}h\int \left[
\int_{s}^{jh}du\int p(s,u,x,v)
M_{1}(u,jh,v,z)dv\right] \Phi (jh,t,z,y)dz\\
\nonumber && =\sqrt{h}\sum_{0\leq jh\leq t-h}h\int \left[
\int_{s}^{t}du\chi \lbrack s,jh]\int p(s,u,x,v)
M_{1}(u,jh,v,z)dv\right] \Phi (jh,t,z,y)dz
\\
\nonumber &&=\sqrt{h}%
\int_{s}^{t}du\int p(s,u,x,v)\sum_{\left| \nu \right| =3}\frac{\chi
_{\nu }(u,v)}{\nu !} \times D_{v}^{\nu }\left[ \sum_{0\leq jh\leq
t-h}h\chi \lbrack s,jh]\int
\widetilde{p}(u,jh,v,z)\Phi (jh,t,z,y)dz\right] dv
\\
\nonumber &&
=\sqrt{h}%
\int_{s}^{t}du\int p(s,u,x,v) \times \sum_{\left| \nu \right|
=3}\frac{\chi _{\nu }(u,v)}{\nu !}D_{v}^{\nu }p_{\Delta
}(u,t,v,y)dv
\\
 &&=\sqrt{h}(p\otimes \mathcal{F}_{1})[p_{\Delta
}](s,t,x,y).  \label{eq:P12}
\end{eqnarray}
Here, $\chi \lbrack s,jh]$ denotes the indicator of the interval
$[s,jh]$. Using similar arguments as in the proof of (\ref{eq:P12})
one can show that
\begin{eqnarray} \nonumber && h\left[
\widetilde{\pi }_{2}\otimes _{h}^{\prime }\Phi (s,t,x,y)+(p\otimes
\Re _{2})\otimes _{h}^{\prime }\Phi (s,t,x,y)\right. \left.
+p\otimes _{h}(\Re _{3}\otimes _{h}^{\prime }\Phi )(s,t,x,y)\right]
\\ &&=h(p\otimes \mathcal{F}_{2})[p_{\Delta }](s,t,x,y)+hp\otimes \Pi
_{1}\otimes _{h}^{\prime }\Phi (s,t,x,y). \label{eq:P13}
\end{eqnarray}
For the first two terms in the right hand side of (\ref{eq:P05}) we
obtain from (\ref{eq:P12}) and (\ref{eq:P13})
\begin{eqnarray} \nonumber &&
\sqrt{h}\left[ \widetilde{\pi }_{1}\otimes _{h}^{\prime }\Phi
(0,T,x,y)+(p\otimes \Re _{1})\otimes _{h}^{\prime }\Phi
(0,T,x,y)\right] \\ \nonumber &&  \qquad  +h\left[ \widetilde{\pi
}_{2}\otimes _{h}^{\prime }\Phi (0,T,x,y)+p\otimes (\Re _{2}\otimes
_{h}^{\prime }\Phi )(0,T,x,y)+p\otimes _{h}(\Re _{3}\otimes
_{h}^{\prime }\Phi )(0,T,x,y)\right]
\\ &&
=\sqrt{h}(p\otimes \mathcal{F}_{1})[p_{\Delta
}](0,T,x,y)+h(p\otimes \mathcal{F}_{2})[p_{\Delta
}](s,t,x,y)+hp\otimes \Pi _{1}\otimes _{h}^{\prime }\Phi
(s,t,x,y).  \label{eq:P14}
\end{eqnarray}
Using (\ref{eq:P12}) we get
\begin{eqnarray*} &&
h\left[ \widetilde{\pi }_{1}\otimes _{h}^{\prime }\Phi +p\otimes
\left( \Re _{1}\otimes _{h}^{\prime }\Phi \right) \right] \otimes
_{h}^{\prime }(\Re _{1}\otimes _{h}^{\prime }\Phi )(0,T,x,y)
\\ &&
=h(p\otimes \mathcal{F}_{1}[p_{\Delta }])\otimes _{h}^{\prime
}(\Re _{1}\otimes _{h}^{\prime }\Phi )(0,T,x,y)
\\ &&
=hp\otimes \mathcal{F}_{1}\left[ p_{\Delta }\otimes _{h}^{\prime
}(\Re _{1}\otimes _{h}^{\prime }\Phi )\right] (0,T,x,y).
\end{eqnarray*}
Note that
\begin{eqnarray*}
hp\otimes \Pi _{1}\otimes _{h}^{\prime }\Phi
(s,t,x,y)&=&h\int_{s}^{t}du\int
p(s,u,x,v)\sum_{\left| \nu \right| =3}\frac{\chi _{\nu }(u,v)}{\nu !}%
D_{v}^{\nu }[\widetilde{\pi }_{1}\otimes _{h}^{\prime }\Phi
](u,t,v,y) \\ &=&hp\otimes \mathcal{F}_{1}[\widetilde{\pi
}_{1}\otimes _{h}^{\prime }\Phi ](s,t,x,y).
\end{eqnarray*}
For the proof of (\ref{eq:P06}) it remains to show that
\begin{eqnarray} \nonumber &&
hp\otimes \mathcal{F}_{1}[\widetilde{\pi }_{1}\otimes _{h}^{\prime
}\Phi +p_{\Delta }\otimes _{h}^{\prime }(\Re _{1}\otimes
_{h}^{\prime }\Phi )](0,T,x,y)
\\ &&
=h\left( p\otimes \mathcal{F}_{1}[p\otimes
\mathcal{F}_{1}[p_{\Delta }]]\right) (0,T,x,y)+O(h^{1+\delta
}\zeta _{\sqrt{T}}(y-x)). \label{eq:P15}
\end{eqnarray}
We will show that
\begin{eqnarray}
&&hp\otimes \mathcal{F}_{1}[(p-p_{\Delta })\otimes _{h}^{\prime
}(\Re _{1}\otimes _{h}^{\prime }\Phi )](0,T,x,y)=O(h^{1+\delta
}\zeta
_{\sqrt{T}%
}(y-x)),  \label{eq:P16}
\\ && \nonumber
hp\otimes \mathcal{F}_{1}[p\otimes _{h}^{\prime }(\Re _{1}\otimes
_{h}^{\prime }\Phi )](0,T,x,y) -hp\otimes \mathcal{F}_{1}[p\otimes
(\Re _{1}\otimes _{h}^{\prime }\Phi )](0,T,x,y)\\ && \qquad
=O(h^{1+\delta }\zeta _{\sqrt{T}}(y-x)). \label{eq:P17}
\end{eqnarray}
Claim (\ref{eq:P15}) follows from (\ref{eq:P16}), (\ref{eq:P17}) and
(\ref {eq:P12}). The estimate (\ref{eq:P17}) can be shown similarly
as in the proof of (\ref {eq:P04i}). An additional singularity
arising from the derivatives in the operator $\mathcal{F}_{1}[\cdot
]$ can be treated by using the additional factor $h$
in (%
\ref{eq:P17}). To estimate (\ref{eq:P16}) note that from the
definition of
$\Re _{1}$ and $\Phi $%
\begin{equation}
\left| (\Re _{1}\otimes _{h}^{\prime }\Phi )(jh,T,x,y)\right| \leq
C(\varepsilon )h^{-\varepsilon }(T-jh)^{\varepsilon -1}\phi
_{\sqrt{T-jh}}. \label{eq:P18}
\end{equation}
Then we use the following estimate which can be proved by the same
method as in the treatment of $T_1$, where an estimate for
$(p-p^{d})(0,jh,x,y)$ was obtained.
\begin{equation}
\left| (p-p_{\Delta })(u,jh,v,z)\right| \leq Ch^{1/2}\phi
_{\sqrt{jh-u}%
}(z-v).  \label{eq:P19}
\end{equation}
From (\ref{eq:P18}) and (\ref{eq:P19})
\begin{equation}
\left| (p-p_{\Delta })\otimes _{h}^{\prime }(\Re _{1}\otimes
_{h}^{\prime }\Phi )(u,T,v,y)\right| \leq C(\varepsilon
)h^{1/2-\varepsilon }(T-u)^{\varepsilon }\phi _{\sqrt{T-u}}(y-v)
\label{eq:P20}
\end{equation}
For an estimate (\ref{eq:P16}) it is enough to estimate a typical
summand of the sum of the detailed representation of the left hand
side of (\ref{eq:P16}). E.g., for $\left| \nu \right| =3$ \ we have
to estimate
\begin{eqnarray*}
&&h\int_{0}^{T}du\int p(0,u,x,v)\frac{\chi _{\nu }(u,v)}{\nu
!}D_{v}^{\nu }[\sum_{\{j:u\leq jh\leq T-h\}}h\int (p-p_{\Delta
})(u,jh,v,z)
\\ && \qquad
\times (\Re _{1}\otimes _{h}^{\prime }\Phi )(jh,T,z,y)dz]dv\\
&&=h\int_{0}^{T/2}...+h\int_{T/2}^{T}...\\ &&=I+II.
\end{eqnarray*}
For an estimate of $II$ we apply integration by parts and transfer
three derivatives to $p(0,u,x,v)\frac{\chi _{\nu }(u,v)}{\nu !}$.
Using (\ref {eq:P20}) we obtain the following estimate
\begin{eqnarray}
\nonumber \left| II\right| &\leq& C(\varepsilon )h^{3/2-\varepsilon
}\int_{T/2}^{T}\frac{%
(T-u)^{\varepsilon }}{u^{3/2}}du\phi _{\sqrt{T}}(y-x)\\ &\leq&
C(\varepsilon )h^{3/2-\varepsilon }T^{\varepsilon }\phi
_{\sqrt{T}%
}(y-x).  \label{eq:P21}
\end{eqnarray}
For the treatment of $I$ we consider two cases: a) $jh-u\geq T/4$
and b) $jh-u\leq T/4$ $\Longrightarrow $ $T-jh\geq T/4$. Similarly
as in (\ref{MR-008}) the difference $h(p-p_{\Delta })$ can be
represented as
\begin{eqnarray} \nonumber &&
h(p-p_{\Delta })(u,jh,v,z)=h(p\otimes H-p\otimes _{h}^{\prime
}H)(u,jh,v,z)
\\ \nonumber && \qquad
+h(p\otimes H-p\otimes _{h}^{\prime }H)\otimes _{h}^{\prime }\Phi
_{1}(u,jh,v,z)
\\ \nonumber &&
=h\int_{u}^{j^{\star }h}d\tau \int p(u,\tau ,v,z^{\prime })H(\tau
,jh,z^{\prime },z)dz^{\prime }
\\ \nonumber && \qquad
+h\sum_{i=j^{\star }}^{j-1}\int_{ih}^{(i+1)h}d\tau \int (\lambda
(\tau ,z^{\prime })-\lambda (ih,z^{\prime }))dz^{\prime
}+h(p\otimes H-p\otimes _{h}^{\prime }H)\otimes _{h}^{\prime }\Phi
_{1}(u,jh,v,z)
\\ &&
=I^{\prime }+II^{\prime }+III^{\prime },  \label{eq:P21a}
\end{eqnarray}
where $\lambda (\tau ,z^{\prime })=p(u,\tau ,v,z^{\prime })H(\tau
,jh,z^{\prime },z),\Phi _{1}(ih,jh,z^{\prime }z)=H(ih,jh,z^{\prime
}z)+H\otimes _{h}^{\prime }H(ih,jh,z^{\prime }z)+...,j^{\star
}=j^{\star
}(u)=[%
\frac{u}{h}]+1$. Here $[x]$ is equal to the integral part for
noninteger $x$ and equal to $x-1$ for integer $x$. For $I^{\prime
},$ case a), we have $jh-\tau >T/5$ for $n$ large enough. With the
substitution $v+v^{\prime }=z^{\prime }$ we obtain
\begin{eqnarray} \nonumber &&
\left| D_{v}^{\nu }h\int_{u}^{j^{\star }h}d\tau \int p(u,\tau
,v,z^{\prime })H(\tau ,jh,z^{\prime },z)dz^{\prime }\right|
\\ \nonumber &&
=\left| D_{v}^{\nu }h\int_{u}^{j^{\star }h}d\tau \int p(u,\tau
,v,v+v^{\prime })H(\tau ,jh,v+v^{\prime },z)dv^{\prime }\right|
\\ \nonumber &&
\leq Ch\int_{u}^{j^{\star }h}\frac{d\tau }{(jh-\tau )^{2}}\phi
_{\sqrt{jh-u}%
}(z-v)\leq Ch^{2}T^{-2}\phi _{\sqrt{jh-u}}(z-v)
\\  &&
\leq Cn^{-2}\phi _{\sqrt{jh-u}}(z-v)=CT^{2}h^{2}\phi
_{\sqrt{jh-u}}(z-v) \label{eq:P22}
\end{eqnarray}
For the proof of (\ref{eq:P22}) we used the following estimate from
Friedman (1964) (Theorem 7, page 260) $$ \left| D_{v}^{\nu }p(u,\tau
,v,v+v^{\prime })\right| \leq C(\tau -u)^{-d/2}\exp \left(
\frac{C\left| v^{\prime }\right| }{\tau -u)}\right) .
$$
For $\int I^{\prime }(u,jh,v,z)(\Re _{1}\otimes _{h}^{\prime }\Phi
)(jh,T,z,y)dz,$ case b), it is enough to estimate for $\left| \nu
\right| =3$%
 \begin{eqnarray} \nonumber &&
h\int_{u}^{j^{\star }h}d\tau \int [p(u,\tau ,v,v+v^{\prime
})(\sigma _{lk}(\tau ,v+v^{\prime })-\sigma _{lk}(\tau
,z))]D_{v^{\prime }}^{\nu
}%
\frac{\partial ^{2}\widetilde{p}(\tau ,jh,v+v^{\prime
},z)}{\partial v_{l}^{\prime }\partial v_{k}^{\prime }}dv^{\prime
}
\\ && \qquad
\times (\Re _{1}\otimes _{h}^{\prime }\Phi )(jh,T,z,y)dz.
\label{eq:P23}
\end{eqnarray}
For an estimate of this term we transfer five derivatives from
$\widetilde{p}$ to $(\Re _{1}\otimes _{h}^{\prime }\Phi )(jh,T,z,y)$
and we use the following estimate for $\left|
\mu \right| =5$%
\begin{eqnarray*} &&
\left| D_{v^{\prime }}^{\mu }\widetilde{p}(\tau ,jh,v+v^{\prime
},z)+D_{z}^{\mu }\widetilde{p}(\tau ,jh,v+v^{\prime },z)\right|
\\ && \qquad
\leq C(jh-\tau )^{-d/2}\phi _{\sqrt{jh-\tau }}(z-v-v^{\prime }).
\end{eqnarray*}
We obtain that (\ref{eq:P23}) does not exceed \begin{eqnarray}
\nonumber C(j^{\star }h-\tau )hT^{-7/2}\phi
_{\sqrt{T-u}}(y-v)&\leq& Ch^{1+\delta }T^{1-\delta }n^{-(1-\delta
-7\varkappa /2)}\phi _{\sqrt{T-u}}(y-v) \\ &=&o(h^{1+\delta
}T^{1-\delta })\phi _{\sqrt{T-u}}(y-v). \label{eq:P24}
\end{eqnarray}
We used that for any $0<$ $\delta <1$ it holds that $\varkappa
<\frac{2-2\delta }{7}$, see condition (B2). For an estimate of $\int
II^{\prime }(u,jh,v,z)(\Re _{1}\otimes _{h}^{\prime }\Phi
)(jh,T,z,y)dz$ we use the decomposition (\ref{MR-009}). For getting
an estimate for the terms in $\ II^{\prime }$ that contain the first
derivatives $\lambda ^{\prime }(ih,z^{\prime })$ we use the identity
(\ref{MR-011}) and similar arguments as already used in the
estimation of $\int I^{\prime }(u,jh,v,z)(\Re _{1}\otimes
_{h}^{\prime }\Phi )(jh,T,z,y)dz$. The estimate for terms
in $%
II^{\prime }$ containing second derivatives $\lambda ^{^{\prime
\prime }}(ih,z^{\prime })$ follows from (\ref{MR-015}) and
(\ref{MR-015a}). Finally, for $III^{\prime }$ the same estimates
hold because of smoothing properties of the convolution ...$\otimes
_{h}^{^{\prime }}\Phi _{1}(u,jh,v,z).$ This implies (\ref{eq:P16})
and, hence, the expansion (\ref {eq:P06}).
\bigskip

\noindent \textit{Asymptotic replacement of }$p_{\Delta }$
\textit{by }$p$. Now, we compare $hp\otimes
\mathcal{F}_{2}[p_{\Delta }](0,T,x,y)$ with $hp\otimes
\mathcal{F}_{2}[p](0,T,x,y)$. Note that for $2\varkappa <\delta
<\frac{2%
}{5},\left| \nu \right| =4$%
\begin{eqnarray} \nonumber
&&\left| h\int_{0}^{h^{\delta }}du\int p(0,u,x,z)\chi _{\nu
}(u,z)D_{z}^{\nu }p(u,T,z,y)dz\right| \leq Ch^{1+\delta
}(T-h^{\delta })^{-2}\phi
_{\sqrt{T}%
}(y-x)
\\ && \qquad
\leq Ch^{1+\delta }\frac{n^{2\varkappa }}{(Tn^{\varkappa
}-n^{\varkappa }h^{\delta })^{2}}\leq Ch^{1+(\delta -2\varkappa
)}T^{2\varkappa }\phi
_{%
\sqrt{T}}(y-x),  \label{eq:P25}\\ && \left| h\int_{T-h^{\delta
}}^{T}du\int D_{z}^{\nu }[p(0,u,x,z)\chi _{\nu
}(u,z)]p(u,T,z,y)dz\right| \leq Ch^{1+(\delta -2\varkappa
)}T^{2\varkappa }\phi _{\sqrt{T}}(y-x).  \label{eq:P26}
\end{eqnarray}
The same estimates hold for $p_{\Delta }(u,T,z,y).$ Hence, it
suffices to consider $u\in \lbrack h^{\delta },T-h^{\delta }]$. We
now treat
\begin{eqnarray} \nonumber &&
h\int_{h^{\delta }}^{T-h^{\delta }}du\int p(0,u,x,z)\chi _{\nu
}(u,z)D_{z}^{\nu }(p-p_{\Delta })(u,T,z,y)dz
\\ &&
=h\int_{h^{\delta }}^{T/2}...+h\int_{T/2}^{T-h^{\delta }}...=I+II.
\label{eq:P27}
\end{eqnarray}
By using (\ref{eq:P19}) we get
\begin{eqnarray} \nonumber
\left| II\right| &=&\left| h\int_{T/2}^{T-h^{\delta }}du\int
D_{z}^{\nu }[p(0,u,x,z)\chi _{\nu }(u,z)](p-p_{\Delta
})(u,T,z,y)dz\right| \\ \nonumber &\leq& Ch^{3/2}n^{\varkappa }\phi
_{\sqrt{T}}(y-x)=Ch^{3/2-\varkappa }T^{\varkappa }\phi
_{\sqrt{T}}(y-x) \\&=&Ch^{1+\gamma }\phi _{\sqrt{T}}(y-x),\gamma >0
\label{eq:P28}
\end{eqnarray}
For $u\in \lbrack h^{\delta },T/2]$ it holds that
\begin{eqnarray} \nonumber
\left| I\right| &=&\left| h\int_{h^{\delta }}^{T/2}du\int D_{z}^{\nu
}[p(0,u,x,z)\chi _{\nu }(u,z)](p-p_{\Delta })(u,T,z,y)dz\right|
\\
&\leq& Ch^{3/2-\delta }\phi _{\sqrt{T}}(y-x).  \label{eq:P29}
\end{eqnarray}
Note that the condition $\delta <\frac{2}{5}$ implies that
$3/2-\delta
>1$. It follows from (\ref{eq:P25})-(\ref{eq:P29}) that
\begin{equation}
hp\otimes \mathcal{F}_{2}[p_{\Delta }](0,T,x,y)-hp\otimes \mathcal{F}%
_{2}[p](0,T,x,y)=O(h^{1+\gamma }\phi _{\sqrt{T}}(y-x)).
\label{eq:P30}
\end{equation}
For the proof of
\[
hp\otimes \mathcal{F}_{1}[p\otimes \mathcal{F}_{1}[p-p_{\Delta
}]]=O(h^{1+\delta }\phi _{\sqrt{T}}(y-x))
\]
we consider a typical summand for fixed $\nu ,\left| \nu \right| =3,$%
\begin{equation}
h\int_{0}^{T}du\int p(0,u,x,z)\chi _{\nu }(u,z)D_{z}^{\nu }\left[
\int_{u}^{T}d\tau \int p(u,\tau ,z,v)\chi _{\nu }(\tau
,v)D_{v}^{\nu }(p-p_{\Delta })(\tau ,T,v,y)dv\right] dz.
\label{eq:P31}
\end{equation}
As before it  suffices to consider the integrals for $u\in \lbrack
h^{\delta },T-h^{\delta }].$ The integral in (\ref{eq:P31}) is a sum
of four integrals
\begin{eqnarray} \nonumber
I_{1}&=&h\int_{h^{\delta }}^{T/2}du\int ...D_{z}^{\nu
}\int_{u}^{(T+u)/2}d\tau \int ...,
\\ \nonumber
I_{2}&=&h\int_{h^{\delta }}^{T/2}du\int ...D_{z}^{\nu
}\int_{(T+u)/2}^{T}d\tau \int ...,
\\ \nonumber
I_{3}&=&h\int_{T/2}^{T-h^{\delta }}du\int ...D_{z}^{\nu
}\int_{u}^{(T+u)/2}d\tau \int ...,
\\
I_{4}&=&h\int_{T/2}^{T-h^{\delta }}du\int ...D_{z}^{\nu
}\int_{(T+u)/2}^{T}d\tau \int ... .  \label{eq:P32}
\end{eqnarray}
Note that in the integrand in $I_{2}$ it holds that $\tau -u\geq
T/4$. By applying integration by parts w.r.t.\ $v$ and
(\ref{eq:P19}) we get
\begin{equation}
\left| I_{2}\right| \leq Ch^{3/2-\varkappa }T^{\varkappa }\phi
_{\sqrt{T}%
}(y-x)).  \label{eq:P33}
\end{equation}
Furthermore,  in the integrand in $I_{4}$ it holds that $u\geq
T/2,\tau -u\geq h^{\delta }/2,T-u\geq h^{\delta }$. Using
integration by parts w.r.t.\ $z$ we obtain
\begin{equation}
\left| I_{4}\right| \leq Ch^{3/2-\delta }T^{-1/2}\phi
_{\sqrt{T}}(y-x))\leq CT^{\varkappa /2}h^{3/2-\varkappa /2-\delta
}\phi _{\sqrt{T}}(y-x)), \label{eq:P34}
\end{equation}
where, by our choice of $\delta ,$ $3/2-\varkappa /2-\delta >1.$ For
an estimate of $I_{3}$ we use the representation
\begin{eqnarray} \nonumber &&
(p-p_{\Delta })(\tau ,T,v,y)=(p\otimes H-p\otimes _{h}^{\prime
}H)(\tau ,T,v,y)
\\ \nonumber && \qquad
+(p\otimes H-p\otimes _{h}^{\prime }H)\otimes _{h}^{\prime }\Phi
_{1}(\tau ,T,v,y)
\\ \nonumber &&
=\int_{\tau }^{j^{\star }h}ds\int p(\tau ,s,v,w)H(s,T,w,y)dw
\\ \nonumber && \qquad
+\frac{h}{2}[p\otimes _{h}^{\prime }(H_{1}+A_{0})](\tau ,T,v,y)
\\ \nonumber && \qquad
+\frac{1}{2}\sum_{i=j^{\star
}}^{n-1}\int_{ih}^{(i+1)h}(t-ih)^{2}\int_{0}^{1}(1-\gamma
)\sum_{k=1}^{4}\int p(\tau ,s,v,w)A_{k}(s,T,w,y)\mid _{s=ih+\gamma
(t-ih)}dwd\gamma dt,
\\ && \qquad
+(p\otimes H-p\otimes _{h}^{\prime }H)\otimes _{h}^{\prime }\Phi
_{1}(\tau ,T,v,y),  \label{eq:P35}
\end{eqnarray}
where $j^{\star }=j^{\star }(\tau )=[\tau /h]+1$. As above $[x]$
denotes the integer part for nonintegers $x$ and it is equal to
$x-1$ for integers $x$. The quantities $H_{1}$ and
$A_{k},k=0,1,2,3,4,$ have been defined (\ref{MR-014a}) and
\[
\Phi _{1}(ih,i^{\prime }h,z,z^{\prime })=H(ih,i^{\prime
}h,z,z^{\prime })+H\otimes _{h}^{\prime }H(ih,i^{\prime
}h,z,z^{\prime })+... .
\]
To estimate $D_{v}^{\nu }(p-p_{\Delta })(\tau ,T,v,y)$ we note that
\begin{eqnarray} \nonumber &&
h\left| D_{v}^{\nu }\int_{\tau }^{j^{\star }h}ds\int p(\tau
,s,v,w)H(s,T,w,y)dw\right|\\
 \nonumber && =h\left| D_{v}^{\nu
}\int_{\tau }^{j^{\star }h}ds\int p(\tau ,s,v,v+w^{\prime
})H(s,T,v+w^{\prime },y)dw^{\prime }\right|
\\ &&
\leq Ch\int_{\tau }^{j^{\star }h}\frac{ds}{(T-s)^{2}}\phi
_{\sqrt{T-\tau
}%
}(y-v)\leq Ch^{2-2\delta }\phi _{\sqrt{T-\tau }}(y-v).
\label{eq:P36}
\end{eqnarray}
Furthermore,
\begin{eqnarray} \nonumber &&
\left| \frac{h^{2}}{2}D_{v}^{\nu }[p\otimes _{h}^{\prime
}H_{1}](\tau ,T,v,y)\right| =\left| \frac{h^{2}}{2}D_{v}^{\nu
}\sum_{\tau \leq jh\leq T-h}h\int p(\tau
,jh,v,w)H_{1}(jh,T,w,y)dw\right|
\\
\nonumber && \leq \frac{h^{2}}{2}\left| \sum_{\tau \leq jh\leq
T-h^{\delta }/2}h\int D_{v}^{\nu }[p(\tau ,jh,v,v+w^{\prime
})H_{1}(jh,T,v+w^{\prime },y)]dw^{\prime }\right|
\\ \nonumber && \qquad
+\frac{h^{2}}{2}\left| C\sum_{i,k=1}^{d}\sum_{T-h^{\delta }/2<jh\leq
T-h}h\int D_{w^{\prime }}^{\nu +e_{i}+e_{k}}[p(\tau
,jh,v,v+w^{\prime })] \widetilde{p}(jh,T,v+w^{\prime },y)dw^{\prime
}\right|
\\ &&
\leq Ch^{2-2\delta }\phi _{\sqrt{T-\tau }}(y-v)+Ch^{2-5\delta
/2}\phi
_{%
\sqrt{T-\tau }}(y-v).  \label{eq:P37}
\end{eqnarray}
Because of the structure of the operator $A_{0}$ it is
enough to estimate for fixed $i,l,k$%
\begin{equation}
\frac{h^{2}}{2}D_{v}^{\nu }\sum_{\tau \leq jh\leq T-h}h\int
D_{w^{\prime
}}^{e_{k}}p(\tau ,jh,v,v+w^{\prime })\frac{\partial ^{2}\widetilde{p}%
(jh,T,v+w^{\prime },y)}{\partial w_{i}^{\prime }\partial
w_{l}^{\prime
}}%
dw^{\prime } . \label{eq:P38}
\end{equation}
With the same decomposition as in (\ref{eq:P37}) we obtain that
(\ref{eq:P38} ) does not exceed
\begin{equation}
Ch^{2-5\delta /2}\phi _{\sqrt{T-\tau }}(y-v).  \label{eq:P39}
\end{equation}
From (\ref{eq:P37}) and (\ref{eq:P39}) we obtain that
\begin{equation}
\frac{h^{2}}{2}\left| D_{v}^{\nu }[p\otimes _{h}^{\prime
}(H_{1}+A_{0})](\tau ,T,v,y)\right| \leq Ch^{1+\gamma }\phi
_{\sqrt{T-\tau }%
}(y-v)  \label{eq:P40}
\end{equation}
for some $\gamma >0.$ It remains to estimate the last summand in
(\ref {eq:P35}). It follows from the structure of the operators
$A_{k},k=1,2,3,4$, that it is enough to estimate
\begin{equation}
\sum_{i=j^{\star
}}^{n-1}\int_{ih}^{(i+1)h}(t-ih)^{2}\int_{0}^{1}(1-\gamma
)\sum_{k=1}^{4}\int D_{v}^{\nu }\left[ p(\tau ,s,v,v+w^{\prime
})\frac{%
\partial ^{4}\widetilde{p}(s,T,v+w^{\prime },y)}{\partial w_{i}^{\prime
}\partial w_{l}^{\prime }\partial w_{p}^{\prime }\partial
w_{q}^{\prime
}}%
\right] \mid _{s=ih+\gamma (t-ih)}dw^{\prime }d\gamma dt
\label{eq:P41}
\end{equation}
for fixed $i,j,p,q$. As above, we obtain that (\ref {eq:P41}) does
not exceed
\begin{eqnarray} \nonumber &&
Ch^{2}\phi _{\sqrt{T-\tau
}}(y-v)\int_{0}^{1}z^{2}\int_{0}^{1}(1-\gamma )\sum_{i=j^{\star
}}^{n-1}h\frac{1}{[(ih-\tau )+\gamma
hz]^{3/2}}\frac{1}{%
[(n-\gamma z)h-ih]^{2}}d\gamma dz
\\ &&
=Ch^{2}\phi _{\sqrt{T-\tau
}}(y-v)\int_{0}^{1}z^{2}\int_{0}^{1}(1-\gamma )\sum_{\{i:j^{\star
}h\leq ih\leq \tau +h^{\delta }/4\}}...
\\ \nonumber && \qquad
+Ch^{2}\phi _{\sqrt{T-\tau
}}(y-v)\int_{0}^{1}z^{2}\int_{0}^{1}(1-\gamma )\sum_{\{i:\tau
+h^{\delta }/4<ih\leq T-h\}}...\\ &&=I^{\prime \prime }+II^{\prime
\prime }.  \label{eq:P42}
\end{eqnarray}
Now,
\begin{eqnarray} \nonumber  \left| I^{\prime \prime }\right| &\leq& Ch^{2}h^{-5\delta
/2}\phi
_{\sqrt{%
T-\tau }}(y-v)\int_{0}^{1}z^{2}\int_{0}^{1}(1-\gamma
)\sum_{\{i:j^{\star }h\leq ih\leq \tau +h^{\delta
}/4\}}h\frac{1}{[(ih-\tau )+\gamma hz]}d\gamma dz
\\ \nonumber &\leq& Ch^{2-\varepsilon }h^{-5\delta /2}\phi _{\sqrt{T-\tau }%
}(y-v)\int_{0}^{1}z^{2-\varepsilon }\int_{0}^{1}\frac{(1-\gamma
)}{\gamma ^{\varepsilon }}\sum_{\{i:j^{\star }h\leq ih\leq \tau
+h^{\delta
}/4\}}h\frac{%
1}{[(ih-\tau )+\gamma hz]^{1-\varepsilon }}d\gamma dz
\\ &\leq& C(\varepsilon )h^{2-\varepsilon -5\delta /2}\phi
_{\sqrt{T-\tau }}(y-v). \label{eq:P43}
\end{eqnarray}
Using inequality $(h-\gamma z)h-ih=(n-i)h-\gamma zh\geq h(1-\gamma
z)\geq h(1-\gamma )$ we obtain that
\begin{eqnarray} \nonumber
\left| II^{\prime \prime }\right| &\leq& Chh^{-5\delta /2}\phi
_{\sqrt{T-\tau }%
}(y-v)\int_{0}^{1}z^{2}\int_{0}^{1}d\gamma \sum_{\{i:\tau
+h^{\delta }/4<ih\leq T-h\}}h
\\
&\leq& Ch^{1-5\delta /2}\phi _{\sqrt{T-\tau }}(y-v).  \label{eq:P44}
\end{eqnarray}
Now from (\ref{eq:P35}), (\ref{eq:P36}), (\ref{eq:P40}),
(\ref{eq:P41}), (%
\ref{eq:P43}) and (\ref{eq:P43}) we obtain that
\begin{equation}
\left| D_{v}^{\nu }(p\otimes H-p\otimes _{h}^{\prime }H)(\tau
,T,v,y)\right| \leq Ch^{\gamma }\phi _{\sqrt{T-\tau }}(y-v)
\label{eq:P45}
\end{equation}
for some positive $\gamma .$ The last summand in the right hand
side of (\ref {eq:P35}) admits the same estimate (\ref{eq:P45})
because of the smoothing properties of the operation $\otimes
_{h}^{\prime }.$ Hence,
\begin{equation}
\left| D_{v}^{\nu }(p-p_{\Delta })(\tau ,T,v,y)\right| \leq
Ch^{\gamma }\phi _{\sqrt{T-\tau }}(y-v).  \label{eq:P45a}
\end{equation}
Making the change of variables $v=z+v^{\prime }$ into (\ref{eq:P31})
we get that the integral w.r.t. $v$ is equal to
\begin{equation}
D_{z}^{\nu }\left[ \int_{u}^{T}d\tau \int p(u,\tau ,z,z+v^{\prime
})\chi _{\nu }(\tau ,v)D_{v}^{\nu }(p-p_{\Delta })(\tau
,T,z+v^{\prime },y)dv^{\prime }\right] .  \label{eq:P46}
\end{equation}
Taking into account (\ref{eq:P45a}) and applying integration by
parts
in (%
\ref{eq:P46}) we obtain that (\ref{eq:P46}) does not exceed
\begin{equation}
Ch^{\gamma }\int_{u}^{T}\frac{d\tau }{(\tau -u)^{3/2}}\phi
_{\sqrt{T-u}%
}(y-z)\leq \frac{Ch^{\gamma }}{\sqrt{T-u}}\phi _{\sqrt{T-u}}(y-z).
\label{eq:P47}
\end{equation}
From (\ref{eq:P31}) and (\ref{eq:P47}) we obtain that
\[
\left| I_{3}\right| \leq Ch^{1+\gamma }\phi _{\sqrt{T}}(y-x).
\]
The estimate for $I_{1}$ can be proved analogously to the estimate
for $I_{3}$. Thus, we proved that
$$
hp\otimes \mathcal{F}_{1}[p\otimes \mathcal{F}_{1}[p-p_{\Delta
}]]=O(h^{1+\delta }\phi _{\sqrt{T}}(y-x)) . $$The estimate
\[
h^{1/2}p\otimes \mathcal{F}_{1}[p-p_{\Delta }]=O(h^{1+\delta }\phi
_{\sqrt{T}%
}(y-x))
\]
can be proved by using the same decomposition of $p-p_{\Delta }.$
This completes the proof of Theorem 1.

\bigskip

\bigskip

{\large \bf References.}

\begin{enumerate}

\item
 Bally V., Talay D. (1996 a). The law of the Euler scheme for
stochastic differential equations: I. Convergence rate of the
distribution function. \textit{Probability Theory and Related
Fields}, \textbf{104}, 43-60.

\item
 Bally V., Talay D. (1996 b). The law of the Euler scheme for
stochastic differential equations: II. Convergence rate of the
density. \textit{Monte Carlo Methods Appl.}, \textbf{2}, 93-128.

\item Bertail, P. and Cl\'{e}men\c{c}on (2004). Edgeworth expansions
of suitably normalized sample mean statistics for atomic Markov
chains. \textit{Probability Theory and Related Fields},
\textbf{130}, 388-414.

\item Bertail, P. and Cl\'{e}men\c{c}on (2006). Regenerative block bootstrap for Markov
chains. \textit{Bernoulli}, \textbf{12}, 689-712.

\item
 Bhattacharya R. and Rao R. (1976). \textsl{Normal approximations and
asymptotic expansions.} John Wiley \& Sons, New York.

\item Bolthausen, E. (1980). The Berry-Esseen theorem for
functionals of discrete Markov chains. \textit{Z. Wahrsch. verw.
Geb.}, \textbf{54}, 59-73.

\item Bolthausen, E. (1982). The Berry-Esseen theorem for
strongly mixing Harris recurrent Markov chains. \textit{Z. Wahrsch.
verw. Geb.}, \textbf{60}, 283-289.

\item
 Friedman A. (1964). \textsl{Partial differential equations of parabolic type.}
Prentice-Hall, Englewood Cliffs, New Jersey.

\item Fukasawa, M. (2006a). Edgeworth expansion for ergodic
diffusions. \textit{Preprint.}

\item Fukasawa, M. (2006b). Regenerative block bootstrap for ergodic
diffusions. \textit{Preprint.}

\item
 G\"{o}tze, F. (1989). Edgeworth expansions in functional limit theorems.
\textit{Ann. Probab.}, \textbf{17}, 4, 1602-1634.

\item
 G\"{o}tze, F. and Hipp, C. (1983). Asymptotic expansions for sums
 of weakly dependent random vectors.
\textit{Z. Wahrsch. verw. Geb.}, \textbf{64}, 211-239.
\item
 Guyon J. (2006). Euler scheme and tempered distributions. \textit{Stoch.
Proc. Appl. } \textbf{116}, 877-904. \item
 Jacod J. (2004). The Euler scheme for Levy driven stochastic
differential equations: limit theorems. \textit{Ann. Probab.},
\textbf{32}, 1830-1872. \item
 Jacod J., Protter P. (1998). Asymptotic error distributiions for the
Euler method for stochastic differential equations. \textit{Ann.
Probab.}, \textbf{26}, 267-307. \item
 Jacod J., Kurtz T., Meleard S., Protter P. (2005). The approximate
Euler method for Levy driven stochastic differential
 equations. \textit{Ann. de l'I.H.P.} , \textbf{41}, 523-558.
 \item Jensen, J.L. (1989). Asymptotic expansions for strongly
 mixing Harris recurrent Markov chains. \textit{Scand. J. Statist.}, \textbf{16},
 47-63.
\item
 Konakov V., Molchanov S. (1984). On the convergence of Markov chains
to diffusion processes. \textit{Teoria veroyatnostei i
matematiceskaya statistika}, \textbf{31}, 51-64 (in russian)
[English translatiion in \textit{Theory Probab. Math. Stat.} (1985),
\textbf{31}, 59-73]. \item
 Konakov V., Mammen E. (2000). Local limit theorems for transition
densities of Markov chains converging to diffusions. \textit{Probab.
Theory Relat. Fields.} \textbf{117}, 551-587. \item
 Konakov V., Mammen E. (2002). Edgeworth type expansions for Euler
schemes for stochastic differential equations. \textit{Monte Carlo
Methods Appl.}, \textbf{8}, 271-286. \item
 Konakov V., Mammen E. (2005) . Edgeworth-type expansions for
transition densities of Markov chains converging to diffusions.
\textit{Bernoulli}, \textbf{11}, 4, 591-641.
\item Kusuoka, S. and Yoshida, N. (2000).
Malliavin calculus, geometric mixing, and expansion of diffusion
functionals. \textit{Probability Theory and Related Fields},
\textbf{116}, 457-484.
\item
 Ladyzenskaya O.A., Solonnikov V.A., Ural'ceva N. (1968). \textsl{Linear and
quasi-linear equations of parabolic type.} Amer. Math. Soc.,
Providence, Phode Island.

\item Malinovskii, V.K. (1987).
Limit theorems for Harris Markov chains, 1. . \textit{Theory Probab.
Appl.}, \textbf{31}, 269-285.

\item Mykland, P.A.(1992).
Aymptotic expansions and bootstrapping distributions for dependent
variables: A martingale approach. \textit{Ann. Statist.},
\textbf{20}, 623-654.

\item
 McKean H.P., Singer I.M. (1967). Curvature and the eigenvalues of the
Laplacian. \textit{J. Diff. Geometry}, \textbf{1}, 43-69.

\item
 Protter P., Talay D. (1997). The Euler scheme for Levy driven
stochastic differential equations. \textit{Ann.Probab.},
\textbf{25}, 393-323.

\item
 Skorohod A.V. (1965). \textsl{Studies in the theory of random processes.}
Addison-Wesley. Reading, Massachussetts. [English translation of \
Skorohod A. V. (1961). \textsl{Issledovaniya po teorii sluchainykh
processov.} Kiev University Press].

\item Stroock D.W., Varadhan S.R. (1979). \textsl{Multidimensional
diffusion processes.} Springer, Berlin, Heidelberg, New York.

\item Yoshida, N. (2004).
Partial mixing and Edgeworth expansion. \textit{Probability Theory
and Related Fields}, \textbf{129}, 559-624.

\end{enumerate}

\end{document}